\documentclass[matann2]{svjour1}

\usepackage{amsmath}
\usepackage{amsfonts}
\usepackage{latexsym}
\usepackage{mathrsfs}

\usepackage{times} 
\usepackage{txfonts} 

\usepackage{ytableau} 

\usepackage{array}
\usepackage{amscd}
\usepackage{a4}
\usepackage[mathcal]{euscript}
\usepackage{amssymb}
\usepackage{epsf, epic, eepic} 
\usepackage{pstricks}
\usepackage[alphabetic]{amsrefs} 
\usepackage{theorem}

\usepackage{relsize}%

\newtheorem{definitionnumberless}{Definition.}

\spnewtheorem*{remark*}{Remark}{\it}{\rm}
\spnewtheorem*{remarks*}{Remarks}{\it}{\rm}
\spnewtheorem{vexercise}{Exercise}{\it}{\rm}
\spnewtheorem*{firstproof*}{First proof}{\it}{\rm}
\spnewtheorem*{secondproof*}{Second proof}{\it}{\rm}
\input{boxedeps} 
\SetEPSFDirectory{} 
\HideDisplacementBoxes
\SetRokickiEPSFSpecial  
\input xy 
\xyoption{all} 
\newcmykcolor{cyan}{1 0 0 0}
\DeclareMathAlphabet{\ams}{U}{msb}{m}{n}
\DeclareMathAlphabet{\goth}{U}{euf}{m}{n}
\def\so{\text{SO}}\def\sl{\text{SL}}
\def\gl{\text{GL}}
\def\m{\text{M}}\def\sp{\text{Sp}}
\def\mso{\text{MSO}}\def\msp{\text{MSp}}

\def\id{\text{id}}

\def\dom{\text{dom}\,}

\def\ov{\overline}

\def\JJ{\mathcal J}\def\HH{\mathcal H}
\def\RR{\mathcal R}\def\LL{\mathcal L}
\def\DD{\mathcal D}
\def\aa{\alpha}\def\bb{\beta}\def\ss{\sigma}

\def\R{\ams{R}}
\def\C{\ams{C}}
\def\G{\ams{G}}
\def\M{\ams{M}}
\def\0{\mathbf{0}}
\def\1{\mathbf{1}}
\def\quo{/\kern -.45em\sim}
\newpsobject{showgrid}{psgrid}{subgriddiv=1,griddots=10,gridlabels=6pt,gridcolor=red}
\def\Langle{\langle\kern -2pt\langle}
\def\Rangle{\rangle\kern -1.9pt\rangle}
\def\irr{\text{Irr}}
\def\sgl{S\kern-1pt(G,L)}
\def\swl{S\kern-1pt(W,L)}
\DeclareMathOperator{\im}{im} 

\newcommand{\ra}{\rightarrow}

\newcommand{\Symn}{S_{\kern-.3mm n}}
\newcommand{\SymX}{S_{\kern-.3mm X}}
\newcommand{\Symthree}{S_{\kern-.3mm 3}}
\newcommand{\Symfour}{S_{\kern-.3mm 4}}
\newcommand{\Symtwo}{S_{\kern-.3mm 2}}
\newcommand{\Symone}{S_{\kern-.3mm 1}}
\newcommand{\Symm}{S_{\kern-.3mm m}}
\newcommand{\Symr}{S_{\kern-.5mm r}}

\newcommand{\cork}[1]{|\kern0.75pt{#1}\kern1pt|}

\newcommand{\hs}{H\kern-0.5pt S}
\newcommand{\htt}{H\kern-0.5pt T}
\newcommand{\hc}{H\kern-0.5pt C}

\title{The sympathetic sceptic's guide to semigroup representations}

\author{Brent Everitt\thanks{Based on lectures given at York in the
    Spring of 2016. I am grateful to Majed Albaity, Sanaa Bajri, 
    Michael Bate, Steve Donkin, James East, John Fountain and
    particularly to Vicky Gould for her constant encouragement.}}

\institute{
{\sc Brent Everitt:} Department of Mathematics, University of York, York
YO10 5DD, United Kingdom. \email{brent.everitt@york.ac.uk}. 
}

\titlerunning{The sympathetic sceptics guide to semigroup representations}
\authorrunning{Brent Everitt}

\begin{document}

\maketitle

\begin{abstract}
This is an elementary, examples driven, introduction to the representation
theory of finite semigroups.
We 
illustrate the Clifford-Munn correspondence 
between the representations of a semigroup and the representations of
its maximal subgroups. The emphasis throughout is on naturally
occurring examples. 
\end{abstract}


\section*{Introduction} 

This is an elementary introduction to the representation
theory of finite semigroups.
As the title suggests, it is not necessarily intended for semigroup
theorists.

We start with a quick primer on the semigroups that will interest us
-- the inverse and regular monoids -- and spend a certain amount of
energy selling these objects to the general mathematical
public. 
Section \ref{section:representations} introduces
from scratch 
linear actions of semigroups on vector spaces, where the emphasis is
on those aspects of the theory that are in common with group
representations. By this point the symmetric group $\Symn$ will have
appeared a number of times,
so we divert to describe
its ``atomic'' representations. There are two fundamental
constructions, reduction and induction, that connect group theory and
semigroup theory, at least when it comes to representations. These are
described in Sections
\ref{section:reduction}-\ref{section:induction}. 
Section \ref{section:clifford:munn} contains, what is, from the point
of view of these notes, the 
fundamental theorem
of semigroup representation theory: the Clifford-Munn
correspondence.
It gives a mechanism for producing the atomic
representations of semigroups using only knowledge from group
theory. The last section is essentially a gratuitous excuse to draw
pictures of our favourite polytope, the permutohedron, dressed up as a
worked example of the representations of an interesting Renner monoid.

Throughout, three running examples, $\Symn$ (a group), $I_n$ (an
inverse monoid) and $T_n$ (a regular monoid) are used as
illustration. 
By the end of Section \ref{section:induction} the emphasis 
will have completely moved to inverse monoids. We also start with
actions on vector spaces over an arbitrary field $k$, but in later
sections we
retreat to the relative safety of representations over $\C$. We borrow
heavily from a number of sources
-- full attributions are given in the Notes and References section at
the end. 


\section{Semigroups}
\label{section:semigroup_basics}

A semigroup is a set equipped with an associative multiplication. This
leaves us with quite a bit of scope! In this section we feel our way towards a
manageable class of semigroups to study. Our guiding principle will
be the role of inverses in semigroup theory.

We start with three finite examples that are the most typical of their
type. Throughout, we write $[n]$ for the set $\{1,2,\ldots,n\}$.
\begin{description}
\item[--]  \emph{The symmetric group $\Symn$}: consisting of all
  bijections $g:[n]\rightarrow [n]$ with multiplication the usual
  composition of maps. 
\item[--] \emph{The symmetric inverse monoid $I_n$}: consisting of all
  \emph{partial\/} bijections $s:[n]\rightarrow [n]$, i.e. bijections $s:X\rightarrow Y$
  where $X,Y\subseteq [n]$. The multiplication is composition of partial
  maps as shown in Figure \ref{fig:partialmapcomposition}.
\begin{figure}[h]
  \centering
\begin{pspicture}(0,0)(12,6.5)
\rput(6,3.25){\BoxedEPSF{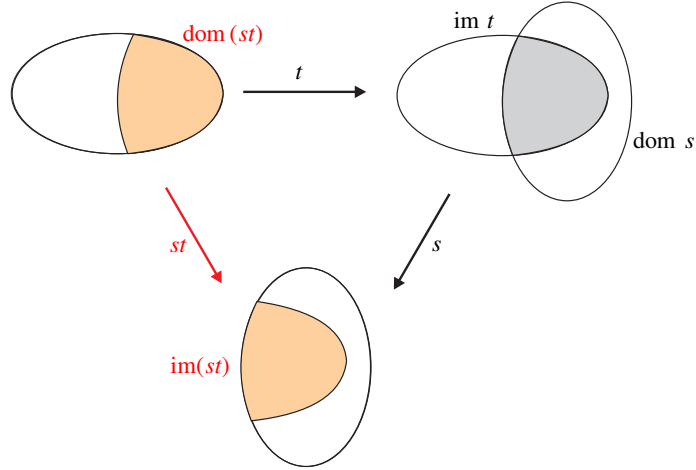 scaled 450}}
\rput(4.1,3.1){${\red st}$}\rput(4.7,5.9){${\red \dom(st)}$}
\rput(4.4,1.5){${\red \im(st)}$}
\rput(8,6.1){$\im\,t$}\rput(10.5,4.5){$\dom\,s$}
\rput(5.7,5.4){$t$}\rput(7.5,3.1){$s$}
\end{pspicture}
\caption{composition of partial maps of $[n]$.}
  \label{fig:partialmapcomposition}
\end{figure}
All our functions, actions, etc, will be on the left,
so the partial map $st$ has domain the $t$-preimage of $\im
t\cap\dom s$ and image the $s$-image of $\im t\cap \dom s$, and is the usual
composition of $t$ followed by $s$ between these two sets. 
If
$\im t\cap \dom s$ is empty, then $st$ is the unique bijection
$\varnothing\rightarrow\varnothing$, which we will call the zero map
$0$. 
\item[--] \emph{The full transformation monoid  $T_n$}: consisting of all mappings
  $s:[n]\rightarrow [n]$ with multiplication the usual composition of maps. 
\end{description}

\paragraph{Inverses in semigroup theory.} Naively, a semigroup is a
group, except without inverses. But to completely rule out inverses in a
semigroup is unnecessarily defeatist. 
The elements of $\Symn$ are ``global'' symmetries of the set $[n]$ -- with
global inverses -- while the elements of $I_n$ are ``local'' symmetries
of $[n]$, with local inverses to match. 
Our three running examples motivate
three ways in which inverses arise:

\smallskip

(i). The symmetric group $\Symn$ is a group, obviously, so for every
$g\in\Symn$ there is a unique 
  $h\in\Symn$ with $gh=\id=hg$. Write $h=g^{-1}$ as usual.

\smallskip

(ii). If $s:X\rightarrow Y$ is an element of  $I_n$ then there is a unique 
  $X\leftarrow Y:t$, that is the inverse of $s$, but only defined on the set
  $Y$. Indeed, $st=\id_Y$ and $ts=\id_X$, where $\id_X:X\ra X$ and
  $\id_Y:Y\ra Y$ are partial identities, and in particular,
  idempotents (i.e: $\id_X\id_X=\id_X$ and $\id_Y\id_Y=\id_Y$).
 
  As a working definition of the local inverse of $s$, we could
  take it to be an element $t$ such that $st$ and $ts$ are idempotents,
  but not necessarily the global idempotent $\id$. It
  turns out that this isn't quite satisfactory, as any map defined
  on some subset of the image of $s$, and equal to the inverse of $s$ on this
  subset, also has this property.

  Instead, we have $s\,\id_X =s=\id_Y s$. Together with $st=\id_Y$ and
  $ts=\id_X$ we get that $t$ satisfies $sts=s$; similarly
  $tst=t$. 

A semigroup with the property that for every $s$
  there is a \emph{unique\/} $s^*$ satisfying 
  \begin{equation}\label{eq:inverse_semigroup}
    ss^*s=s\text{ and }s^*ss^*=s^*
  \end{equation}
  is called an inverse semigroup; a semigroup with an
  identity $\id$ is a monoid, and an inverse semigroup with an identity
  $\id$ is an \emph{inverse monoid\/}. $I_n$ is thus the most inverse
  monoid-like of the inverse monoids. We will sometimes call $s^*$ an
  inverse ``in the sense of semigroup theory'' and reserve the notation
  $s^{-1}$ for inverses in a group.

\smallskip

\begin{figure}
  \centering
\begin{pspicture}(0,0)(14,3.5)
\rput(0,0.35){
\rput(3.5,1.5){\BoxedEPSF{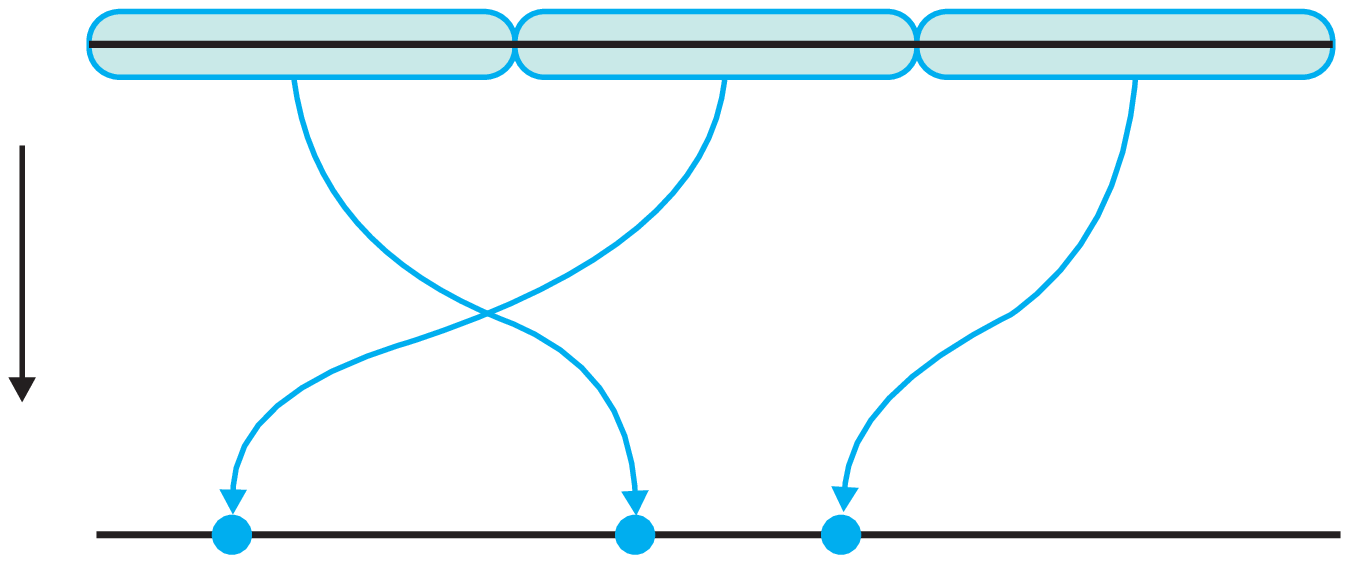 scaled 450}}
\rput(10.5,1.5){\BoxedEPSF{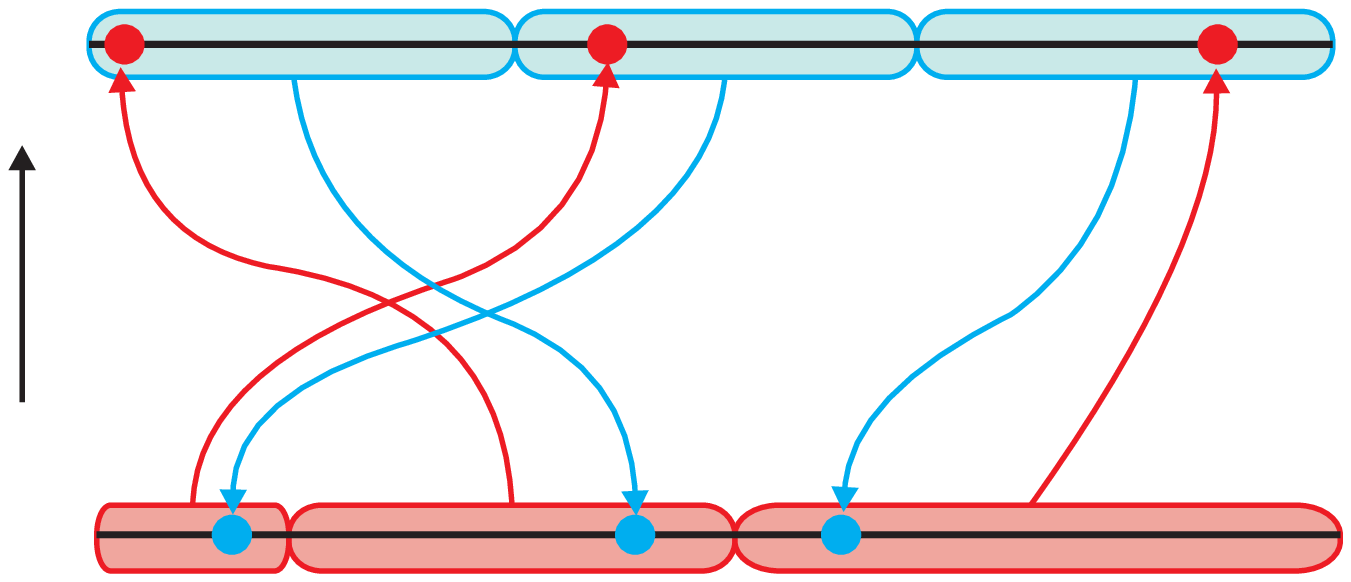 scaled 450}}
\rput(0.3,1.5){$s$}
\rput(0.85,3){$1$}\rput(6.5,3){$n$}
\rput(0.9,0.05){$1$}\rput(6.55,0.05){$n$}
\rput(3.5,3){{\cyan fibers/kernel of $s$}}
\rput(7.3,1.5){$t$}
\rput(7.8,2.9){$1$}\rput(13.4,2.9){$n$}
\rput(7.8,0.05){$1$}\rput(13.4,0.05){$n$}
\rput(10.5,0){{\red fibers/kernel of $t$}}
\rput(10.2,2.9){$x$}\rput(8.5,0.05){$sx$}
}
\end{pspicture}
\caption{inverses in $T_n$.}
  \label{fig:inverses_in_Tn}
\end{figure}

(iii). Definition (\ref{eq:inverse_semigroup}) of inverses opens up new
possibilities. The element $s\in T_n$ shown on the left of Figure
\ref{fig:inverses_in_Tn} has \emph{kernel\/} the 
partition of $[n]$ whose blocks are the \emph{fibers\/} of $s$: the
$s$-preimages of a point in the image of $s$.
Now construct an element $t\in T_n$ in the following way:
\begin{description}
\item[--] partition the image copy of $[n]$ so that each block of the
  partition 
  contains exactly one element of the image of $s$; this partition 
  will be the kernel of $t$.
\item[--] For each block in this new kernel choose an $x$ in the domain
  copy of $[n]$ such that the
  block contains the point $sx$; then
  define $t$ so that it maps this block to $x$; see the right of
  Figure \ref{fig:inverses_in_Tn}.
\end{description}
The $t$ just constructed satisfies $sts=s$ and $tst=t$; conversely,
any $t$ satisfying these relations
must come about in this way.
But this $t$ is clearly not
unique -- there is choice in the partition of the image $[n]$ and for each block
of this partition, choice in the
$x$ so that the block is labelled by $sx$. 

A monoid with the property that for every $s$
there is \emph{some\/} $t$ satisfying 
$sts=s$ and $tst=t$
is called a \emph{regular monoid\/}.

\emph{From now on:} $S$ will be a finite regular monoid. 


\paragraph{The structure of semigroups: Green's relations.} These allow
us to draw strategic pictures of semigroups. Define an equivalence relation
$\LL$ on $S$ by $s\LL t$ when $S\kern-0.5mm s=S\kern-0.35mm
t$, where $S\kern-0.5mm s=\{rs\,:\,r\in
S\}$ is a left ideal (hence the ``$\LL$''). Dually, define $s\RR t$
when $sS=tS$. 

In $\Symn$, and indeed any group, these relations are trivial: any two
elements are $\LL$ and $\RR$-related. In $I_n$ and $T_n$ they take a particularly
simple form:
\begin{description}
\item[--] $s\LL t$ when the fibers of $s$ are equal to the
  fibers of $t$ (or $s$ and $t$ have the same kernel). In $I_n$ this
  is equivalent to $\dom s=\dom t$.  
\item[--] $s\RR t$ when $\im s=\im t$.
\end{description}

\begin{figure}[t]
  \centering
\begin{pspicture}(0,0)(14,7)
\rput(0,1){
\rput(3.5,2.5){\BoxedEPSF{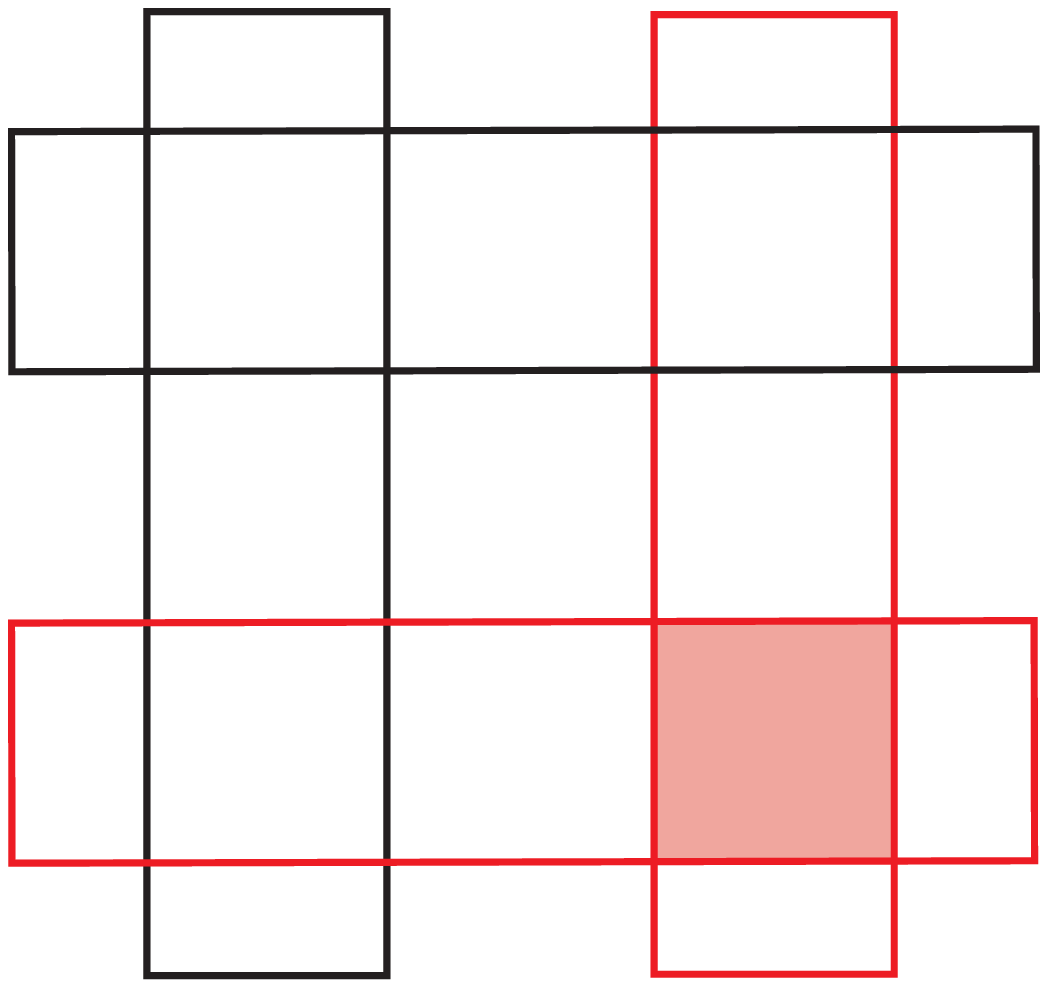 scaled 450}}
\rput(2.3,3.65){$X\stackrel{s}{\ra}Y$}\rput(4.6,3.65){$W\stackrel{t}{\ra}Y$}
\rput(2.3,1.45){$X\stackrel{r}{\ra}Z$}\rput(4.6,1.45){$W\stackrel{q}{\ra}Z$}
\rput(0.55,3.65){$\im=Y$}\rput(0.55,1.45){$\im=Z$}
\rput(2.3,4.95){$\dom=X$}\rput(4.6,4.95){$\dom=W$}
\rput(4.6,0.05){${\red \LL_t}$}\rput(6.1,1.45){${\red \RR_r}$}
}
\rput(10,4){\BoxedEPSF{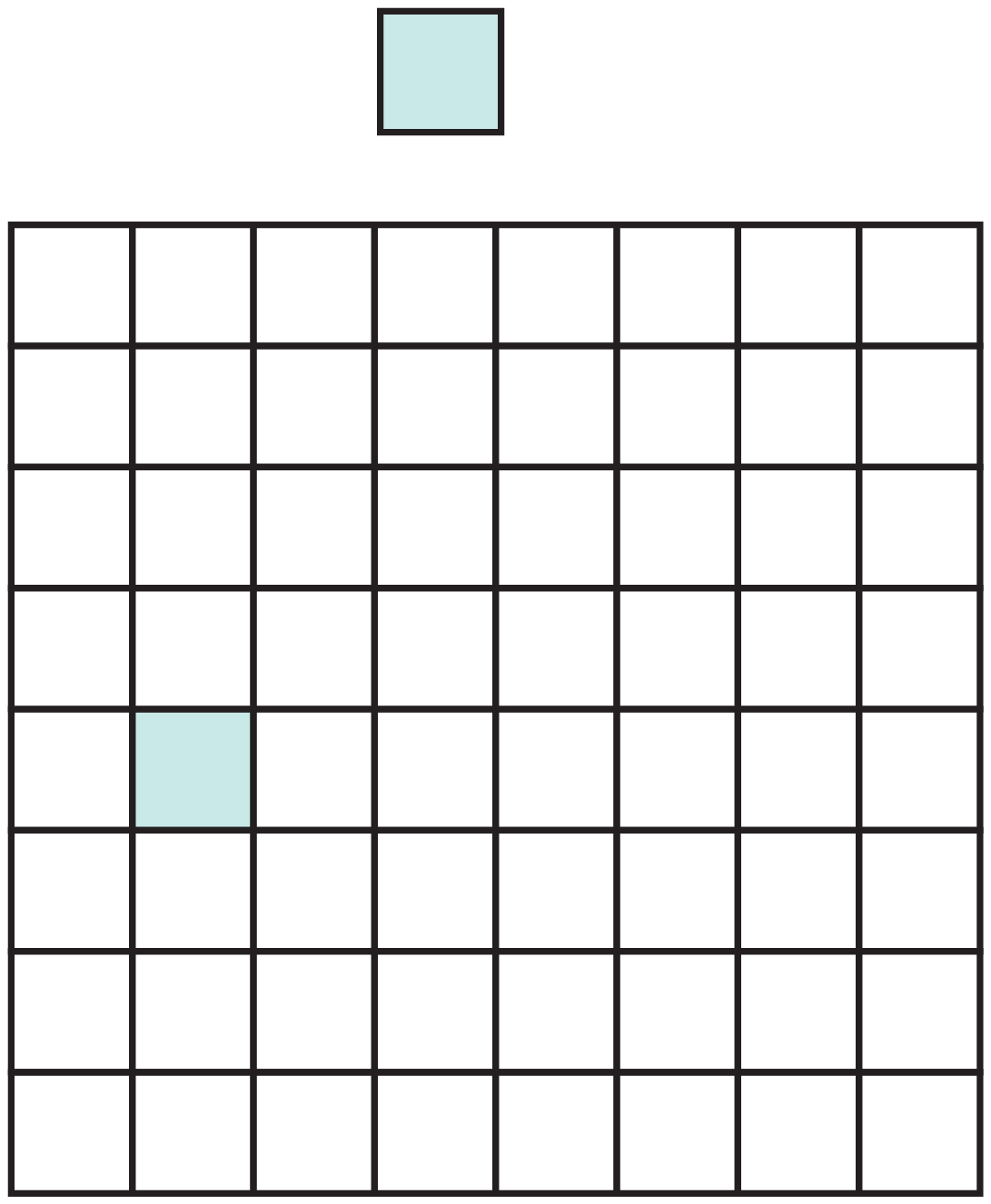 scaled 450}}
\rput(8,6.2){$\langle\LL,\RR\rangle$-class eggbox:}
\rput(4.1,-0.1){\rput{90}(5.9,.8){$\left\{\begin{array}{c}
\vrule width 0 mm height 50 mm depth 0 pt\end{array}\right.$}}
\rput(12.8,3.4){$\left.\begin{array}{c}
\vrule width 0 mm height 50 mm depth 0 pt\end{array}\right\}$}
\rput(10,0.4){$\LL$-classes $=$ domain sets of size $m$}
\rput(0,0.4){
\rput(13.7,3.4){$\RR$-classes}
\rput(13.7,3.1){$=$}
\rput(13.7,2.8){image sets}
\rput(13.7,2.5){of size $m$}
}
\rput(8.5,0.85){$X$}\rput(12.6,3.2){$Y$}
\rput(11.6,6.7){$=$ all bijections $X\stackrel{s}{\ra}Y$}
\end{pspicture}
\caption{In the symmetric inverse monoid $I_n$, the
$\LL$-class $\LL_t$ of any element $t$ that is $\RR$-related to $s$ intersects the
$\RR$-class $\RR_r$ of any element $r$ that is  $\LL$-related to $s$
\emph{(left)\/} and an eggbox grid of a $\langle\LL,\RR\rangle$-class
\emph{(right)} partitioned into mutually intersecting $\LL$ and $\RR$-classes.}
  \label{fig:LR_commute}
\end{figure}

If we consider the equivalence relation $\langle\LL,\RR\rangle$ 
generated by $\LL$ and $\RR$, then something very nice happens. 
The
$\LL$-class of any element $t$ that is $\RR$-related to $s$ intersects the
$\RR$-class of any element $r$ that is  $\LL$-related to $s$. It is particularly
easy to see for $I_n$ as on the left of  Figure
\ref{fig:LR_commute}. The
$\langle\LL,\RR\rangle$-classes are thus partitioned into $\LL$-classes and
partitioned into $\RR$-classes, with any $\LL$-class intersecting any
$\RR$-class and vice-versa. Semigroup theorists call this grid an
``eggbox'' -- see the right of Figure \ref{fig:LR_commute}.

Moreover, pursuing the ideal theme, define a relation $\JJ$ on $S$
by $s\JJ t$ when $S\kern-0.5mm sS=S\kern-0.35mm tS$. Again,
this has a simple form in 
$I_n$ and $T_n$, with $s\JJ t$ when $\im s$ and $\im t$ are sets of
the same size. But the $\langle\LL,\RR\rangle$-class of $I_n$ in
Figure \ref{fig:LR_commute} consists precisely of those partial
bijections whose image has the fixed size $|\im s\,|=|Y|$. Thus
$\JJ=\langle\LL,\RR\rangle$ in $I_n$, and in general for any finite $S$.

Our final relation is $\HH=\LL\cap\RR$, so that $s\HH t$ when they
are both $\LL$ and $\RR$-related. In $I_n$ and $T_n$ 
a pair $s\HH t$ means that $s$ and $t$ have the same fibers (or
domains in $I_n$) and the same 
images. The $\HH$-classes are thus the small boxes in the eggbox grid
with one marked on the right of Figure \ref{fig:LR_commute}.
Write $\LL_s,\RR_s,\JJ_s$ and $\HH_s$ for the equivalence class of
$s\in S$ under these relations. 

The $\JJ$-classes 
are not just floating around in the ether in a disembodied fashion. They
can be compared to each other; in other words, they
form a poset. This is what we mean by ``strategic picture''.

Again we can see this quite naturally by looking at
$I_n$ and $T_n$, where the $\JJ$-classes are parametrised by the
possible sizes of the images: by $\{0,\ldots,n\}$ in 
$I_n$ and by $\{1,\ldots,n\}$ in $T_n$.
Indeed, $S\kern-0.5mm sS$ consists of the maps having image size $\leq
|\im s\,|$, so that $S\kern-0.5mm sS\subseteq S\kern-0.35mm tS$ exactly
when $|\im s\,|\leq|\im t\,|$. We will write $\JJ_m$ for the
$\JJ$-class consisting of those
maps with image size $m$. 

In general, define a partial order on
the 
$\JJ$-classes of a semigroup $S$ by $\JJ_s\leq\JJ_t$ whenever
$S\kern-0.5mm sS\subseteq S\kern-0.35mm tS$. 

\begin{figure}[t]
  \centering
\begin{pspicture}(0,0)(14,9)
\rput(0,-0.5){
\rput(3.5,5){\BoxedEPSF{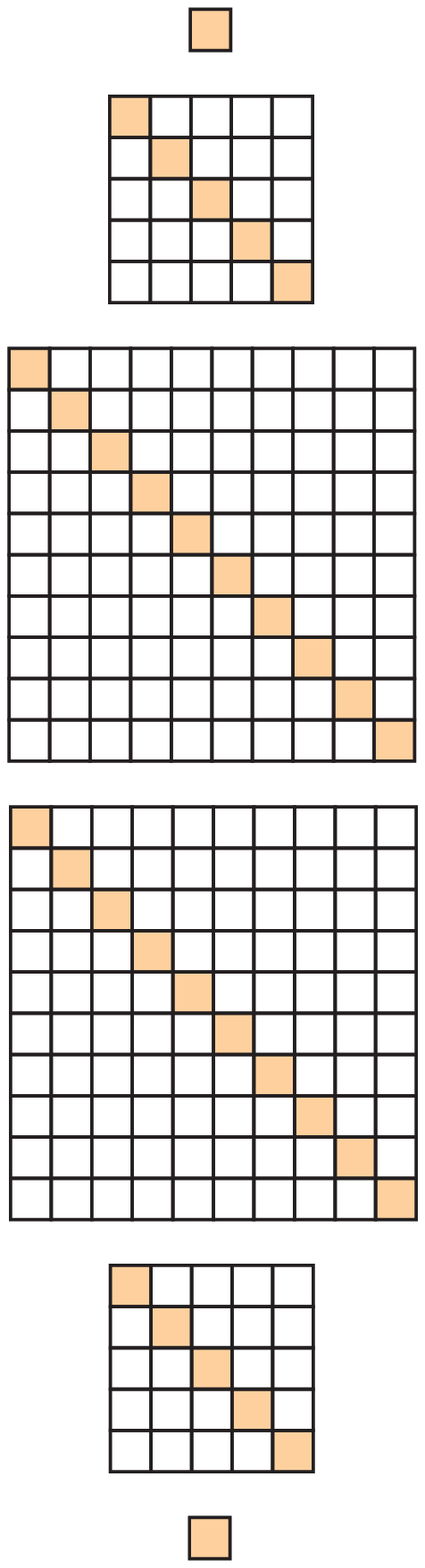 scaled 450}}
\rput(10.5,5){\BoxedEPSF{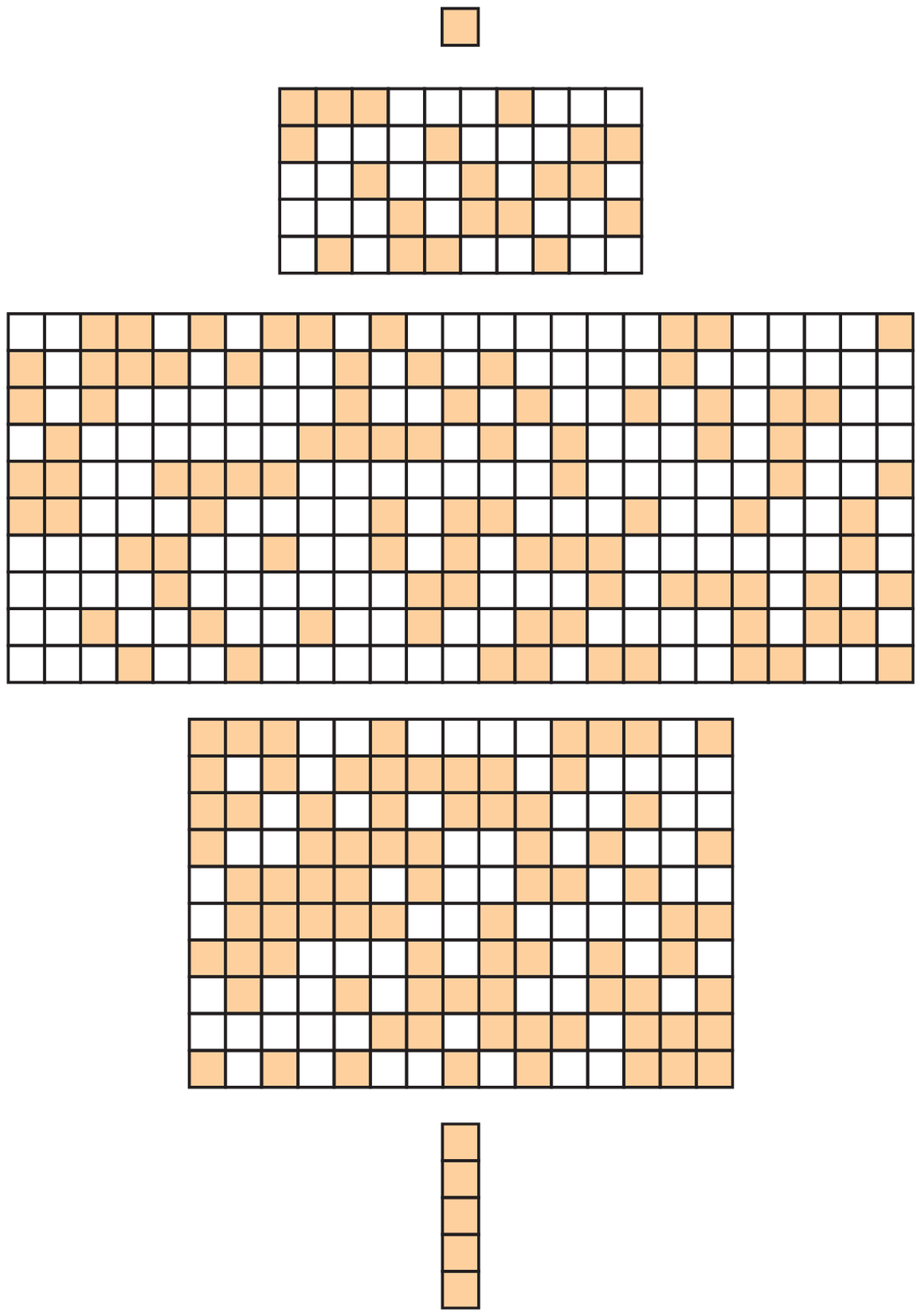 scaled 450}}
\rput(4,8.9){$\JJ_5$}\rput(4.5,8){$\JJ_4$}\rput(5,6.1){$\JJ_3$}
\rput(5,3.9){$\JJ_2$}\rput(4.5,2){$\JJ_1$}\rput(4,1.1){$\JJ_0$}
\rput(10,8.6){$\JJ_5$}\rput(9,7.7){$\JJ_4$}\rput(7.5,6){$\JJ_3$}
\rput(8.5,3.5){$\JJ_2$}\rput(10,1.8){$\JJ_1$}
}
\end{pspicture}
\caption{The $\HH$-classes containing idempotents in $I_5$
  \emph{(left)} and $T_5$ \emph{(right)}.}
  \label{fig:Tn_Hclasses_idempotents}
\end{figure}

\paragraph{Idempotents.} An idempotent is an element $e$
with the property that $e^2=e$. In a group there is precisely one: the identity
$\id$. But in $I_n$ there are others, and in $T_n$ even more again. 

The idempotents in $I_n$ are the maps $\id_X:X\rightarrow X$
that are the identity on some $X\subseteq [n]$; they are the \emph{partial
  identities\/}. In the eggbox on the right of Figure \ref{fig:LR_commute},
 $\id_X$ lives on the
diagonal in the row and
column labelled by $X$.
Moreover, $\id_X$ is the only
idempotent in its 
row and column. In fact, this is true for any inverse semigroup: each
$\RR$-class and each $\LL$-class contains a unique idempotent.

\begin{figure}[b]
  \centering
\begin{pspicture}(0,0)(14,3)
\rput(7,1.5){\BoxedEPSF{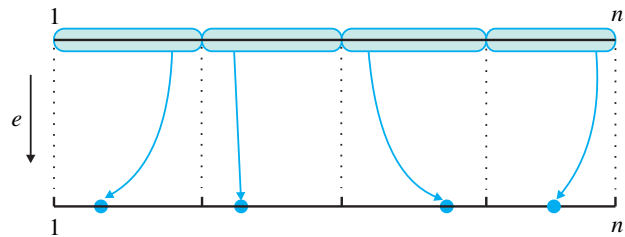 scaled 450}}
\rput(3.5,2.9){$1$}\rput(10.9,2.9){$n$}
\rput(3.5,0.1){$1$}\rput(10.9,0.1){$n$}
\rput(3,1.5){$e$}
\end{pspicture}
\caption{An idempotent in $T_n$.}
  \label{fig:Tn_idempotents}
\end{figure}

To find idempotents in $T_n$, fix any partition of $[n]$; this will be
the fibers/kernel of $e$. In each fiber fix a point, and then define $e$ 
to map each fiber to the point chosen in it -- see Figure
\ref{fig:Tn_idempotents}. (A slicker way to say it is that $e$
restricts to the identity on its image.)
If the fibers -- and hence the $\LL$-class -- are fixed, there is still
wiggle-room in the 
choice of point in each one. So a given $\LL$-class may contain several
idempotents. Dually, fixing some image points (and hence the
$\RR$-class) there are many partitions of $[n]$ with a unique image
point in each block of the partition, and so several idempotents in a
given $\RR$-class.  This behaviour is typical of regular, non-inverse
semigroups. Figure
\ref{fig:Tn_Hclasses_idempotents} compares the
$\HH$-classes containing idempotents in $I_5$ and $T_5$. 

In any case, in both $I_n$ and $T_n$ -- and in a regular monoid in
general -- each $\HH$-class contains at most one idempotent. 

\begin{figure}
  \centering
\begin{pspicture}(0,0)(14,7.5)
\rput(7,3.75){\BoxedEPSF{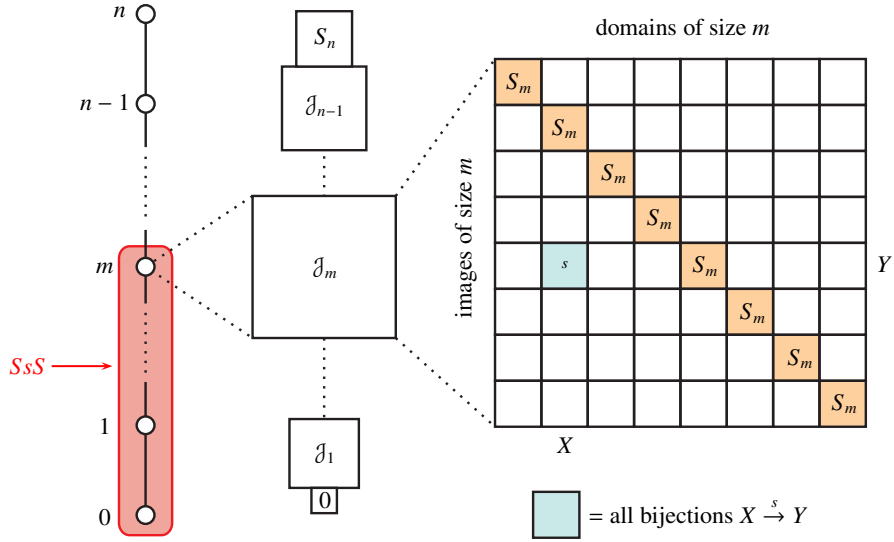 scaled 450}}
\rput(7.95,3.85){${\scriptstyle s}$}
\rput(-7.5,-0.7){
\rput(9.6,7.9){$n$}\rput(9.4,6.7){$n-1$}
\rput(9.4,4.5){$m$}\rput(9.4,2.4){$1$}\rput(9.4,1.2){$0$}
}
\rput(0.9,2.5){$\red S\kern-0.5mmsS$}
\psline[linewidth=0.75pt,linecolor=red]{->}(1.2,2.5)(2,2.5)
\rput(4.8,3.8){$\JJ_m$}\rput(4.8,1.35){$\JJ_1$}
\rput(4.8,0.725){$0$}\rput(4.8,5.95){$\JJ_{n-1}$}
\rput(4.8,6.85){$\Symn$}
\rput(9.7,0.55){$=$ all bijections $X\stackrel{s}{\ra}Y$}
\rput(7.95,1.45){$X$}\rput(12.15,3.85){$Y$}
\rput(7.35,6.25){$\Symm$}\rput(7.95,5.65){$\Symm$}
\rput(8.6,5.05){$\Symm$}\rput(9.2,4.45){$\Symm$}
\rput(9.8,3.825){$\Symm$}\rput(10.425,3.2){$\Symm$}
\rput(11.05,2.6){$\Symm$}\rput(11.63,2){$\Symm$}
\rput(9.5,7){domains of size $m$}
\rput(-2.85,-2.8){\rput{90}(9.5,7){images of size $m$}}
\end{pspicture}
\caption{Strategic picture of $I_n$. The $\JJ$-class poset
  \emph{(left)\/}, the stacked eggboxes \emph{(middle)\/} and the
  eggbox picture of the $\JJ$-class $\JJ_m$ \emph{(right)\/}; there
  are $\binom{n}{m}$ rows and columns with the maximal subgroups
  $\cong\Symm$ down the diagonal.}
  \label{fig:In_stategic_picture}
\end{figure}

\paragraph{Subgroups.} In any monoid the \emph{units\/} are the elements
that have inverses in the sense of group theory, and these form a
subgroup. In our three examples $\Symn,I_n$ and $T_n$, these are the
bijections $[n]\rightarrow [n]$, so the group of units is $\Symn$
with identity $\id:[n]\rightarrow [n]$. In $\Symn$ this is the whole story,
but in $I_n$ there are other subgroups, \emph{disjoint from the units\/}, and in $T_n$ even
more again. 

For $X\subseteq [n]$ fixed, the bijections $X\rightarrow X$ form a subgroup of
$I_n$ isomorphic to $\Symm$, where $m=|X|$. This subgroup is precisely the
diagonal $\HH$-class containing the idempotent $\id_X:X\ra X$. 

In general, if
$\HH_e$ is the $\HH$-class containing the idempotent $e$, then this is
a subgroup of $S$ with identity $e$. Moreover, any subgroup of $S$ is a
subgroup of an $\HH_e$ for some $e$, hence these are \emph{maximal
  subgroups\/} of $S$. We write $G_e$ for $\HH_e$ from now on, to
stress its group structure. 

\begin{figure}[h]
  \centering
\begin{pspicture}(0,0)(14,3)
\rput(3.5,1.5){\BoxedEPSF{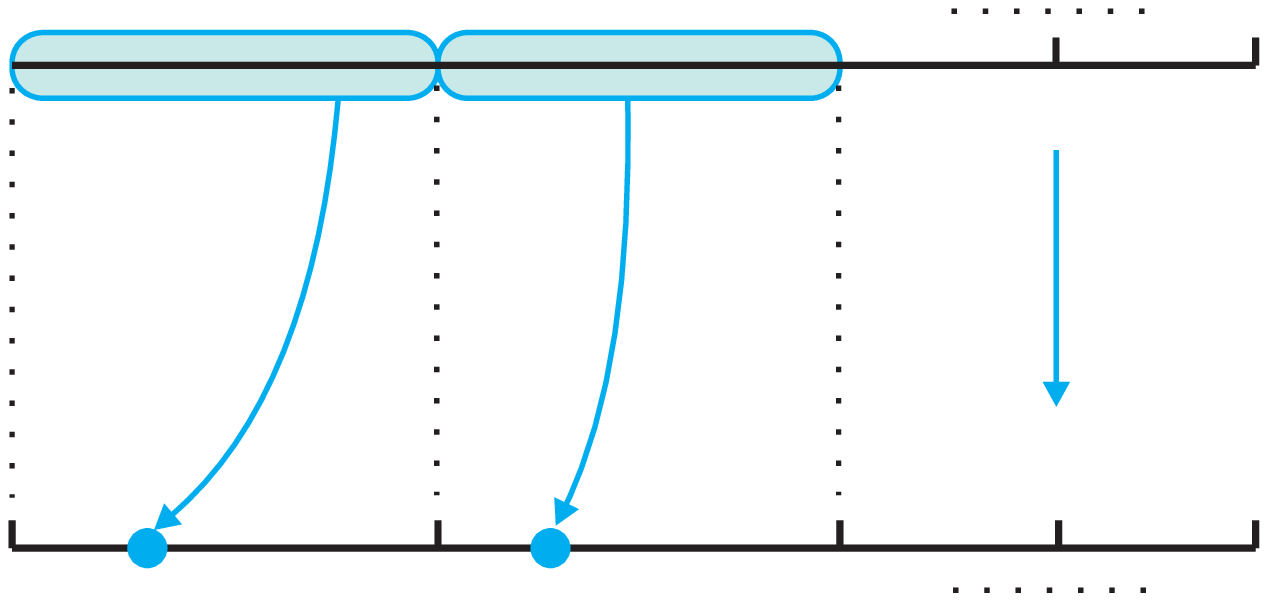 scaled 450}}
\rput(1.6,2.95){$X_1$}\rput(3.5,2.95){$X_2$}
\rput(1.25,0.1){$y_1$}\rput(3.1,0.1){$y_2$}
\rput(5.6,1.5){$e$}
\rput(10.5,1.5){\BoxedEPSF{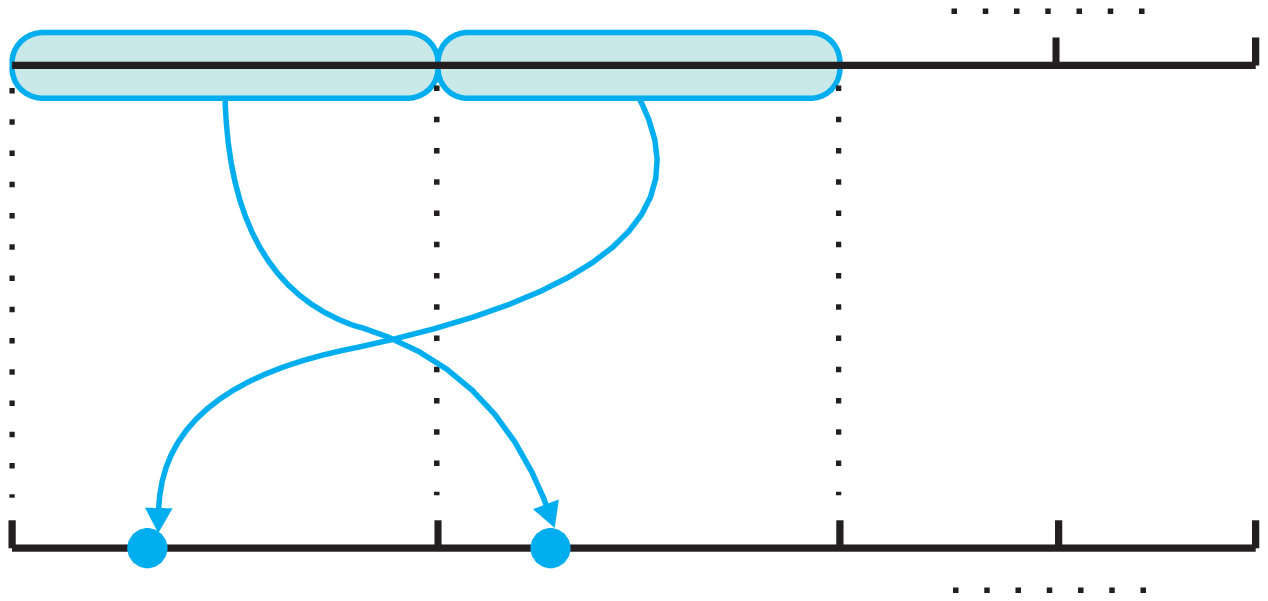 scaled 450}}
\rput(7,0){
\rput(1.6,2.95){$X_1$}\rput(3.5,2.95){$X_2$}
\rput(1.25,0.1){$y_1$}\rput(3.1,0.1){$y_2$}
}
\end{pspicture}
\caption{Maximal subgroup of $T_n$ with identity the idempotent $e$
  \emph{(left)\/} and a typical element \emph{(right)\/}.}
  \label{fig:Tn_maximal_subgroup}
\end{figure}

$T_n$ has many more $\HH$-classes containing idempotents, hence many
more maximal subgroups. The $\HH$-class containing the idempotent $e$
on the left of
Figure \ref{fig:Tn_maximal_subgroup} consists of the maps with fibers
$X_1,\ldots,X_m$ and image points $y_1,\ldots,y_m$, 
and where these maps give a bijection
$\{X_1,\ldots,X_m\}\ra\{y_1,\ldots,y_m\}$. The maximal subgroups of
$T_n$ are thus symmetric groups $\Symm$ as well, but in a slightly
different way to $I_n$. 
Figure \ref{fig:Tn_Hclasses_idempotents} therefore also shows these
$\Symm$ subgroups (shaded), for $0\leq m\leq 5$, in $I_5$ and $T_5$.

\begin{figure}[t]
  \centering
\begin{pspicture}(0,0)(14,6.5)
\rput(0,0.25){
\rput(7,3){\BoxedEPSF{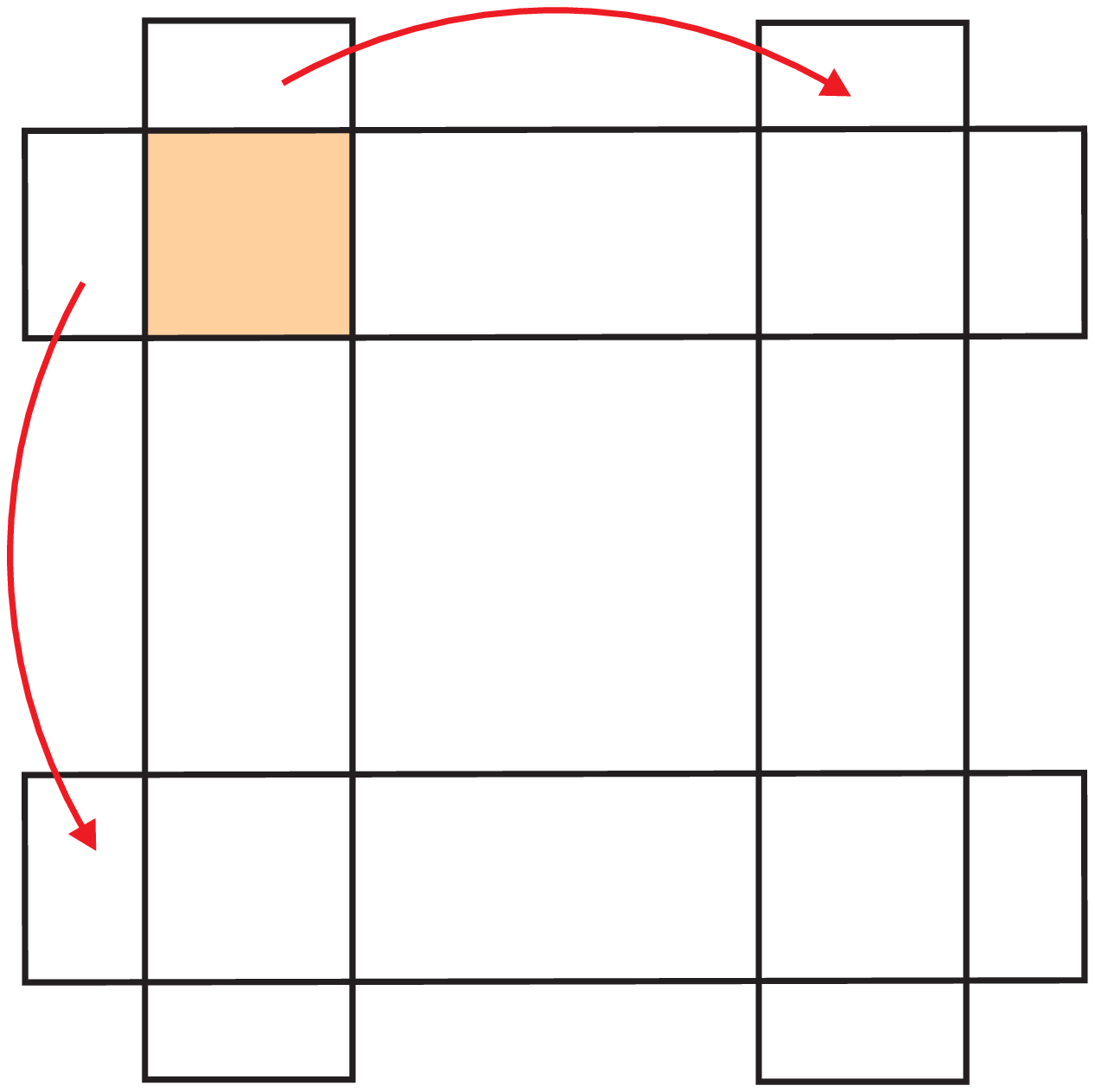 scaled 450}}
\rput(0,-0.55){
\rput(5.4,5.3){$X\stackrel{e}{\ra}X$}
\rput(5.4,1.9){$X\stackrel{s}{\ra}Y$}
\rput(8.65,5.3){$Z\stackrel{t}{\ra}X$}
\rput*(4.1,3.5){${\red s(-)}$}
\rput*(7,6.4){${\red (-)t}$}
\rput(5.4,4.4){$G_e$}
\rput(3.95,5.3){$\RR_e$}
\rput(5.4,6.6){$\LL_e$}
\rput(5.4,2.6){$\HH_s$}
\rput(8.6,4.4){$\HH_t$}
\rput(7,6){{\red bijection}}
\rput(3.1,3.5){{\red bijection}}
}
}
\end{pspicture}
\caption{Green's lemma.}
  \label{fig:Green_lemma}
\end{figure}

\paragraph{The strategic picture for $I_n$:} 
is given in Figure
\ref{fig:In_stategic_picture}. 
The $\JJ$-class poset is on the left -- the image sizes $\{0,\ldots,n\}$
with their usual total order -- and the stacked eggboxes are in the
middle. The maximal $\JJ$-class consists of all the bijections with
image size $n$, so is the symmetric group $\Symn$, and the minimal $\JJ$-class
has single element the zero map $0:\varnothing\rightarrow\varnothing$.  
The class $\JJ_m$ has rows and columns indexed by the $\binom{n}{m}$
subsets of size $m$, with the blue box in Figure \ref{fig:In_stategic_picture}
containing the bijections
$s:X\rightarrow Y$. The idempotents $\id_X:X\ra X$ lie down the
diagonal, with the maximal subgroups consisting of all the bijections
$X\ra X$ for fixed $|X|=m$, and thus $\cong\Symm$.

\paragraph{Green's lemma.} The maximal subgroups can be used to parametrise the
$\LL$ and $\RR$ classes containing them -- indeed, the
$\HH$-classes in $\RR_e$ are like right cosets of the subgroup $G_e$
and the $\HH$-classes in $\LL_e$ are like left cosets. 

It is easy to
see in $I_n$: let $e$ be the idempotent 
$\id_X:X\ra X$, contained in the maximal subgroup $G_e$ of all
bijections $X\ra X$ (see Figure \ref{fig:Green_lemma}). An $s\in\LL_e$
is a bijection $s:X\ra Y$. For any $g\in G_e$, the composition
$$
X\stackrel{g}{\ra}X\stackrel{s}{\ra}Y
$$
is a bijection $X\ra Y$, and all such bijections arise in this way via
some $g$. Put another way, left multiplication by $s$ is a bijection
$s(-):G_e\ra\HH_s$. 
Thus:
\begin{equation}
  \label{eq:1}
\text{every element of the $\HH$-class $\HH_s$ can be uniquely
expressed as $sg$ for some $g\in G_e$} 
\end{equation}
(so $\HH_s=sG_e$ is the 
left coset in $\LL_e$ of the maximal subgroup $G_e$).
Dually, if $t\in\RR_e$ is
some bijection $t:Z\ra X$ then every element of $\HH_t$ has a unique
expression as $gt$ for some $g\in G_e$, and so $\HH_t$ is the right
coset $G_et$ of $G_e$ in $\RR_e$. These observations are called 
Green's lemma.

\begin{figure}
  \centering
\begin{pspicture}(0,0)(14,7)
\rput(0,0){
\rput(7,3.5){\BoxedEPSF{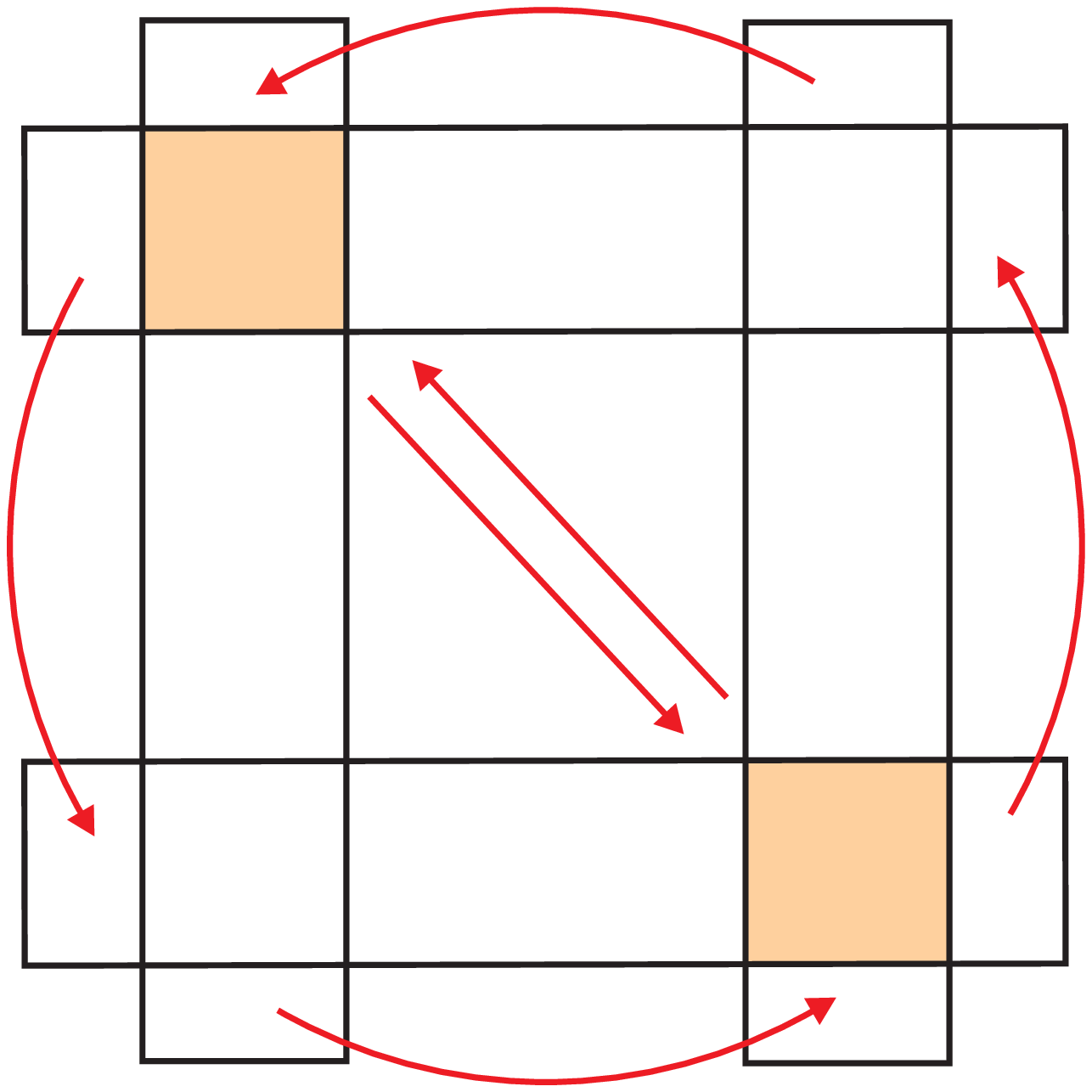 scaled 450}}
\rput(5.4,5.3){$X\stackrel{e}{\ra}X$}
\rput(8.65,1.9){$Y\stackrel{f}{\ra}Y$}
\rput(8.65,5.3){$Y\stackrel{s^*}{\ra}X$}
\rput(5.4,1.9){$X\stackrel{s}{\ra}Y$}
\rput*(7.025,3.5){${\red\cong}$}
\rput*(7,6.4){${\red (-)s}$}
\rput*(7,0.6){${\red (-)s^*}$}
\rput*(4.1,3.5){${\red s(-)}$}
\rput*(9.9,3.5){${\red s^*(-)}$}
\rput(5.4,4.4){$G_e$}\rput(8.6,2.6){$G_f$}
\rput(7.4,4){${\red s^*(-)s}$}
\rput(6.7,3){${\red s(-)s^*}$}
}
\end{pspicture}
\caption{Maximal subgroups are isomorphic.}
  \label{fig:maximal_subgroups_isomorphic}
\end{figure}

An important consequence is that the maximal
subgroups in a fixed $\JJ$-class are isomorphic. 
Again we see it in $I_n$; let $e=\id_X$ and $f=\id_Y$ be idempotents in
the $\JJ$-class 
$\JJ_e=\JJ_f$ and 
let $G_e,G_f$ be the corresponding maximal subgroups -- see Figure
\ref{fig:maximal_subgroups_isomorphic}. Somewhat incidentally,
$G_e\cong\Symm\cong G_f$ with $|X|=m=|Y|$, but this isomorphism
also arises naturally as follows. 
``Complete the square'' of $\HH$-classes that has the maximal
subgroups at its diagonal corners, and fix a representative $s:X\ra Y$
of the $\HH$-class lying the the same column as $G_e$ and row as
$G_f$.  The inverse $s^*:Y\ra X$ then lies in the diagonally
opposite $\HH$-class.

Any $h:Y\ra Y$ in the group $G_f$ can now
be decomposed as:
$$
Y\stackrel{h}{\ra}Y=Y\stackrel{s^*}{\ra}X\stackrel{g}{\ra}X\stackrel{s}{\ra}Y
$$
for some $g\in G_e$, and the map $g\mapsto sgs^*$ is a
homomorphism $G_e\ra G_f$ with inverse the map $h\mapsto s^*hs$. 

\paragraph{The inverse monoids $\sgl$.} 
Mathematics contains many examples of a group
acting on a poset or a lattice -- we now describe an
inverse monoid that wraps up the group, the lattice and the
action into a single object. It turns out that many naturally occurring
inverse monoids arise this way. It is also a
particularly useful format for understanding
their representations -- we will find it essential for the examples of
\S\ref{section:sexy:example}.

Let $G$ be a finite group and
$L$ a finite lattice -- a poset in which every pair of elements $a,b$ has a
greatest lower bound, or \emph{meet\/} $a\wedge b$, and a least upper
bound, or \emph{join\/} $a\vee b$. Suppose that $G$ acts on $L$: 
each $g\in G$ gives rise to a poset map $a\mapsto g\cdot a$, so that if $a\leq
b$ in $L$ then $g\cdot a\leq g\cdot b$. As this must also be true for $g^{-1}$, we
have $a\leq b$ iff $g\cdot a\leq g\cdot b$. 

We form a semigroup $\sgl$ out of this input data: the elements
have \emph{expressions\/} of the form $g_a$ where $g\in G$ and $a\in L$. Two
different expressions can represent the same element:
\begin{equation}
  \label{eq:12}
  g_a=h_b\text{ in }\sgl
\text{ iff }
a=b\text{ and }
g^{-1}h\cdot c=c
\text{ for all }
c\leq a
\end{equation}
Finally, the product is given by
\begin{equation}
  \label{eq:13}
  g_ah_b=(gh)_{h^{-1}\cdot a\wedge b}
\end{equation}
where $gh$ is the product of $g$ and $h$ in $G$. If it seems a little
mysterious, you can think of $g_a$ as the element of the symmetric
inverse monoid on the set $L$ and having domain the interval $L_{\leq
  a}=\{c\in L:c\leq a\}$ with effect the restriction of $g$ to this
interval. Then (\ref{eq:13}) is just the composition of 
partial bijections for $I_{L}$ and (\ref{eq:12}) warns us that
different elements of $G$ can restrict to the same partial bijection
in $\sgl$.

As $L$ is a finite lattice it has a maximum $\1=\bigvee_{a\in L} a$
and a minimum $\0=\bigwedge_{a\in L}a$, hence $\sgl$ has an identity
$\id_{\1}$, where $\id$ is the identity in $G$, and a zero $g_{\0}$,
for any $g\in G$ (as $g\cdot\0=\0=h\cdot\0$ for any other $h$, we have
by (\ref{eq:12}) that $g_{\0}=h_{\0}$). More significantly, $g_a$ has the
semigroup inverse $g_a^*=g^{-1}_{g\cdot a}$ so that $\sgl$ is an inverse
monoid. 

The Green's relation structure follows the dictates of the symmetric
inverse monoid on $L$: we have $g_a\LL h_b$ exactly when $a=b$ and
$g_a\RR h_b$ when $g\cdot a=h\cdot b$. In particular the $\LL$-class of $g_a$
consists of all the $h_a$ as $h\in G$ varies, and the $\RR$-class of all
the $h_{h^{-1}\cdot a}$. 

The $\HH$-class of $g_a$ consists of the $h_a$ for those $h\in G$ such
that $h\cdot a=g\cdot a$. The $\JJ$-classes correspond to the orbits of the
$G$-action on $L$; if $\{a_1,\ldots,a_m\}$ is an orbit, then the
eggbox decomposition of the corresponding $\JJ$-class has rows and
columns indexed by the $a_i$ and the
$\JJ$-class consists of all the $g_a$ where $g\in G$ and $a$ is one of the $a_i$.

To complete the strategic picture, let $\JJ_1$ and $\JJ_2$ be
$\JJ$-classes corresponding to the $G$-orbits $\{a_1,\ldots,a_m\}$ and
$\{b_1,\ldots,b_\ell\}$. Then $\JJ_1\leq\JJ_2$ exactly when $a_i\leq b_j$
in $L$ for some $a_i$ and some $b_j$ (or equivalently
$a_i\leq b_j$ for any $a_i$ and some $b_j$, or, any $b_j$ and some
$a_i$). 

The idempotents of $\sgl$ are the $\id_a$ for $\id$ the
identity of $G$, and the units are the $g_{\1}$ with $\1$ the maximum of
$L$. By (\ref{eq:12})-(\ref{eq:13}) the units form an isomorphic copy of $G$ in $\sgl$.

The maximal subgroup $G_a$ containing the idempotent $\id_a$ consists of the
$g_a$ with $g\cdot a=a$, subject to our running ambiguity
(\ref{eq:12}). It turns out that the ambiguity can be easily ironed out:
let $G^a$ be those elements of $G$ with $g\cdot a=a$ and let
$G^{\leq a}$ be those elements of $G$ with $g\cdot c=c$ for all $c\leq
a$. Then $G^{\leq a}$ is a normal subgroup of $G^a$ and there is an isomorphism
\begin{equation}
  \label{eq:17}
G_a\cong G^a/G^{\leq a}
\end{equation}
from the maximal subgroup $G_a$ to this
sub-quotient of the group of units $G$. 

Let $s=g_a$ be an element in the $\LL$-class of the maximal subgroup
$G_a$ and $s(-):G_a\rightarrow\HH_s$ the bijection promised by
Green's lemma. We can make the decomposition (\ref{eq:1}) quite
explicit: if $t=h_a$ is another element of $\HH_s$ then 
\begin{equation}
  \label{eq:14}
  t=s\cdot g^{-1}_{b}h_a
\end{equation}
where $b=g\cdot a$ and
$g^{-1}_{b}h_a\in G_a$. 

\begin{example}
We can shoehorn $I_n$ into this setting: $G=\Symn$ and $L$ is the
lattice of subsets of $[n]$ ordered by inclusion with $\Symn$ acting on
$L$ in the obvious way. We leave the reader to show that
$S\kern-1pt(\Symn,L)\cong I_n$ via the map that sends $g_a$ to the partial
permutation obtained by restricting $g$ to the subset $a$. 
\end{example}


\begin{vexercise}
\label{exercise:sgl:order}
If $[G:H]$ is the index of the subgroup $H$ in $G$, show that
${\displaystyle |\sgl|=\sum_{a\in L}[G:G^{\leq a}]}.$
Hence $|I_n|=\sum_{X\subseteq [n]}[\Symn:\SymX]$, where $\SymX$ is the
symmetric group on the set $X$. 
\end{vexercise}

\begin{figure}
  \centering
\begin{pspicture}(0,0)(14,13)
\rput(0,3){
\rput(1,9){(a).}
\rput(3,4.75){\BoxedEPSF{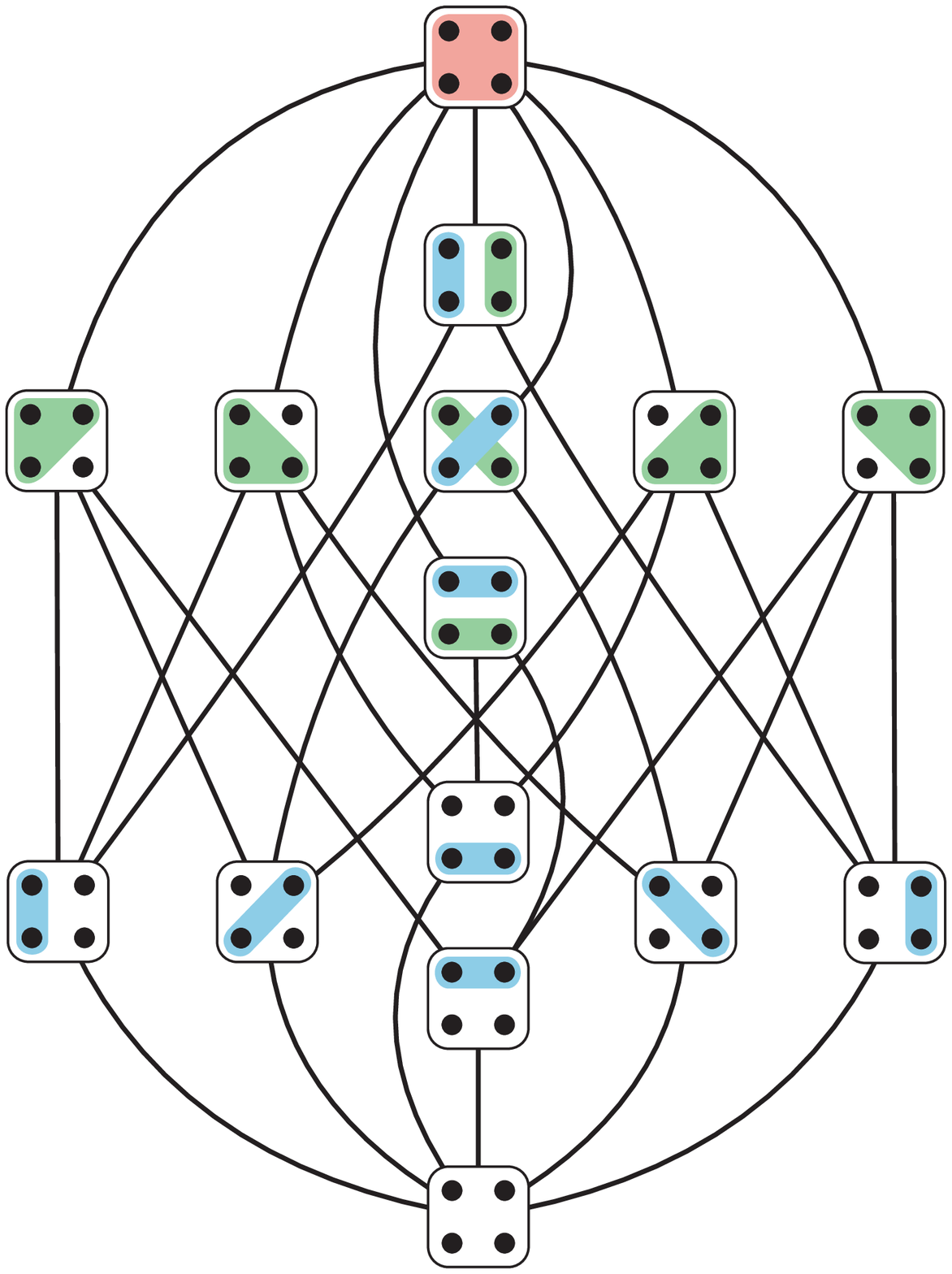 scaled 450}}
}
\rput(-1,-2){
\rput(7,4){(b).}
\rput(8.5,4.75){\BoxedEPSF{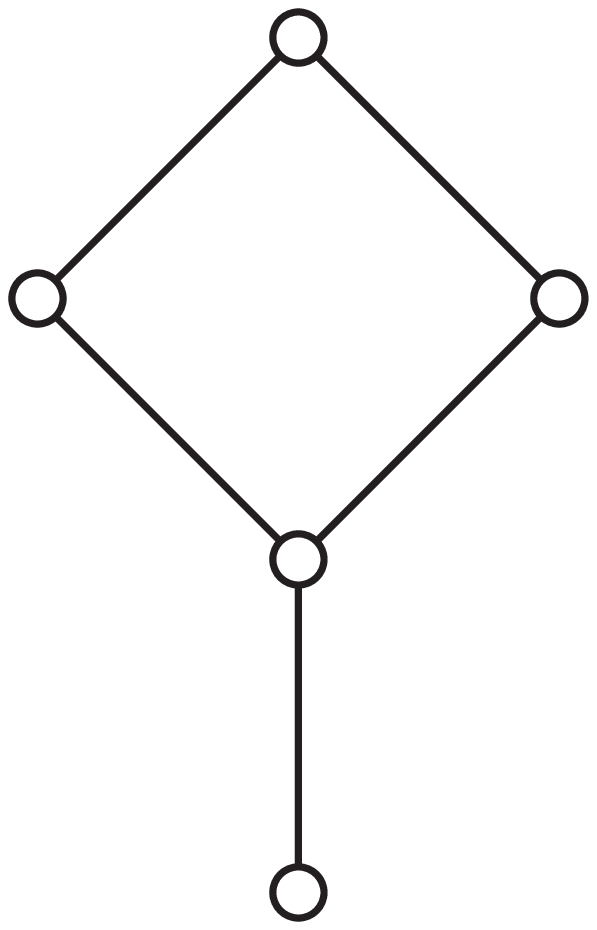 scaled 450}}
\rput(8.5,2.5){$\{a,b,c,d\}$}
\rput(9.25,4.3){$\{ab,c,d\}$}
\rput(6.6,5.5){$\{abc,d\}$}
\rput(10.4,5.5){$\{ab,cd\}$}
\rput(8.5,7){$\{abcd\}$}
}
\rput(8,3){
\rput(3,4.75){\BoxedEPSF{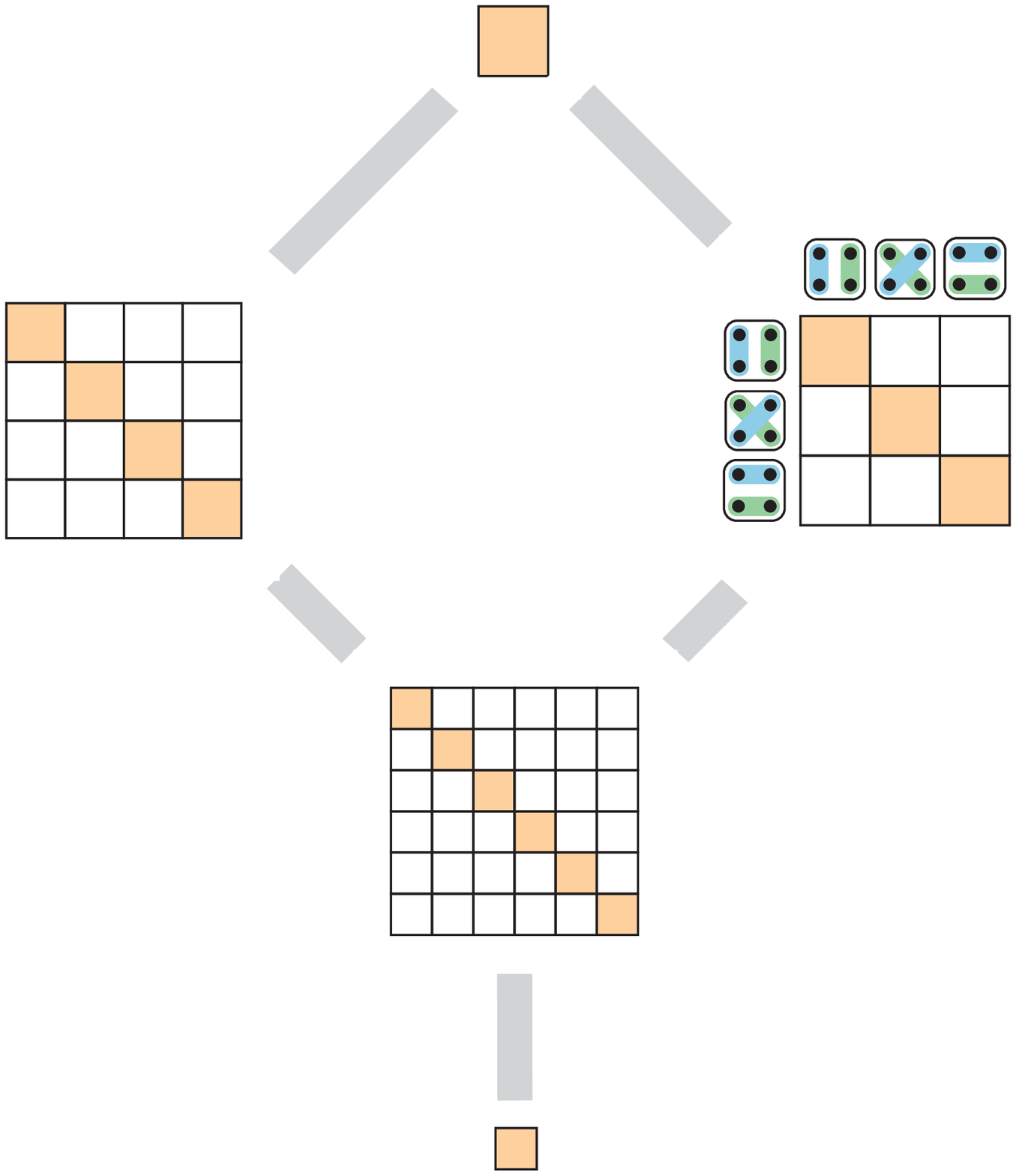 scaled 450}}
\rput(3.1,9.1,){$\{abcd\}$}
\rput(3.1,9.4){$\Symfour$}
\rput(0,8){(c).}
\rput(0.25,7){$\{abc,d\}$}
\rput(0.25,4.8){$\Symthree\times\Symone$}
\rput(3.1,4.2){$\{ab,c,d\}$}
\rput(5,3){$\Symtwo\times\Symone\times\Symone$}
\rput(3.1,0.2){$\{a,b,c,d\}$}
\rput(3.1,-0.2){$\Symone\times\Symone\times\Symone\times\Symone$}
\rput(6,4.9){$\Symtwo\times\Symtwo$}
}
\end{pspicture}
\caption{Strategic picture of $\sgl$ when $G=\Symfour$ and
  $L=\Pi(4)$ from Example \ref{exercise:partition}: (a). the partition
  lattice $\Pi(4)$; (b). the poset of
  $\JJ$-classes labelled by the type of partitions; (c).
  the strategic picture.}
  \label{fig:S(G,L)}
\end{figure}

\begin{example}
\label{exercise:partition}
An important lattice in combinatorics is the \emph{partition
  lattice\/} $\Pi(n)$, having elements 
the partitions
$\Lambda=\{\Lambda_1,\ldots,\Lambda_p\}$ of $[n]$ ordered by 
$\{\Lambda_1,\ldots,\Lambda_p\}\leq\{\Delta_1,\ldots,\Delta_q\}$ iff
each $\Lambda_i$ is a subset of some $\Delta_j$. The symmetric group
$\Symn$
acts on $\Pi(n)$ via $g\cdot
\{\Lambda_1,\ldots,\Lambda_p\}=\{g\Lambda_1,\ldots,g\Lambda_p\}$. 
The resulting $\sgl$ is called the monoid of uniform block
permutations and the strategic picture, when $n=4$, is in Figure \ref{fig:S(G,L)}.
The $\JJ$-class poset is the poset of partitions
$\lambda=\{\lambda_1,\ldots,\lambda_p\}$ of the integer $n$ (see the
beginning of the Interlude),
and the corresponding maximal subgroup is isomorphic to the Young subgroup 
$S_{\kern-1pt\lambda_1}\times\cdots\times
S_{\kern-1pt\lambda_p}$. In particular, the order of the monoid of
uniform block permutations is 
$\sum_{\Lambda\in\Pi(n)}[\Symn:S_{\lambda_1}\times\cdots\times S_{\lambda_p}]$.
\end{example}


\section{Representations}
\label{section:representations}

This section contains the basics of representation theory that are
common to all
finite regular monoids. The theme is the extent to which
representations can be decomposed into ``atomic'' pieces. These can
then be reassembled to get a handle on the sociology of the
representations of a semigroup. It turns out
that this is almost always possible for groups and inverse monoids,
but less so for regular, non-inverse monoids. 

Throughout, $k$ is a field and $V$ a finite dimensional vector space over
$k$. Let $End(V)$ be the monoid, under composition, of all vector space
homomorphisms (or linear 
maps) $V\ra V$.

An \emph{$S$-action on $V$\/} or \emph{linear representation of $S$\/}
is a monoid homomorphism 
$$
\varphi:S\ra End(V).
$$ 
We adopt the convention that all monoid homomorphisms send $1$'s to $1$'s,
so that $\varphi(1_S)$ is the identity homomorpism $\id:V\ra V$. In
particular
$\im\varphi\not=\{0\}$, and so our representations are not
\emph{null\/}. If $S$ is a group then necessarily 
$\im\varphi\subset GL(V)$,  the group of vector space isomorphisms
(or invertible linear maps) $V\ra V$. The notion of a semigroup representation is thus
a straight generalisation of that of a group representation. 

We will identify $s\in S$ and $\varphi(s)\in
End(V)$, so that if $v\in V$, we just write $s\cdot v$, or even $sv$, for the effect
of the linear map $\varphi(s)$ on the vector $v$. Mostly we will just write
$V$ for an $S$-representation without explicit reference to the
action. 

The following representation of our three running examples $\Symn,I_n$
and $T_n$ will turn out to display the full range of possible behaviours:

\begin{example}[mapping representations]
\label{example:mapping:representation}
Fix a  basis $\{v_1,\ldots,v_n\}$ for the space $V$ and for $s\in
\Symn,I_n$ or $T_n$ define
\begin{equation}
  \label{eq:2}
s\cdot v_i=v_{s(i)}
\quad (s\in \Symn, T_n)
\quad
\text{ or }
\quad
s\cdot v_i
=
\left\{\begin{array}{ll}
 v_{s(i)}, &\text{ if }i\in\dom s\\
0,&\text{else.}
\end{array}\right.
\quad (s\in I_n)  
\end{equation}
and then extend linearly. To analyse the structure of the mapping
representations, we need to know how to decompose representations in general.
\end{example}

\paragraph{Sub-representations and reducibility.} These allow
us to understand representations in the large. If $V$ is a
representation and $U$ is a subspace left invariant 
by the $S$-action, i.e. $SU=U$, then we call $U$ an
\emph{($S$-)subrepresentation\/} of 
$V$. The quotient space $V/U$ then carries an $S$-action via
$s\cdot(v+U)=sv+U$, well-defined, as $sU=U$. There is then a 1-1
correspondence between the subrepresentations of $V/U$ and the
subrepresentations $W$ of $V$ such that $U\subseteq W\subseteq V$.

If $V$ has a proper, non-zero subrepresentation $U$, then call $V$
\emph{reducible\/}; $V$ is \emph{irreducible\/} if the only
subrepresentations are
$\{0\}$ and $V$. 

A subrepresentation $U$ of $V$ is \emph{maximal\/} when $U\not=V$ but
for any subrepresentation $W$ of $V$ with $U\subseteq W\subseteq V$ we
have either $W=U$ or $W=V$. Because of the 1-1 correspondence
mentioned above, $U$
is maximal exactly when $V/U$ is irreducible.

If $U,W$ are subrepresentations of $V$ such that $V=U\oplus W$ as
vector spaces, then the representation $V$ is the \emph{(internal)
  direct sum\/} of $U$ and $W$. Externally, if $U$ and $W$ are
arbitrary $S$-representations then the vector space direct sum
$U\oplus W$ carries an $S$-action via $s\cdot(u+w)=su+sw$, and
$U\oplus W$ is the \emph{(external) direct sum\/} of $U$ and $W$.

\begin{example}[the permuting coordinates and reflectional representations of
  $\Symn$]
\label{example:permutation:representation}
  The $\Symn$-action in (\ref{eq:2}) is by ``permuting coordinates''
  (or more accurately, permuting basis vectors). In particular,
  for $(i,j)\in\Symn$ the resulting 
  isomorphism $V\ra V$ is the reflection in the hyperplane with
  equation $x_i-x_j=0$.
As $\Symn$ is generated by the transpositions, its image in $GL(V)$
is generated by reflections, i.e: is a \emph{reflection group\/}.

\begin{figure}
  \centering
\begin{pspicture}(0,0)(14,4)
\rput(1.5,0){
\rput(0,2){\BoxedEPSF{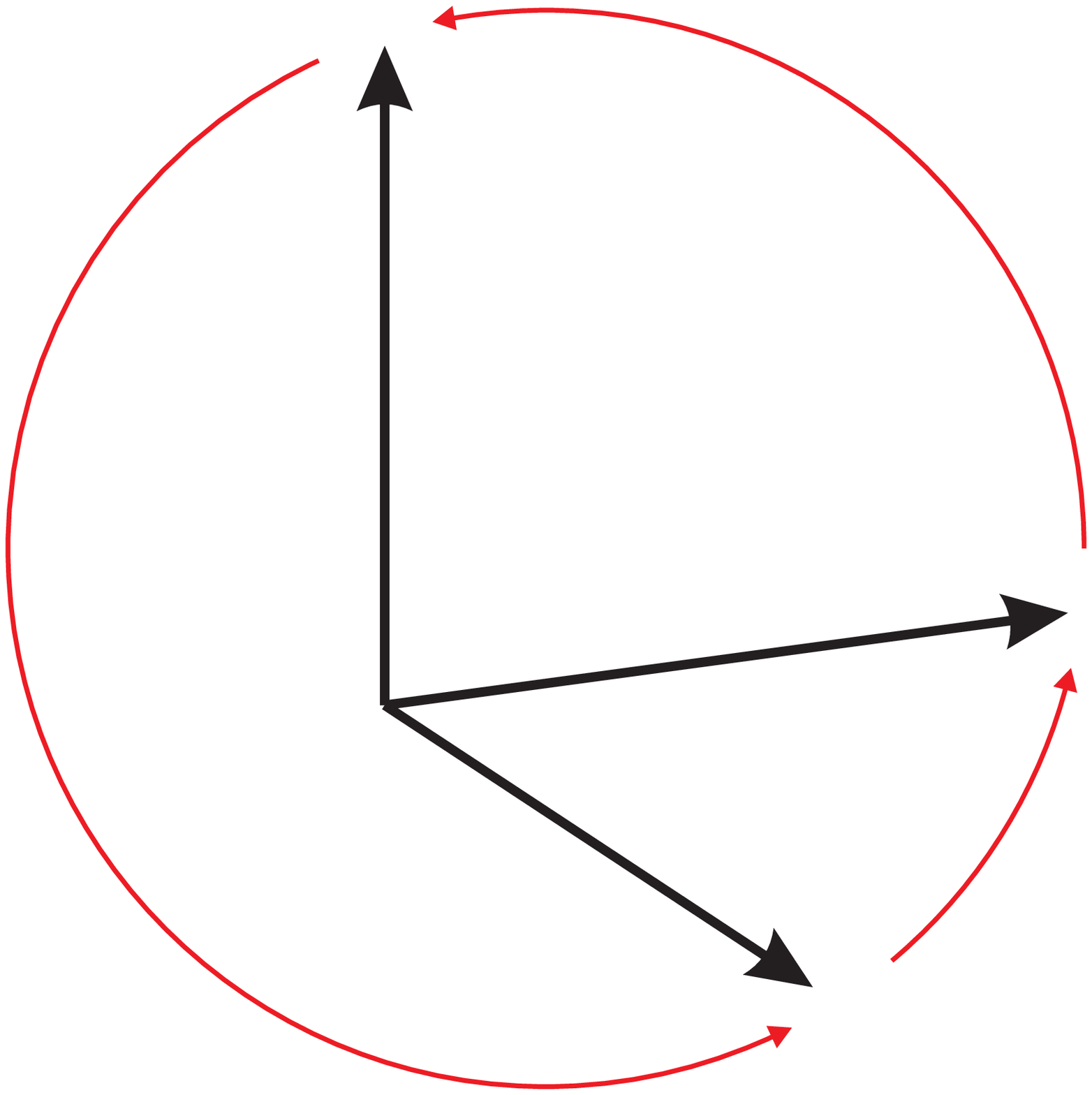 scaled 225}}
\rput(0.5,2.5){${\red s=(1,2,3)}$}
\rput(1.05,0.45){$v_1$}\rput(1.95,1.8){$v_2$}\rput(-0.6,3.85){$v_3$}
}
\rput(0,0){
\rput(5.75,2){\BoxedEPSF{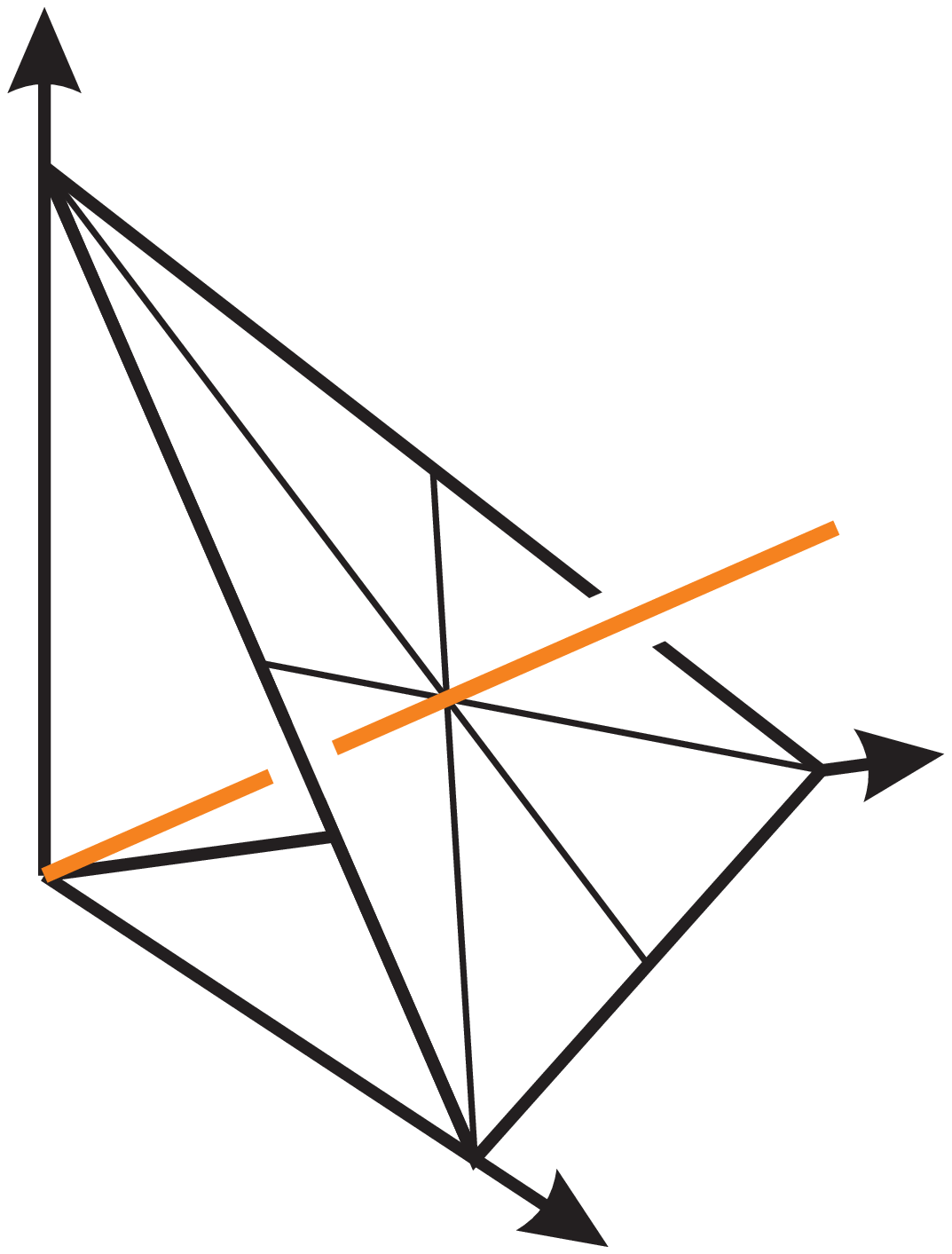 scaled 225}}
\rput(6.85,2.3){$U$}
}
\rput(0,0){
\rput(9,2){\BoxedEPSF{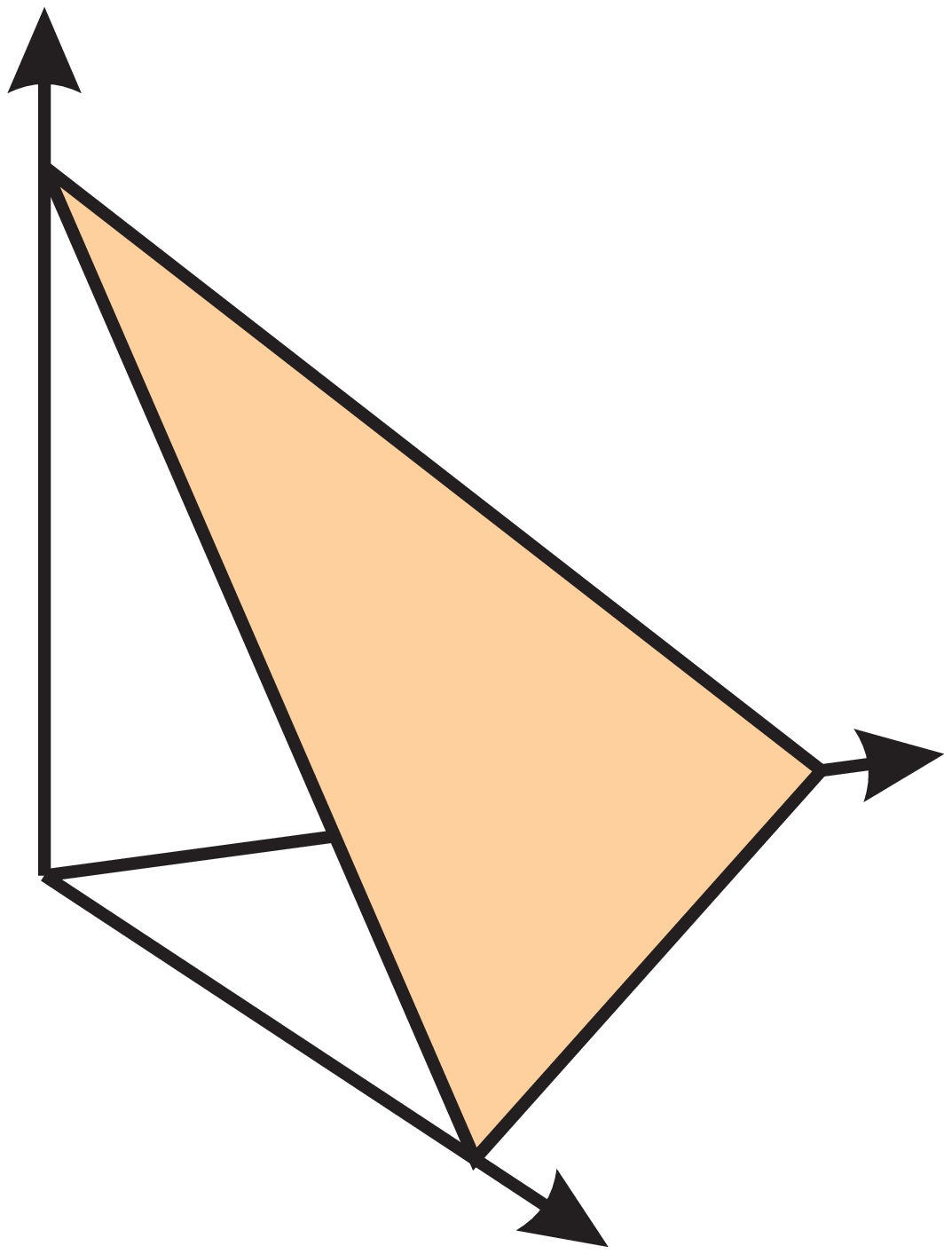 scaled 225}}
\rput(9.5,2.5){$W+\frac{1}{n}u$}
}
\rput(12.5,0){
\rput(0,2){\BoxedEPSF{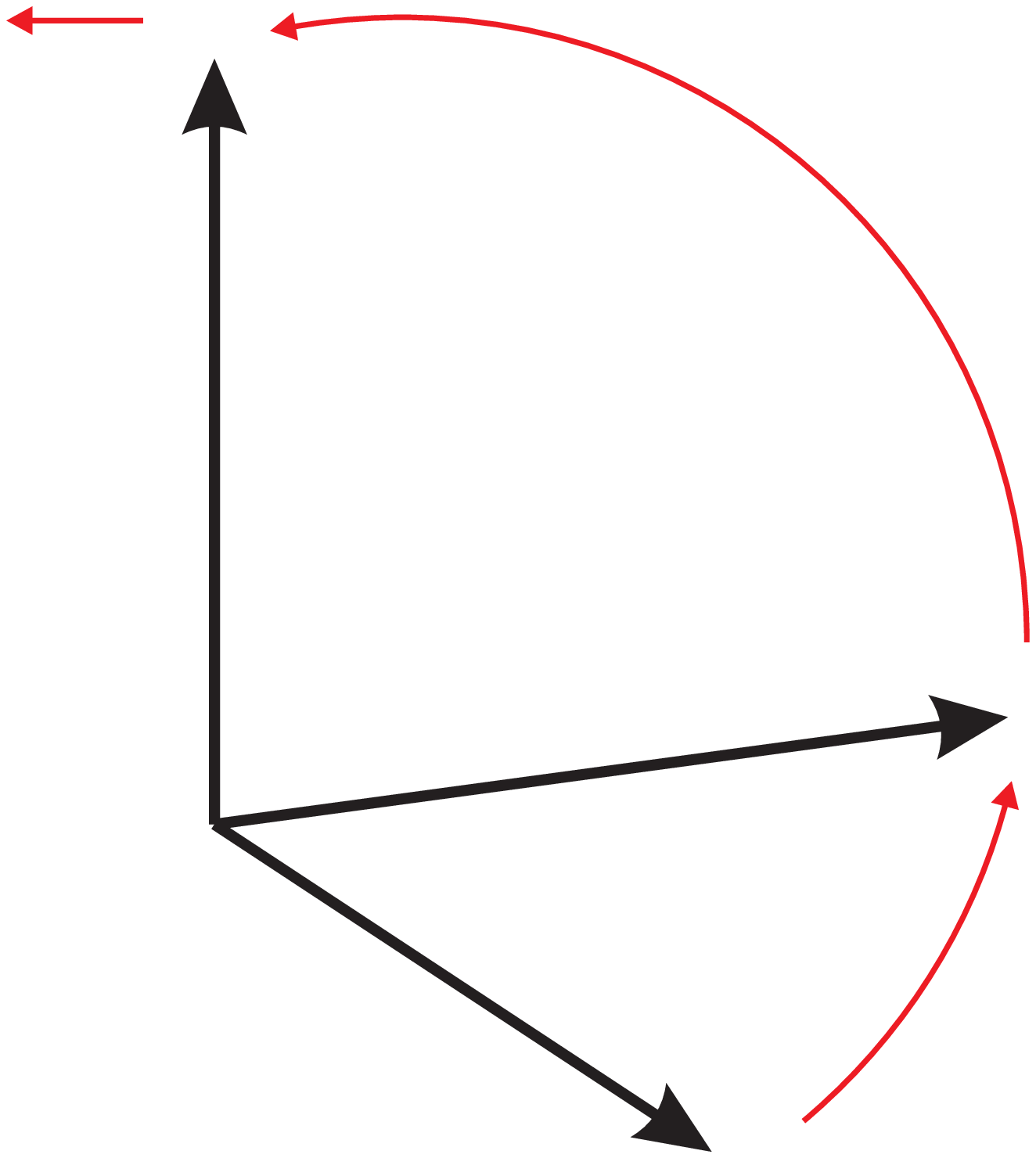 scaled 225}}
\rput(0.1,2.5){${\red s=[1,2,3]}$}
\rput(-0.3,-0.175){\rput(1.05,0.45){$v_1$}
\rput(1.95,1.8){$v_2$}
\rput(-0.6,3.85){$v_3$}  
\rput(-1.3,3.85){$0$}}
}
\end{pspicture}
\caption{\emph{From left to right}: Permutation action of $\Symn$; the
  line $U$ which is the $k$-span of $u=v_1+\cdots+v_n$; the (affine) 
  hyperplane $W+\frac{1}{n}u$ coming from the reflectional representation; the partial
  reflection action of $I_n$. The pictures are for $n=3$
  and the notation for partial
  permutations is described in the Notes and References section.}
  \label{fig:reflection_rep_Sn}
\end{figure}

The vector $u=v_1+\cdots+v_n$  is fixed by any permutation
in $\Symn$, so that if 
$U$ is the $k$-span of $u$ then
$\Symn U=U$, a subrepresentation. As each vector in $U$ is fixed by
every element of $\Symn$,
this is the \emph{trivial
  representation\/} of $\Symn$ (see Figure \ref{fig:reflection_rep_Sn}).

Thus if $n>1$ then the permuting coordinates representation $V$ is
reducible. Moreover, as $U$ is $1$-dimensional it has only the two
subspaces, $\Symn$-invariant or otherwise, namely $\{0\}$ and $U$.
Hence $U$ is irreducible. When $n=1$ we have $V=U$ is irreducible. 

Now let $W$ be the hyperplane with equation $x_1+x_2+\cdots+x_n=0$,
that is, the
set of 
points whose coordinates with respect to the $v_i$ sum to
$0$. Permuting the coordinates of such 
a vector doesn't change the
fact that they add to $0$, hence $W$ is also a
subrepresentation of $V$. Figure \ref{fig:reflection_rep_Sn} has the
plane $W$ when $n=3$, shifted off the origin to make it  easier to
see. For reasons that are maybe a little obscure at the moment, $W$ is
called the \emph{reflectional representation\/} of $\Symn$. 

Moreover, if the characteristic $\text{char}(k)$ of the
field does not divide $\dim V=n$, then $W$ is irreducible. 
For suppose that $X\not=\{0\}$ is a $\Symn$-invariant subspace of $W$
and let
$v\in X$ with $v\not=0$. 
If all the coordinates of $v$ are equal to some $\lambda\in k$, then these
sum to $0$ to give
$n\lambda=0$, hence -- by the restriction on the characteristic --
we must have $\lambda=0$, hence $v=0$, a contradiction. The vector $v$ must
therefore have two
coordinates that are different. For each $1\leq i <n$ we can engineer a $g_i\in\Symn$
such that in the vector
$g_i\cdot v$ it is the $i$-th and $(i+1)$-st coordinates that are different.
Then $(i,i+1)g_i\cdot
v-g_i\cdot v$ is a non-zero multiple of $v_i-v_{i+1}$. As $X$ is
$\Symn$-invariant we conclude that for each $1\leq i <n$
the vector $v_i-v_{i+1}$ is an element of $X$. But these vectors form a basis
for $W$, so $X=W$, and $W$ is irreducible as claimed.

The permuting coordinates representation of $\Symn$ can thus be decomposed 
$V=U\oplus W$ into the trivial and reflectional representations, with
both of these
irreducible.  
\end{example}

\begin{example}[the partial reflectional representation of $I_n$]
\label{example:partial:reflection:representation}
The $I_n$-action (\ref{eq:2}) is by partial permutations of
coordinates and
the image of $I_n$ in $End(V)$ is a \emph{reflection monoid\/}.

Assume that $n>1$.
The line $U$ spanned by $u=v_1+\cdots+v_n$ is no longer $I_n$-invariant: 
if $s\in I_n$ is the partial identity with domain $\{1\}$
then $s\cdot u=v_1\not\in U$. Similarly $W$ is not 
$I_n$-invariant. 

In fact, $V$ itself is irreducible: for suppose $U$ is an
$I_n$-invariant subspace containing $0\not=v\in U$ with
$v=\sum \lambda_iv_i$ and $\lambda_j\not=0$ for
some $j$. For each $1\leq i\leq n$ let $s_i\in I_n$ be the partial
permutation shown in Figure \ref{fig:In_irreducible}. 
\begin{figure}[h]
  \centering
\begin{pspicture}(0,0)(14,3)
\rput(7,1.5){\BoxedEPSF{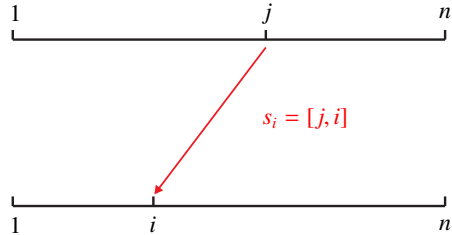 scaled 450}}
\rput(4.2,2.9){$1$}\rput(7.5,2.9){$j$}\rput(9.85,2.9){$n$}
\rput(4.2,0.1){$1$}\rput(6,0.1){$i$}\rput(9.85,0.1){$n$}
\rput(8,1.5){${\red s_i=[j,i]}$}
\end{pspicture}
\caption{$s_i\in I_n$ for $1\leq i\leq n$; the notation $[j,i]$ is
  explained in the Notes and References section.}
  \label{fig:In_irreducible}
\end{figure}
Then $s_i\cdot v=\lambda_jv_i$, hence
$v_i\in U$ for all
$i$, and so $U=V$. 
\end{example}

\begin{example}[the mapping representation of $T_n$]
\label{example:mapping:representation:T_n}
Again suppose that $n>1$. The hyperplane $W$ goes back to being a
subrepresentation of $V$: if $w\in W$ with $w=\sum \lambda_i
v_i$ where $\sum\lambda_i=0$, then $sw$ for $s\in T_n$ is shown in
Figure \ref{fig:Tn_mapping_rep}. In particular the non-zero
coordinates of $sw$ are sums of the coordinates of $w$, and so still
sum to $0$:
\begin{figure}[h]
  \centering
\begin{pspicture}(0,0)(14,4)
\rput(0,0.5){
\rput(7,1.5){\BoxedEPSF{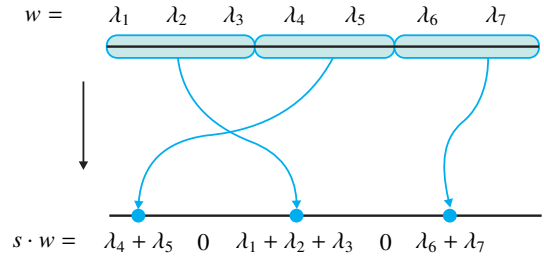 scaled 450}}
\rput(3.5,3){$w=$}
\rput(4.5,3){$\lambda_1$}
\rput(5.25,3){$\lambda_2$}
\rput(6,3){$\lambda_3$}
\rput(6.8,3){$\lambda_4$}
\rput(7.6,3){$\lambda_5$}
\rput(8.55,3){$\lambda_6$}
\rput(9.45,3){$\lambda_7$}
\rput(3.5,0){$s\cdot w=$}
\rput(4.75,0){$\lambda_4+\lambda_5$}
\rput(5.6,0){$0$}
\rput(6.8,0){$\lambda_1+\lambda_2+\lambda_3$}
\rput(8.85,0){$\lambda_6+\lambda_7$}
\rput(8,0){$0$}
}
\end{pspicture}
\caption{The hyperplane $W$ is a $T_n$-subrepresentation.}
  \label{fig:Tn_mapping_rep}
\end{figure}
i.e. $sw=\sum\mu_i v_i$ where
$\mu_i=\sum \lambda_{ij}$, the sum over the $j$ in the
fiber of $i$. As $\sum\mu_i=\sum\lambda_i=0$, we get $sw\in W$.

If $S$ is any monoid, $V$ an $S$-representation and $T$ a submonoid of $S$, then
it is easy to see that restricting the $S$-action to $T$
makes $V$ into a $T$-representation. This observation gives that
$W$ is irreducible: if $X$ is a subrepresentation of $W$ then
it is an $\Symn$-subrepresentation
of the $S_n$-representation $W$. The irreducibility of this
-- when $\text{char}(k)$ does not divide $n$ -- then gives $X=\{0\}$ or
$W$. 

Just as for $I_n$, the line $U$ spanned by $v_1+\cdots+ v_n$ is not a
subrepresentation: for
example when $s$ is the constant map in $T_n$ that
sends all of  $[n]$ to $1$, then $su=nv_1\not\in U$. 

In fact, when $n>2$ we claim that there are \emph{no\/}
$1$-dimensional subrepresentations of the $T_n$-mapping representation
$V$. For, suppose that 
$v\not=0$ 
so that $v=\sum\lambda_iv_i$ with
$\lambda_j\not=0$ for some $j$. If $s_1,s_2\in T_n$ are given by
Figure \ref{fig:Tn_elements}, where
everything apart from $j$ is sent to $n$, with $s_i$ sending $j$ to $i$,
then 
$$
s_1(v)=\lambda_jv_1+\biggl(\sum_{i\not=j}\lambda_i\biggr)v_n
\quad
\text{ and }
\quad
s_2(v)=\lambda_jv_2+\biggl(\sum_{i\not=j}\lambda_i\biggr)v_n.
$$
As these are independent, any non-trivial $T_n$-invariant subspace 
must be at least $2$-dimensional, and the claim follows.
\begin{figure}[h]
  \centering
\begin{pspicture}(0,0)(14,3)
\rput(3.5,1.5){\BoxedEPSF{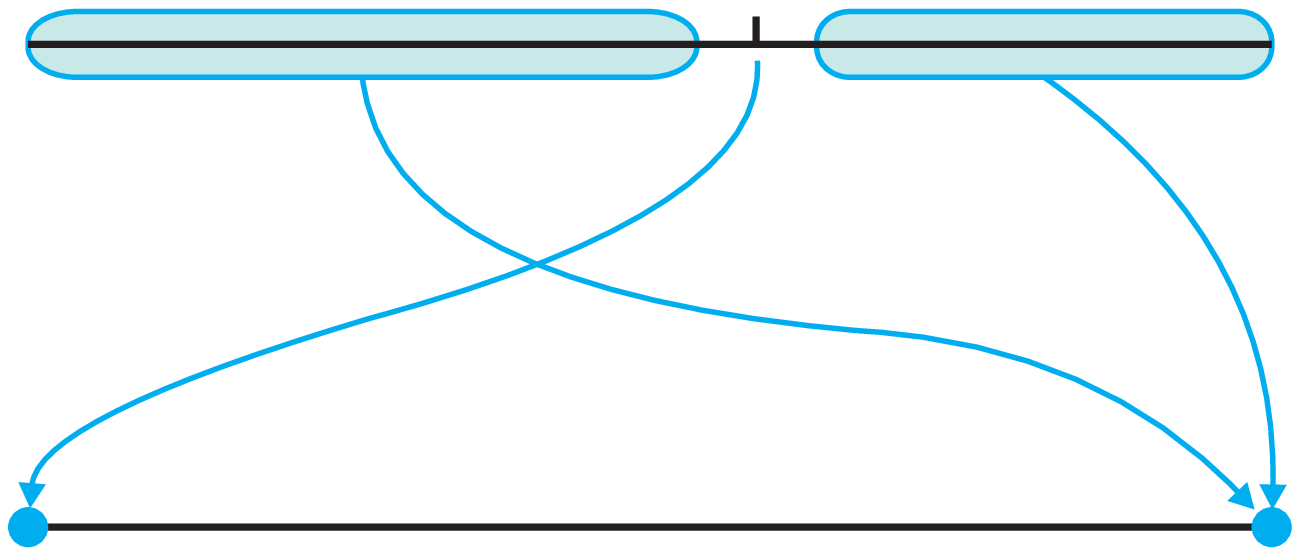 scaled 450}}
\rput(0.65,2.9){$1$}\rput(4,2.9){$j$}\rput(6.35,2.9){$n$}
\rput(7.6,2.9){$1$}\rput(10.95,2.9){$j$}\rput(13.3,2.9){$n$}
\rput(0.65,0.1){$1$}\rput(6.35,0.1){$n$}
\rput(8.35,0.1){$2$}\rput(13.3,0.1){$n$}
\rput(10.5,1.5){\BoxedEPSF{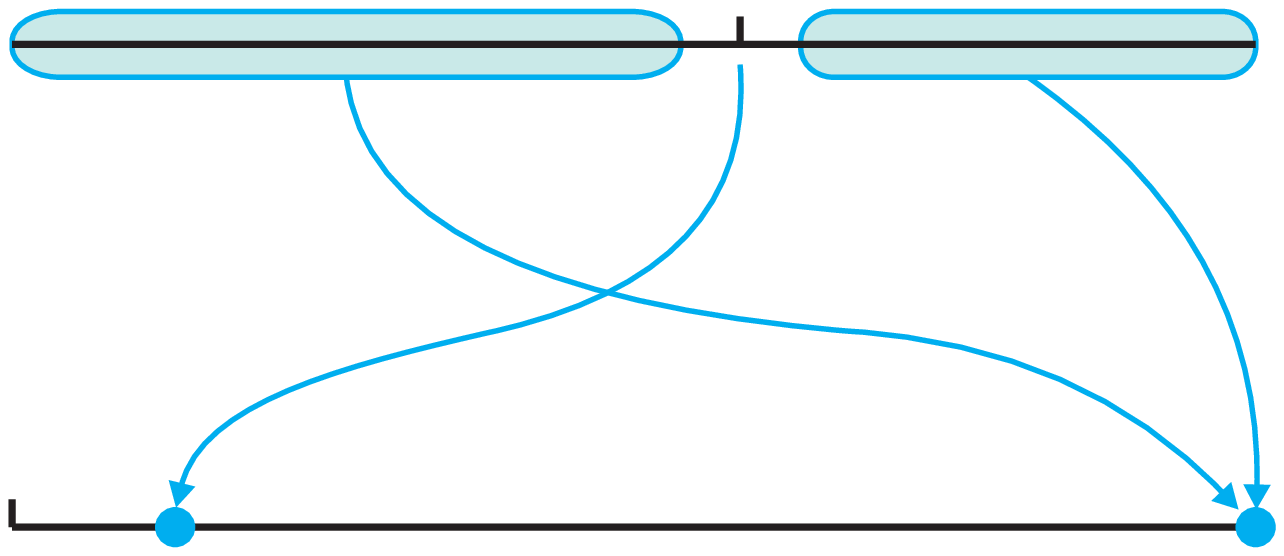 scaled 450}}
\rput(0.65,1.5){$s_1$}\rput(7.6,1.5){$s_2$}
\end{pspicture}
\caption{The $s_1$ and $s_2$ in $T_n$ for Example
  \ref{example:mapping:representation:T_n}.} 
  \label{fig:Tn_elements}
\end{figure}
\end{example}

\paragraph{Semisimplicity.}
  If $S$ is a finite regular monoid and $k$ a field, then the pair
  $(S,k)$ is \emph{semisimple\/} when every $S$-representation $V$
  over $k$ can be
  decomposed 
$$V=\bigoplus V_i$$
with the $V_i$ irreducible
  subrepresentations of $V$.

Such a decomposition is then unique in the following sense: a
\emph{morphism\/} $q:V\ra U$ of 
$S$-representations is a linear map that 
  commutes with the $S$-actions on $V$ and $U$, i.e. for all $s\in S$ the diagram
$$
\begin{pspicture}(0,0)(14,2)
\rput(5,0){
\rput(1,2){$V$}
\psline[linewidth=0.5pt]{->}(1.25,2)(2.75,2)
\rput(2,2.2){${\scriptstyle s\cdot (-)}$}
\rput(3,2){$V$}
\psline[linewidth=0.5pt]{->}(1,1.7)(1,0.3)
\psline[linewidth=0.5pt]{->}(3,1.7)(3,0.3)
\rput(0.8,1){${\scriptstyle q}$}\rput(3.2,1){${\scriptstyle q}$}
\rput(1,0){$U$}
\psline[linewidth=0.5pt]{->}(1.25,0)(2.75,0)
\rput(2,0.2){${\scriptstyle s\cdot (-)}$}
\rput(3,0){$U$}
}
\end{pspicture}
$$
commutes, where the top $s\cdot(-)$ is the $S$-action on $V$
and the bottom is the $S$-action on $U$. Call $q$ an
\emph{isomorphism\/} if it is a bijective morphism. Uniqueness then
means:

\begin{theorem}[Jordan-H\"{o}lder]
\label{theorem:JordanHolder}
 Let $V$ be an $S$-representation and $V=\bigoplus V_i$ with the
 subrepresentations $V_i$ irreducible. If $W$ is an irreducible
 subrepresentation of $V$ then $W$ is isomorphic to one of the $V_i$. 
\end{theorem}

We saw at the end of Example \ref{example:permutation:representation} that
the reflectional representation of $S_n$ can be so decomposed, and indeed:

\begin{theorem}[Mashke]
\label{theorem:Mashke}
  If $S$ is a finite group then $(S,k)$ is semisimple if and only if
  the characteristic
  $\text{char}(k)$ does not divide the order of $S$.
\end{theorem}

In particular, $(\Symn,k)$ is semisimple exactly when $\text{char}(k)$ doesn't
divide $n!$, so that characteristic $0$ representations can always be
decomposed. The situation for inverse monoids is almost as good:

\begin{theorem}[Munn-Oganesyan]
\label{theorem: Munn:Oganesyan}
  If $S$ is a finite inverse monoid then $(S,k)$ is semisimple if
  and only if $\text{char}(k)$ does not divide the order of any
  subgroup of $S$.
\end{theorem}

As any subgroup of $S$ is in turn a subgroup of a maximal subgroup $G_e$ for
some idempotent $e$, it suffices that the characteristic does not
divide the order of any $G_e$.

For our model inverse monoid $I_n$, we have already seen that 
the maximal subgroups
are isomorphic to $\Symm$ for $1\leq m\leq
n$. The pair $(I_n,k)$ is thus semisimple when the
characteristic of $k$ does not divide $m!$ for any $m\leq n$,
i.e. when it does not divide $n!$ We therefore get the same condition for
the semisimplicity of $I_n$ and $\Symn$ representations.

For $T_n$ things are not so good. 
If the mapping representation $V$ of $T_n$
is decomposable $V=\bigoplus V_i$ with the $V_i$
irreducible, then by Theorem
\ref{theorem:JordanHolder}, one of the $V_i$ is isomorphic to the
representation on the hyperplane $W$ with equation $x_1+\cdots
x_n=0$ that we saw above. The
decomposition of $V$ must then be $V\cong
W\oplus W'$ with $W'$ a $1$-dimensional subrepresentation. But we have
seen for $n>2$ that there are no $1$-dimensional subrepresentations of
$V$, and so no such decomposition of $V$ can exist when $n>2$.

Thus the pair $(T_n,k)$ is not semisimple, when $n>2$,
for \emph{any\/} $k$ whose characteristic does not divide $n$. In
particular, not even for $k$ of characteristic $0$.

\paragraph{Here then is what we have learned from the three examples:} 
in characteristic $0$ the partial permuting coordinates (or 
reflectional)
representation of $I_n$ is ``atomic''; the permuting coordinates representation
of $\Symn$ is not atomic but can be decomposed into
pieces that are; the mapping representation of $T_n$ is not atomic and 
cannot even be decomposed into atomic pieces.

\section*{Interlude: the symmetric group}

The moral of \S\S\ref{section:reduction}-\ref{section:clifford:munn}
is that the representations of a 
(finite regular) monoid $S$ are largely driven by the representations
of its maximal subgroups. In every example that
we have seen so far these maximal subgroups have been symmetric
groups, or products of symmetric groups. It seems reasonable then to
understand better the representations of the symmetric
group, at least when $k=\C$.

We do this in a completely self-contained-tailored-to-$\Symn$
way, without any reference to the general theory of 
representations of finite groups. This will make it seem a little like
pulling a rabbit out of a hat; the reader who is interested in the
broader context of these facts should consult the Notes and References
at the end.

By Theorem \ref{theorem:Mashke}, any $\Symn$-representation over $\C$
is a direct sum of irreducible representations; we will thus content
ourselves with describing just these. Despite the comments in the previous
paragraph, we allow ourselves one general fact: the irreducible
representations over $\C$ of a finite group are in 1-1 correspondence
with the conjugacy classes of the group. For $\Symn$, these in turn
are in 1-1 correspondence with the possible cycle structures of
permutations of degree $n$ and \emph{these in turn\/} with the
partitions of the integer $n$: the integer
sequences $\lambda_1\geq\ldots\geq\lambda_p>0$ with
$\sum\lambda_i=n$. Write $\lambda=\{\lambda_1,\ldots,\lambda_p\}\vdash n$. 

\begin{figure}
  \centering
\begin{pspicture}(0,0)(14,2.5)
\rput(3,1.25){\BoxedEPSF{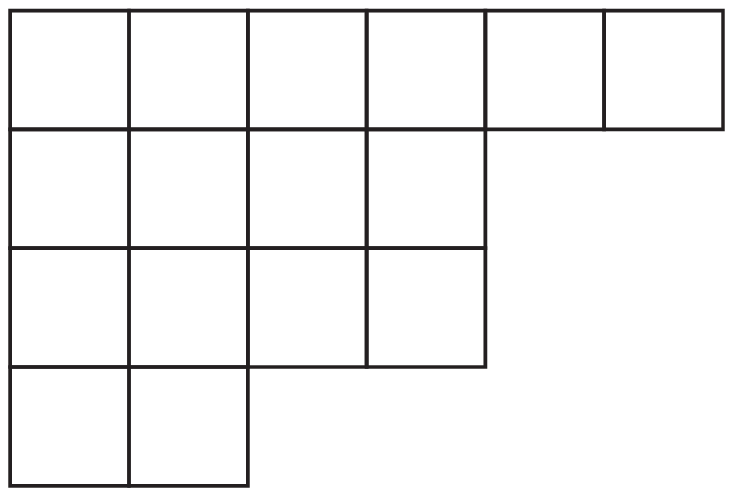 scaled 450}}
\rput(1.1,2.1){$\lambda_1$}\rput(1.1,1.5){$\lambda_2$}
\rput(1.1,1.1){$\vdots$}\rput(1.1,0.5){$\lambda_p$}
\rput(7,1.25){\BoxedEPSF{tableau3.eps scaled 450}}
\rput(5.65,2.05){$2$}\rput(6.15,2.05){$3$}\rput(6.7,2.05){$12$}
\rput(7.25,2.05){$7$}\rput(7.8,2.05){$13$}\rput(8.35,2.05){$15$}
\rput(5.65,1.55){$5$}\rput(6.15,1.55){$1$}\rput(6.7,1.55){$10$}
\rput(7.25,1.55){$14$}
\rput(5.65,1){$4$}\rput(6.15,1){$9$}\rput(6.7,1){$11$}\rput(7.25,1){$16$}
\rput(5.65,0.45){$6$}\rput(6.15,0.45){$8$}
\rput(11,1.25){\BoxedEPSF{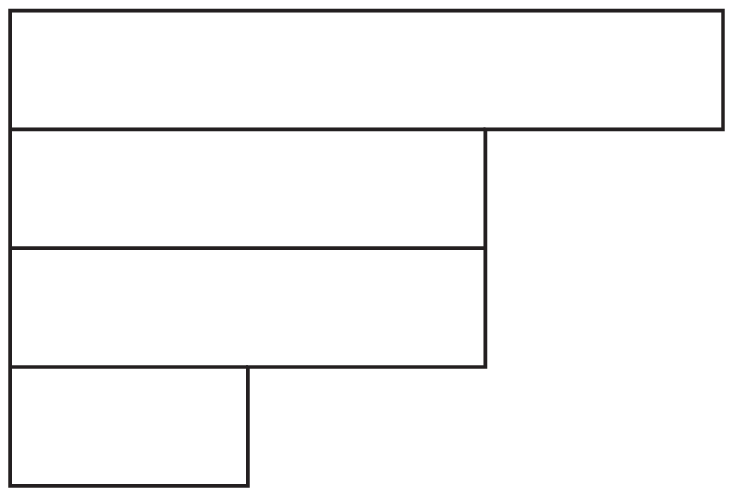 scaled 450}}
\rput(4,0){
\rput(5.65,2.05){$2$}\rput(6.15,2.05){$7$}\rput(6.7,2.05){$12$}
\rput(7.25,2.05){$3$}\rput(7.8,2.05){$13$}\rput(8.35,2.05){$15$}
\rput(5.65,1.55){$1$}\rput(6.15,1.55){$10$}\rput(6.7,1.55){$5$}
\rput(7.25,1.55){$14$}
\rput(5.65,1){$16$}\rput(6.15,1){$11$}\rput(6.7,1){$9$}\rput(7.25,1){$4$}
\rput(5.65,0.45){$8$}\rput(6.15,0.45){$6$}
}
\end{pspicture}
\caption{A Young diagram \emph{(left)}, a tableau $T$
  \emph{(middle)} and the resulting tabloid $\{T\}$ \emph{(right)}
  corresponding to a partition
  $\lambda=\{\lambda_1,\ldots,\lambda_p\}\vdash n$.} 
  \label{fig:Young:diagram}
\end{figure}

Fix $\lambda=\{\lambda_1,\ldots,\lambda_p\}\vdash n$ a partition of
$n$. A Young diagram of shape $\lambda$ illustrates the structure of
$\lambda$, as on the 
left of Figure \ref{fig:Young:diagram}, and a \emph{tableau\/} $T$ is a Young
diagram filled with entries from $[n]$ with no repeats allowed -- as in the
middle of Figure \ref{fig:Young:diagram}. The
tableau is standard, or $T$ is a \emph{standard tableau\/}, when
the entries increase along the rows and down the columns. Finally, a tableau $T$
yields a \emph{tabloid\/} $\{T\}$, which is just a tableau where we no
longer care about the ordering in the rows -- see the right of Figure
\ref{fig:Young:diagram}. 

The symmetric group $\Symn$ acts on the set of tableau of shape
$\lambda$, via $g:T\mapsto gT$ for $g\in\Symn$, where $gT$ is the
tableau that has $g(i)$ in the box in which $T$ has $i$. This action
extends to the set of tabloids of shape $\lambda$ via
$g\cdot\{T\}=\{gT\}$. For a tableau $T$, the column group $c_T$ is
defined 
$$
c_T=\{g\in\Symn:g\text{ preserves each column of }T\}\subseteq\Symn
$$
Let $M^\lambda$ be the $\C$-vector space with basis the tabloids
$\{T\}$ of the fixed shape $\lambda$. Via the action above, $\Symn$
acts on $M^\lambda$ by permuting the basis vectors. For the partition
$\lambda=\{n-1,1\}$ we will see below that $M^\lambda$ is the
permuting coordinates representation of
\S\ref{section:representations}. Now we have other representations of $\Symn$.

In any case, the $M^\lambda$ are in general reducible -- much like the
permuting coordinates representation -- and we will pass to a
particular subrepresentation.
If $T$ is a tableau then let $v_T\in M^\lambda$ be the vector
\begin{equation}
  \label{eq:19}
v_T=\sum_{h\in c_T}\text{sign}(h)\, h\cdot\{T\}
\end{equation}
where $\text{sign}(h)=1$ or $-1$ depending on whether $h$ is an even
or odd permutation. 
We will see below that in general
the vectors $v_T$, as $T$ ranges over the tableau of
shape $\lambda$, are \emph{not\/} independent. (The $v_T$ for $T$
standard \emph{are\/} an 
independent subset, although we won't need this fact here.)
In any case, let
\begin{equation}
  \label{eq:15}
 S^\lambda:=\text{ the subspace of }M^\lambda\text{ spanned by the }v_T. 
\end{equation}
It turns out that $S^\lambda$ is an irreducible subrepresentation of
$M^\lambda$, and 
as $\lambda$ varies over the partitions of $n$, the $S^\lambda$ -- 
called \emph{Specht\/} representations -- give a
complete and non-redundant list of the irreducible
$\Symn$-representations over $\C$. 

\begin{example}
\label{symmetric:example:trivial}
If  $\lambda=\BoxedEPSF{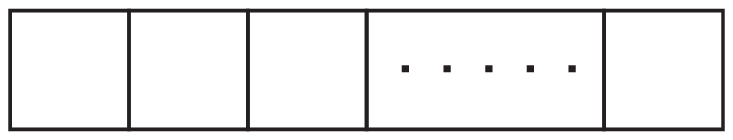 scaled 200}$ 
then there is a single tabloid:
$$
\begin{pspicture}(0,0)(14,1)
\rput(4.8,0.5){$T=$}
\rput(7,0.5){\BoxedEPSF{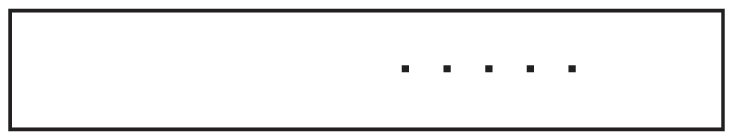 scaled 450}}
\rput(5.65,0.5){$1$}\rput(6.2,0.5){$2$}
\rput(6.75,0.5){$3$}\rput(8.35,0.5){$n$}
\end{pspicture}
$$
and the column group $c_T$ is trivial. There is thus just one vector
$v_T=\{T\}$ in the $1$-dimensional space $M^\lambda$, with
$g\cdot\{T\}=\{T\}$ for all $g\in\Symn$, so that $M^\lambda=S^\lambda$ is the
trivial $1$-dimensional representation of $\Symn$.
\end{example}

\begin{example}
\label{symmetric:example:sign}
At the other extreme we have:
$$
\begin{pspicture}(0,0)(14,3)
\rput(0,0){
\rput(6.15,1.5){$\lambda=$}
\rput(7,1.5){\BoxedEPSF{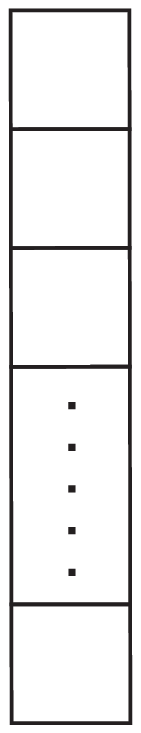 scaled 400}}
}
\end{pspicture}
$$
For any tableau of this shape the column
group $c_T$ is the full symmetric group $\Symn$,
and upto sign, there is just one of the vectors
$$
v_T=\sum_{h\in\Symn}\text{sign}(h)\, h\cdot\{T\}=A-B,
$$
where $A$ is the sum of those terms with $h$ even and $B$ the sum
involving those with $h$ odd. An even permutation $g\in\Symn$
preserves both summands and an odd one swaps them over, so that
$$
g\cdot v_T
=\left\{
\begin{array}{ll}
A-B=v_T,  & g\text{ even},\\
B-A=-v_T,  & g\text{ odd}.\\
\end{array}\right.
=\text{sign}(g)\cdot v_T
$$
The resulting Specht representation $S^\lambda$ is
thus $1$-dimensional (hence irreducible) but
not
the trivial representation; it is called the \emph{sign representation\/} of
$\Symn$. 
\end{example}

\begin{figure}
  \centering
\begin{pspicture}(0,0)(14,3)
\rput(2,2){$\mathlarger{\mathlarger{S}}^{\BoxedEPSF{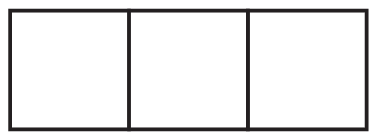
      scaled 80}}$}
\rput(2,1.5){$1$-dimensional}
\rput(2,1.1){trivial representation}
\rput(5,2){$\mathlarger{\mathlarger{S}}^{\BoxedEPSF{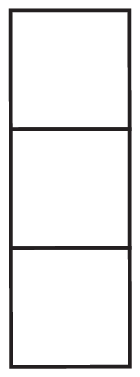
      scaled 80}}$}
\rput(5,1.5){$1$-dimensional}
\rput(5,1.1){sign representation}
\rput(8.5,2){$\mathlarger{\mathlarger{S}}^{\BoxedEPSF{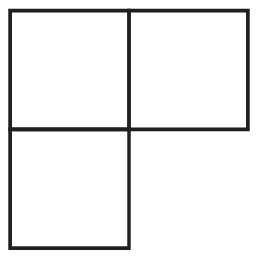
      scaled 80}}$}
\rput(8.5,1.5){$2$-dimensional}
\rput(8.5,1.1){reflectional representation}
\rput(12,1.5){$\BoxedEPSF{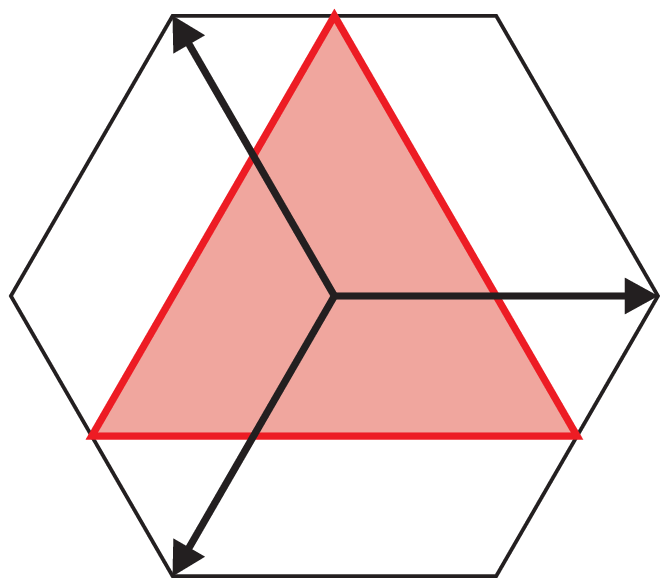 scaled 400}$}
\rput(13.15,0.9){${\red 1}$}\rput(10.85,0.9){${\red
    2}$}\rput(12,2.8){${\red 3}$}
\rput(13.8,1.5){$v_1-v_2$}
\rput(11.1,0.2){$v_2-v_3$}
\rput(11.1,2.8){$v_1-v_3$}
\end{pspicture}
\caption{The three irreducible $\Symthree$ representations over $\C$: the
  trivial representation, the sign representation and the
  reflectional representation. The last corresponds to
the symmetries of an equilateral triangle.} 
  \label{fig:symmetric:three}
\end{figure}

\begin{example}
\label{symmetric:example:reflectional}
If now 
$\lambda=\BoxedEPSF{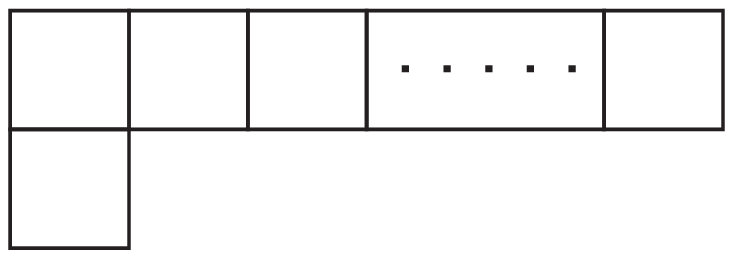 scaled 200}$
then $M^\lambda$ is an $n$-dimensional space with
basis the tabloids:
$$
\begin{pspicture}(0,0)(14,2)
\rput(2,0){
\rput(-1.5,0){
\rput(1.75,1){$v_1=$}
\rput(2.65,1.3){$2$}\rput(3.2,1.3){$3$}\rput(5.35,1.3){$n$}
\rput(2.65,0.75){$1$}
\rput(4,1){\BoxedEPSF{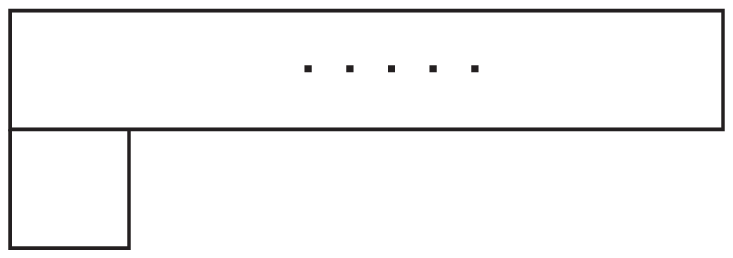 scaled 450}}
}
\rput(5,1){$\cdots$}
\rput(4.25,0){
\rput(1.75,1){$v_n=$}
\rput(2.65,1.3){$2$}\rput(3.2,1.3){$3$}\rput(5.15,1.3){$n-1$}
\rput(2.65,0.75){$n$}
\rput(4,1){\BoxedEPSF{tableau2a.eps scaled 450}}
}}
\rput(6.3,1){,}
\rput(7.5,1,1){,}
\end{pspicture}
$$
and the $\Symn$-action is $g\cdot v_i=v_{g\cdot i}$; it is
thus the permuting coordinates representation of Example
\ref{example:mapping:representation}. 
If
$$
\begin{pspicture}(0,0)(14,1.75)
\rput(-1.8,-0.125){
\rput(1.75,1){$T=$}
\rput(2.65,1.3){$i$}\rput(2.65,0.75){$j$}
\rput(4,1){\BoxedEPSF{tableau1.eps scaled 450}}
\rput(6.5,1){then}
\rput(1,0){
\rput(6.8,1){$v_T=$}
\rput(9,1){\BoxedEPSF{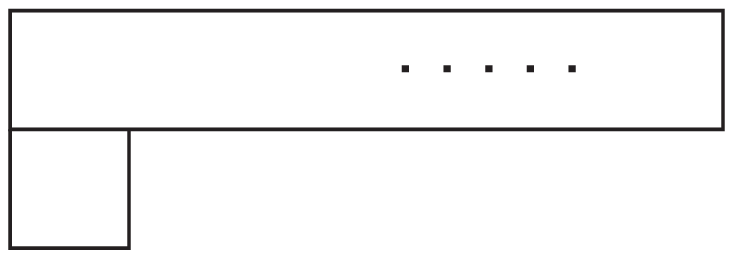 scaled 450}}
\rput(7.65,1.3){$j$}\rput(7.65,0.75){$i$}
\rput(11,1){$-$}
\rput(13,1){\BoxedEPSF{tableau2.eps scaled 450}}
\rput(11.65,1.3){$i$}\rput(11.65,0.75){$j$}
}
}
\end{pspicture}
$$
i.e. $v_T=v_i-v_j\in M^\lambda$, and the $v_T$, as $T$ ranges over the
tableau of shape $\lambda$, give the vectors
$\{v_i-v_j\}_{1\leq i\not=j\leq n}$, which span the hyperplane in
$M^\lambda$ having equation $\sum x_i=0$ (and as promised, the
$v_T$, as $T$ ranges over all tableau, form a dependent set). The restriction of the 
$\Symn$-action on $M^\lambda$ to this hyperplane then gives that $S^\lambda$
the reflectional 
representation of Example \ref{example:permutation:representation}. 
\end{example}


When $n=3$ the possible $\lambda$ are 
$\begin{pspicture}(0,0)(0.15,0.1)
\BoxedEPSF{tableau5a.eps scaled 150}
\end{pspicture}$,
$\begin{pspicture}(0,0)(0.15,0.1)
\BoxedEPSF{tableau6a.eps scaled 150}
\end{pspicture}$ and 
$\begin{pspicture}(0,0)(0.15,0.1)
\BoxedEPSF{tableau1a.eps scaled 150}
\end{pspicture}$ and the $S^\lambda$ are the
three examples in Figure \ref{fig:symmetric:three}.

\begin{vexercise}
\label{exercise:exterior:powers;reflectional}
Show that the exterior power $\bigwedge^p S^{\BoxedEPSF{tableau1.eps scaled 100}}\cong
S^{\BoxedEPSF{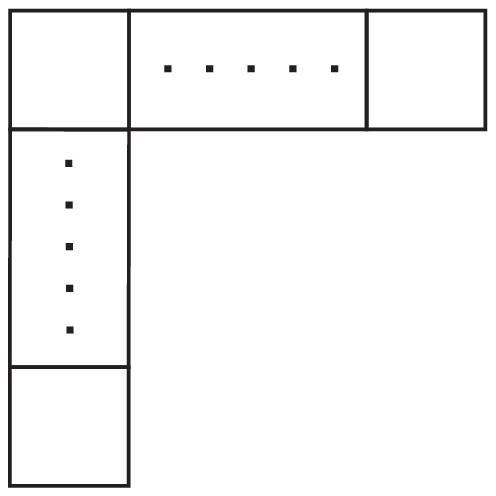 scaled 100}}$, where there are $n-p$ boxes
in the first row of the Young 
diagram of the second Specht representation. (The exterior powers of the reflectional
representation are thus also irreducible).
\end{vexercise}


\section{Reduction}
\label{section:reduction}


This section and the next give two constructions for shuttling back
and forth between representations of a finite regular monoid $S$ and
representations of the maximal subgroups $G_e$ of $S$. The first of
these -- reduction -- squashes $S$-representations down to
$G_e$-representations; the second, induction, blows up 
$G_e$-representations into $S$-representations. In Section
\ref{section:clifford:munn} we will see that with a little care in the
choice of $e$, these constructions turn out to be inverses of each
other. Throughout, the underlying field $k=\C$.


We start with two examples that illustrate all the key features:

\begin{example}
\label{example:I_n:reduction}
Let $S=I_n$, the symmetric invere monoid, and let $V$ be the partial
reflectional representation of Example
\ref{example:partial:reflection:representation}; we saw in that
example that $V$ is irreducible with basis $\{v_1,\ldots,v_n\}$.

Now let $e$ be an idempotent in the $\JJ$-class $\JJ_m$ in the
strategic picture for $I_n$ of Figure
\ref{fig:In_stategic_picture}. This idempotent is the identity map
$\id_X:X\ra X$ on some subset $X\subseteq [n]$ of size $m$ and is the
identity of the maximal subgroup $G_e$ consisting of all bijections
$X\ra X$, which in turn is a copy of the symmetric group $\Symm$. 

To squash $V$ down to a representation of $G_e\cong\Symm$, we take its
image under $e$: let $eV:=e\cdot V=\{ev\,:\,v\in V\}$. Then,
as $e\cdot v_i\not=0$ exactly when $i\in\dom(e)=X$, 
in which case $e\cdot v_i=v_i$, the
space
$eV$ has basis the $v_i$ for $i\in X$. Define an action of $G_e$ on
$eV$ by:
\begin{equation}
  \label{eq:8}
  g\cdot(ev)=(ge)\cdot v,
\quad\quad
(g\in G_e)
\end{equation}
observing, as $e$ is an identity for $G_e$, that $(ge)\cdot
v=(eg)\cdot v=e\cdot(gv)\in eV$. Indeed, $eV$ with this action is just
the 
permuting coordinates 
representation of $\Symm$ given in Example
\ref{example:permutation:representation}. 
We saw there that this representation is reducible when $m\geq 2$,
irreducible when $m=1$, and if $e\in \JJ_0$ is the zero map then $eV=0$.

Upto isomorphism of representations, $eV$ doesn't depend on the choice
of the idempotent $e$ in $\JJ_m$. For suppose that $f=\id_Y:Y\ra Y$ is
another idempotent in $\JJ_m$, with $Y$ a subset of size $m$,  and $f$
the identity of the maximal subgroup $G_f$. We know from Figure
\ref{fig:maximal_subgroups_isomorphic} that $G_f\cong S_Y\cong
S_X\cong G_e$ via the map $h\mapsto s^*hs$, where $s$ is some
bijection $X\ra Y$ and $s^*$ is its semigroup inverse. Defining a
$G_f$-action on $fV$ as in (\ref{eq:8}) gives a
representation of $G_f$. 

\begin{figure}
  \centering
\begin{pspicture}(0,0)(14,5.5)
\rput(-4.8,0){
\rput(6.4,4){
\rput(0,0){$V=$ partial reflectional}
\rput(0,-0.35){representation of $I_n$}
}
\rput(6.4,2.9){
\rput(-0.2,0){$eV=$ permuting coordinates}
\rput(0,-0.35){representation of $\Symm$}
\rput(0,-0.7){($0\leq m\leq n$)}
}
\rput(-1.5,-1.75){
\rput(10,4.5){\BoxedEPSF{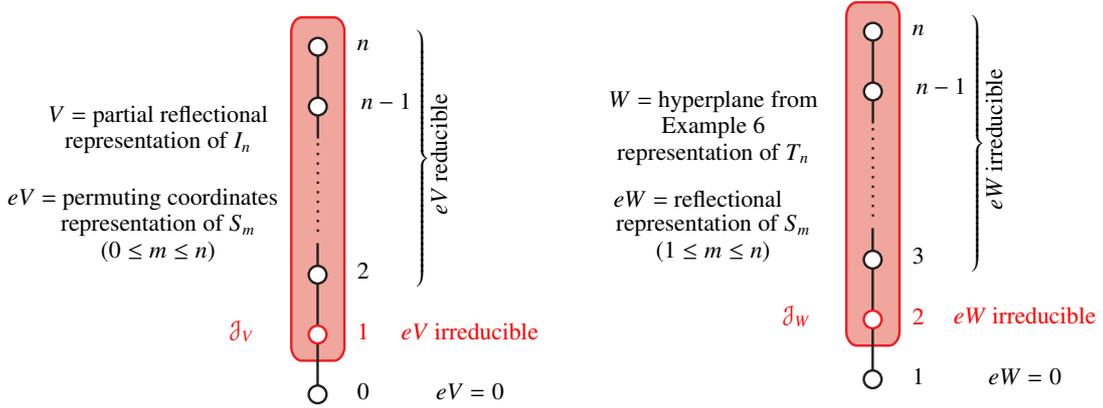 scaled 450}}
\rput(10.6,6.7){$n$}\rput(10.9,5.95){$n-1$}
\rput(10.6,3.7){$2$}\rput(10.6,2.9){${\red 1}$}\rput(10.6,2.1){$0$}
\rput(11.3,5.2){$\left.\begin{array}{c}
\vrule width 0 mm height 34 mm depth 0 pt\end{array}\right\}$}
\rput(12,2.9){{\red $eV$ irreducible}}
\rput(-0.8,0){\rput{90}(12.4,5.25){$eV$ reducible}}
\rput(9,2.9){${\red \JJ_V}$}
\rput(12,2.1){$eV=0$}
}}
\rput(2.5,0.2){
\rput(6.4,4){
\rput(0,0){$W=$ hyperplane from}
\rput(0,-0.35){Example \ref{example:mapping:representation:T_n}} 
\rput(0,-0.7){representation of $T_n$}
}
\rput(6.4,2.7){
\rput(-0.2,0){$eW=$ reflectional}
\rput(0,-0.35){representation of $\Symm$}
\rput(0,-0.7){($1\leq m\leq n$)}
}
\rput(-1.5,-1.75){
\rput(10,4.5){\BoxedEPSF{fig15.eps scaled 450}}
\rput(10.6,6.7){$n$}\rput(10.9,5.95){$n-1$}
\rput(10.6,3.7){$3$}\rput(10.6,2.9){${\red 2}$}\rput(10.6,2.1){$1$}
\rput(11.3,5.2){$\left.\begin{array}{c}
\vrule width 0 mm height 34 mm depth 0 pt\end{array}\right\}$}
\rput(12,2.9){{\red $eW$ irreducible}}
\rput(-0.8,0){\rput{90}(12.4,5.25){$eW$ irreducible}}
\rput(9,2.9){${\red \JJ_W}$}
\rput(12,2.1){$eW=0$}
}}
\end{pspicture}
\caption{Reducing irreducible representations of $I_n$ 
  \emph{(left)} and $T_n$
  \emph{(right)}. The apexes $\JJ_V$, in red, are at the bottom of the red intervals.}
  \label{fig:I_n:reduction}
\end{figure}

The spaces $fV$ and $eV$ are incidentally isomorphic as they both have
dimension $m$; but the map $f\cdot v\mapsto (s^*f)\cdot v$ naturally
gives an
isomorphism $fV\ra eV$ that respects the actions of $G_f$ and $G_e$:
firstly $s^*f=es^*$, so that $(s^*f)\cdot v=(es^*)\cdot
v=e\cdot(s^*v)\in eV$. Then
$$
\begin{pspicture}(0,0)(14,2.5)
\rput(5,0.25){
\rput(1,2){$fV$}
\psline[linewidth=0.5pt]{->}(1.3,2)(2.7,2)
\rput(2,2.2){${\scriptstyle h(-)}$}
\rput(3,2){$fV$}
\psline[linewidth=0.5pt]{->}(1,1.7)(1,0.3)
\psline[linewidth=0.5pt]{->}(3,1.7)(3,0.3)
\rput(0.8,1){${\scriptstyle\cong}$}\rput(3.2,1){${\scriptstyle\cong}$}
\rput(1,0){$eV$}
\psline[linewidth=0.5pt]{->}(1.3,0)(2.7,0)
\rput(2,0.2){${\scriptstyle s^*hs(-)}$}
\rput(3,0){$eV$}
}
\end{pspicture}
$$
commutes, and so the representations $fV$ and $eV$ are indeed
isomorphic as claimed.
The results of the example are summarised on the left of Figure
\ref{fig:I_n:reduction}. 
\end{example}

\begin{example}
\label{example:T_n:reduction}
 The calculations, if not necessarily the results, are similar for the
 mapping representation of $T_n$ from Example
 \ref{example:mapping:representation:T_n}. Now $V$ is reducible, so we
 start instead with the hyperplane $W$ consisting of the $w=\sum \lambda_i
v_i$ with $\sum\lambda_i=0$; this is an irreducible representation of $T_n$.

Let $e$ be the idempotent in $\JJ_m$ (the maps $[n]\ra [n]$ having
image size $m$) given in Figure \ref{fig:T_n:reduction:idempotent}. 
\begin{figure}[h]
  \centering
\begin{pspicture}(0,0)(14,3)
\rput(7,1.5){\BoxedEPSF{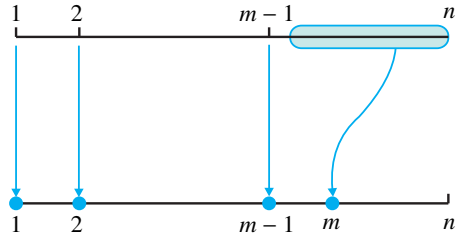 scaled 450}}
\rput(4.2,2.9){$1$}\rput(5,2.9){$2$}\rput(7.5,2.9){$m-1$}\rput(9.9,2.9){$n$}
\rput(4.2,0.1){$1$}\rput(5,0.1){$2$}\rput(7.5,0.1){$m-1$}
\rput(8.35,0.1){$m$}\rput(9.9,0.1){$n$}
\end{pspicture}
\caption{$e\in\JJ_m\subset T_n$ for Example \ref{example:T_n:reduction}.}
  \label{fig:T_n:reduction:idempotent}
\end{figure}
Then the maximal subgroup $G_e$ consists of all bijections from
the fibres of $e$ to the image of $e$ -- again, isomorphic to the
symmetric group $\Symm$. The space $eV$ has basis $\{v_1,\ldots,v_m\}$
and the subspace $eW$ is the hyperplane in $eV$ whose coordinates add to $0$
with respect to this basis. Via the action (\ref{eq:8}) the space
$eV$ is again the permuting coordinates representation of $\Symm$ and
$eW$ is now
the reflectional representation -- see Figure \ref{fig:I_n:reduction} (right).
\end{example}

\paragraph{The general picture is as follows:} let $S$ be a finite regular monoid
and $V$ an \emph{irreducible\/} representation of $S$. Choose an
idempotent in each $\JJ$-class of $S$; the spaces $eV$, equipped with the action
(\ref{eq:8}), are then representations of the maximal subgroups $G_e$
for the various choices of $e$.
It doesn't matter, upto
isomorphism of the resulting representations, which
idempotent in a given $\JJ$-class is chosen. 

\begin{figure}
  \centering
\begin{pspicture}(0,0)(14,5.5)
\rput(0,0){
\rput(2.5,2.75){\BoxedEPSF{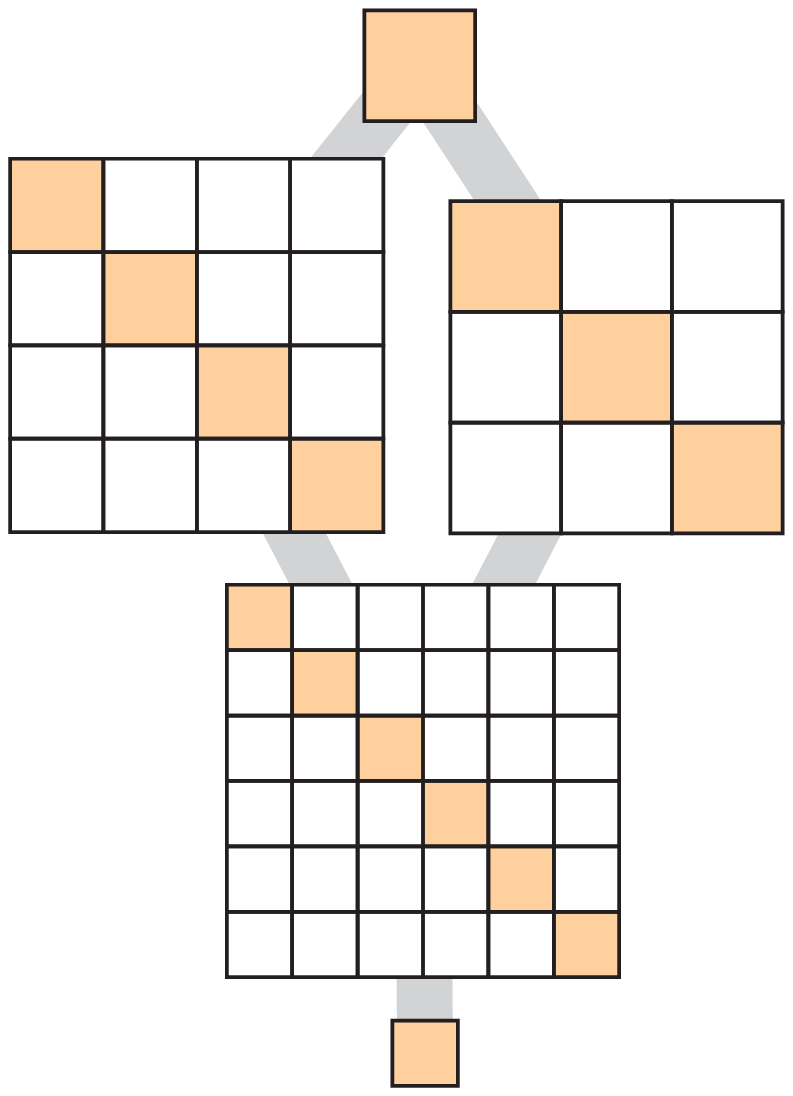 scaled 450}}
\rput(2.6,4.975){${\scriptstyle e_1}$}
\rput(1.375,3.9){${\scriptstyle e_2}$}
\rput(4,3.05){${\scriptstyle e_3}$}
\rput(2.775,1.55){${\scriptstyle e_4}$}
\rput(2.65,0.425){${\scriptstyle e_5}$}
}
\rput(0,0){
\rput(7,2.75){\BoxedEPSF{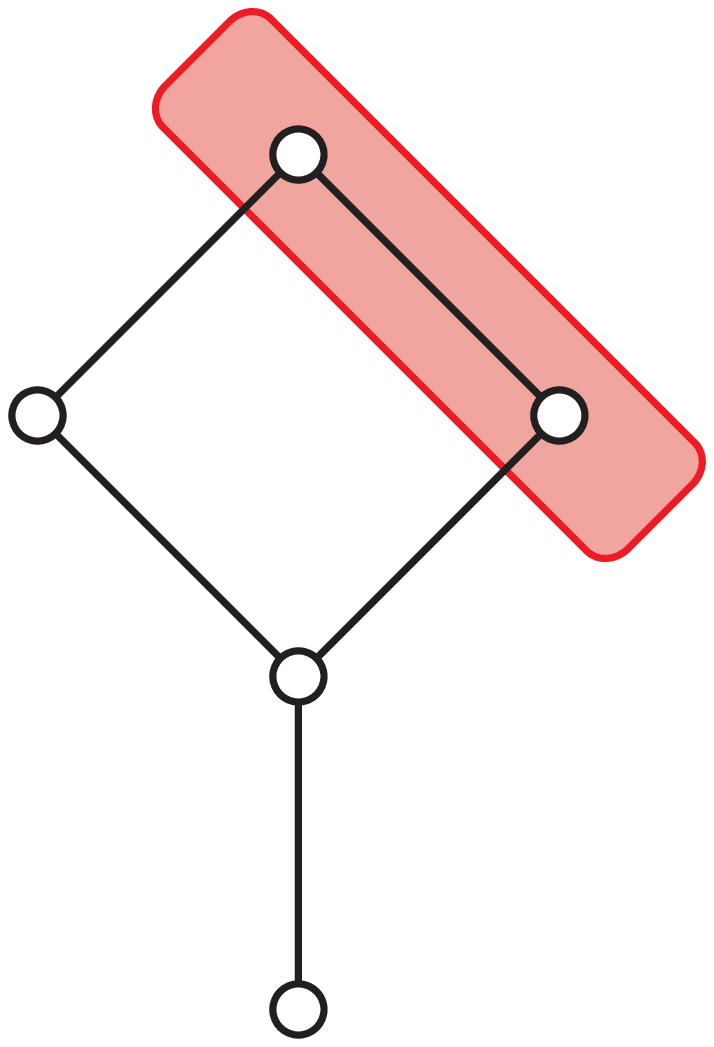 scaled 450}}
\rput(6.7,4.7){${e_1V\not=0}$}
\rput(6.3,3.25){${e_2V=0}$}
\rput(8.65,3.25){${\red e_3V\not=0}$}
\rput(7.45,2){${e_4V=0}$}
\rput(6.7,0.25){${e_5V=0}$}
}
\rput(0,0){
\rput(11.5,2.75){\BoxedEPSF{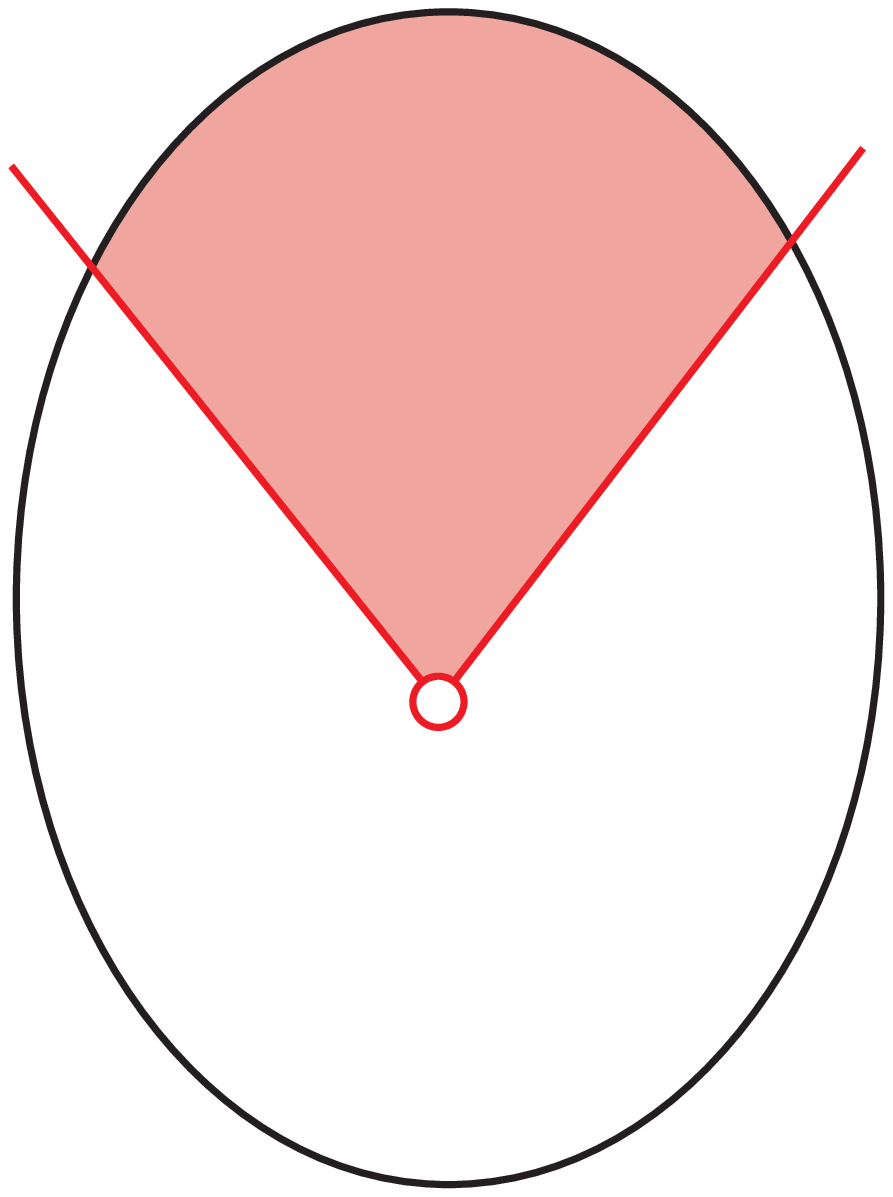 scaled 450}}
\rput(11.5,0.7){$\JJ$-class poset}
\rput(11.5,1.95){${\red \JJ_V}$}\rput(11.5,1.6){{\red apex}}
\rput*(12.675,2.3){{\red $eV$ irreducible}}
\rput(10.3,2.5){$eV=0$}\rput(11.5,4){$eV\not=0$}
\rput(13,5.5){{\red interval}}\rput(13,5.2){${\red \geq \JJ_V}$}
}
\end{pspicture}
\caption{Schematic of reduction: the example of Figure 
\ref{fig:S(G,L)}
  \emph{(left and middle\/)} and the generic set-up \emph{(right)\/.}
}
  \label{fig:reduction:general}
\end{figure}

In particular, whether $eV=0$, or not, is a property of the $\JJ$-class
containing $e$. The $\JJ$-classes for which
$eV\not=0$ form an \emph{interval\/} in the poset of $\JJ$-classes:
there is a 
$\JJ$-class $\JJ_V$ such that  
$$
eV\not=0\text{ exactly when }e\in\JJ_s\text{ with }\JJ_s\geq\JJ_V.
$$
The $\JJ$-class $\JJ_V$ is called the \emph{apex\/} of the
representation $V$, although ``trough'' would probably be a better
name. 
Figure \ref{fig:reduction:general} shows the idea for the $\JJ$-class
poset of Figure \ref{fig:S(G,L)} (left and middle) and generically (right).

\paragraph{} In Example \ref{example:I_n:reduction}, the partial reflectional
representation has apex the $\JJ$-class $\JJ_1$ consisting of the
partial bijections on sets of size $1$; the interval of $\JJ$-classes
$\geq\JJ_1$ is marked on the left in Figure
\ref{fig:I_n:reduction}. For the irreducible representation $W$ of $T_n$ in
Example \ref{example:T_n:reduction}, the apex is the $\JJ$-class
$\JJ_2$ of maps with image size $2$ -- see the right of Figure
\ref{fig:I_n:reduction}. 

In both examples, starting with an irreducible $S$-representation $V$,
we get an irreducible $G_e$-representation \emph{when $e$ is in the
  apex $\JJ$-class $\JJ_V$\/}. 
If $f$ is an idempotent lying in a $\JJ$-class strictly greater
than $\JJ_V$, then the resulting $G_f$-representation may or may not
be irreducible. 

\begin{definitionnumberless}[reduced representations]
\label{definition:reduced:representation}
Let $V$ be an irreducible representation of $S$ with apex $\JJ_V$ and
$e\in\JJ_V$ an idempotent. Then the reduced $G_e$-representation is
given by
\begin{equation}
  \label{eq:9}
  V\downarrow G_e:=eV
\end{equation}
together with the $G_e$-action (\ref{eq:8}).
\end{definitionnumberless}

\begin{vexercise}
\label{exercise:reduction:I_n}
We can verify the general picture for $I_n$ and then for
inverse monoids of the form $S=\sgl$. The reader 
might want to leave this exercise until after they have read 
\S\ref{section:clifford:munn}. 
  \begin{enumerate}
  \item Let $S=I_n$ and let $V$ be a representation of $S$. Let 
$I$ be the set of $\JJ$-classes $\JJ_e$ such that $eV\not=0$. 
Show that $I\not=\varnothing$ and if $\JJ_e\in
I$ and $\JJ_e\leq \JJ_f$ in the $\JJ$-class poset, then
$\JJ_f\in I$. Hence $I$ forms a (closed) interval in the $\JJ$-class
poset. Denote the minimum element by $\JJ_v$.
  \item Let $e$ be an idempotent in this minimal $\JJ$-class $\JJ_V$
    and let $T=\{s_i\}$ be a collection of representatives for the
    $\HH$-classes in the $\LL$-class $\LL_e$. Suppose that we have 
$0\not=eW\not=eV$ with $eW$ a $G_e\cong\Symm$-subrepresentation of
$eV$, and consider the subspace
\begin{equation}
  \label{eq:22}
U=\sum_{s_i\in T} s_i e W  
\end{equation}
of $V$. Show that $U$ is an $I_n$-subrepresentation of $V$. 
  \item If now $V$ is an irreducible $I_n$-representation, then use
    the arguments of \S\ref{section:clifford:munn} 
(about the composition $\irr_m(I_n)\rightarrow\irr(\Symm)\rightarrow\irr_m(I_n)$)
to show that $V=\bigoplus_{s_i\in
      T}s_ieV$. Deduce that $0\not=U\not=V$ for the $U$ of
    (\ref{eq:22}) and hence that $eV$ is 
    irreducible as a $\Symm$-representation.
  \item Repeat the whole thing for an inverse monoid of the form
    $S=\sgl$. 
  \end{enumerate}
\end{vexercise}


\section{Induction}
\label{section:induction}

Induction is the opposite of reduction: it takes as
input a representation of a maximal subgroup and spits out a
representation of the whole semigroup. 

We start with the general construction when $k=\C$. Let $e$ be an idempotent in
the semigroup 
$S$, with $G_e$ the maximal subgroup having identity $e$, and let $V$
be a representation of the group $G_e$. 

The induction of $V$ to an $S$-representation is controlled by the
$\HH$-classes that are in the $\LL$-class $\LL_e$ containing
$G_e$. Choose a transversal for these $\HH$-classes, i.e. a set
$T=\{s_i\}$ with exactly one $s_i$ in each $\HH$-class of $\LL_e$; choose
$e$ itself as the representative in the $\HH$-class $G_e$. The transversal is just
scaffolding for the construction -- the 
resulting $S$-representation is independent of the choice of $T$. 

For each $s_i\in T$, let $V_i$ be an isomorphic copy of the space $V$
having the
elements 
$$
V_i=\{s_i\otimes v : v\in V\}
$$
and vector space operations given by $\lambda(s_i\otimes v)+\mu (s_i\otimes
u)=s_i\otimes (\lambda 
v+\mu u)$. The ``$s_i\,\otimes$" notation is just a device to tell us
which particular copy of $V$ we are working in; other than that it serves no 
purpose and just comes along for the ride in the vector space
operation on $V_i$ (although see the comments in the Note and References section)

\begin{figure}
  \centering
\begin{pspicture}(0,0)(14,8.5)
\rput(0,0){
\rput(7,4.25){\BoxedEPSF{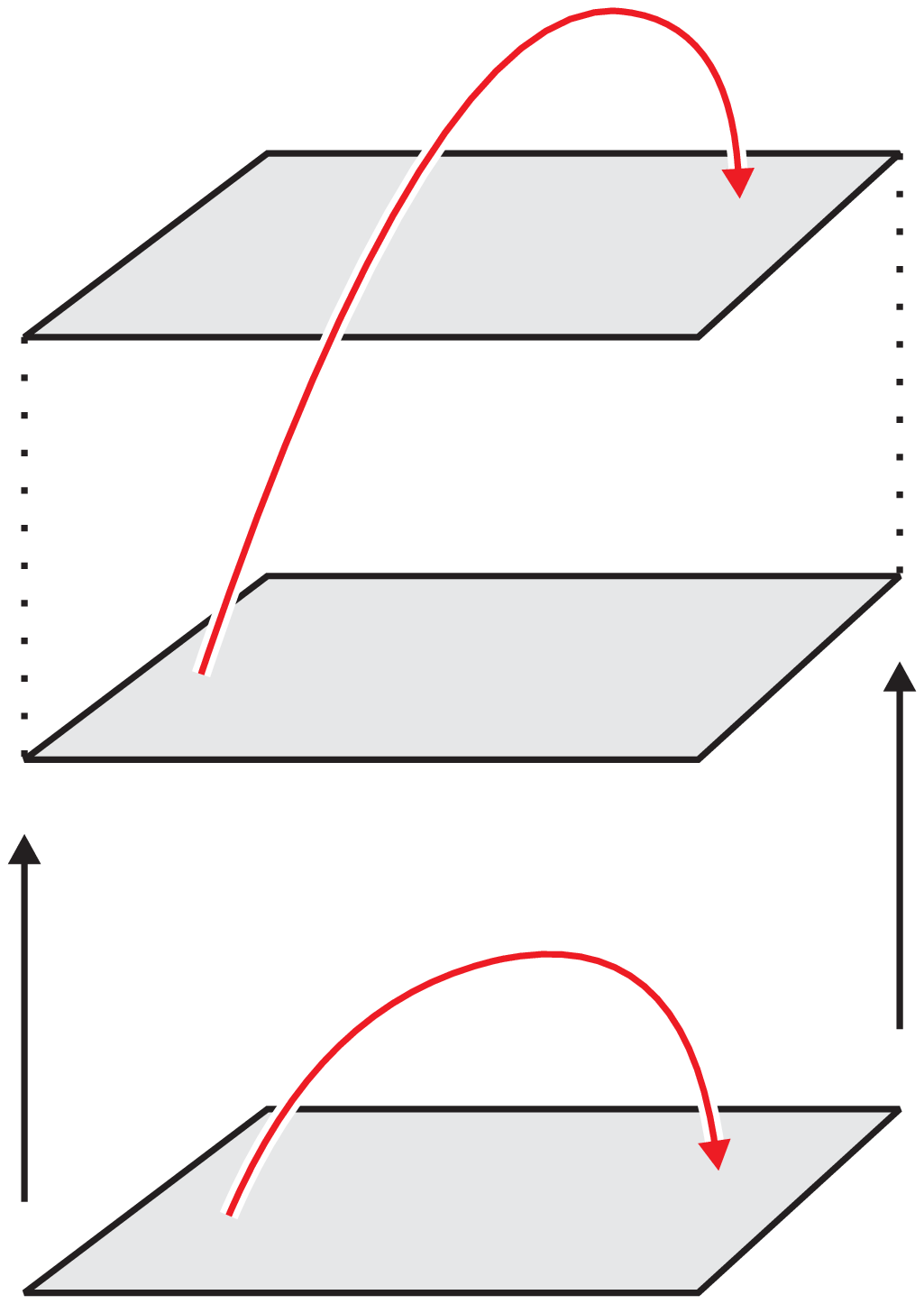 scaled 550}}
\rput(11,1){$G_e$-representation}
\rput(11.5,6){$S$-representation}
\rput(11.5,5.6){$U=\bigoplus_T V_i$}
\rput(10.2,5.75){$\left.\begin{array}{c}
\vrule width 0 mm height 45 mm depth 0 pt\end{array}\right\}$}
}
\rput(5.55,0.6){${\blue v}$}
\rput(8.6,0.9){${\blue g\cdot v}$}
\rput(5.45,3.9){${\blue s_i\otimes v}$}
\rput(5.45,6.5){${\blue s_j\otimes v}$}
\rput(8.45,4.3){${\blue s_i\otimes g\cdot v}$}
\rput(8.45,6.9){${\blue s_j\otimes g\cdot v}$}
\rput(7.2,2.1){${\red g}$}\rput(7.8,8){${\red t}$}
\rput(4,0.2){$V$}\rput(4,3.6){$V_i$}\rput(4,6.2){$V_j$}
\end{pspicture}
\caption{Schematic of the first step of induction.}
  \label{fig:induction}
\end{figure}

Let $U$ be the space 
\begin{equation}
  \label{eq:5}
U=\bigoplus_{s_i\in T} V_i  
\end{equation}
and define an $S$-action on $U$ by 
\begin{equation}
  \label{eq:4}
  t\cdot(s_i\otimes v)
=
\left\{
\begin{array}{ll}
s_j\otimes g\cdot v,  &\text{ if }ts_i\in\LL_e,\text{ hence }ts_i=s_jg\\
0, &\text{ if }ts_i\not\in \LL_e,
\end{array}
\right.
\end{equation}
for $t\in S$.
The action of $t$ thus kills the vector $s_i\otimes v$, unless $ts_i$
is also in the $\LL$-class $\LL_e$, in which case by 
(\ref{eq:1}) there are unique $s_j\in T$ and $g\in G_e$ with
$ts_i=s_jg$. The vector $v$ is then moved in $V$ by the action of $g$,
with the resulting image 
transferred to the corresponding element of $V_j$ -- see Figure
\ref{fig:induction}.

\begin{example}[trivial $S_{\kern-.4mm 1}$ to trivial $I_n$]
\label{induction:example:trivial:to:trivial}
Let $S=I_n$ and $e$ be the zero map $0:\varnothing\ra\varnothing$. The
subgroup $G_e$ is the trivial group (or $S_{\kern-.4mm 1}$!) with the single
element $0$. If $V$ is the trivial representation of $S_{\kern-.4mm
  1}$ then (slightly confusingly)
$0\cdot v=v$ for all $v\in V$, and the transversal $T$ consists of the
one element $\{0\}$, so $U=V_0$ is a single copy of $V$. Finally, for any
$t\in I_n$ we have $t\cdot(0\otimes v)=0\otimes 0\cdot v=0\otimes v$,
and so $U$ is the trivial $1$-dimensional representation of $I_n$.  
\end{example}

\begin{example}
At the opposite end of the strategic picture for $I_n$ we have the
idempotent $e=\id:[n]\ra[n]$, the identity of the group of units
$\Symn$. If $V$ is any representation of $S_n$, then again we have a
transversal $T$
containing the single element $\{e\}$ and so $U=V_e$. Moreover $t\cdot U=0$
unless $t\in\Symn$ is also a unit, in which case it has effect
that of the representation $V$. 

Every $\Symn$-representation is thus
also an $I_n$-representation, just by making that part of $I_n$ not in
$\Symn$ (i.e. the stacked $\JJ$-classes $0,\JJ_1,\ldots,\JJ_{n-1}$ in
the middle of Figure \ref{fig:In_stategic_picture}) 
act as the zero map. It is easy to check that in general any
$G$-representation is also an $S$-representation in this way, when
$G$ is the group 
of units of $S$. 
\end{example}

\begin{example}[trivial $S_{\kern-.4mm 1}$ to partial reflectional
  $I_n$]
\label{example:induction100}
Moving up one rung from the bottom in the strategic picture for $I_n$
in Figure \ref{fig:In_stategic_picture},
let $e$ be the partial identity $1\mapsto 1$ with domain and image
$\{1\}$, so that $G_e$ is again the trivial group $\{e\}$. Also again, let
$V$ be the trivial $1$-dimensional representation of $G_e$, with basis
vector $v$, and action $e\cdot v=v$.  

\begin{figure}
  \centering
\begin{pspicture}(0,0)(14,5.5)
\rput(0,0){
\rput(5.5,2.75){\BoxedEPSF{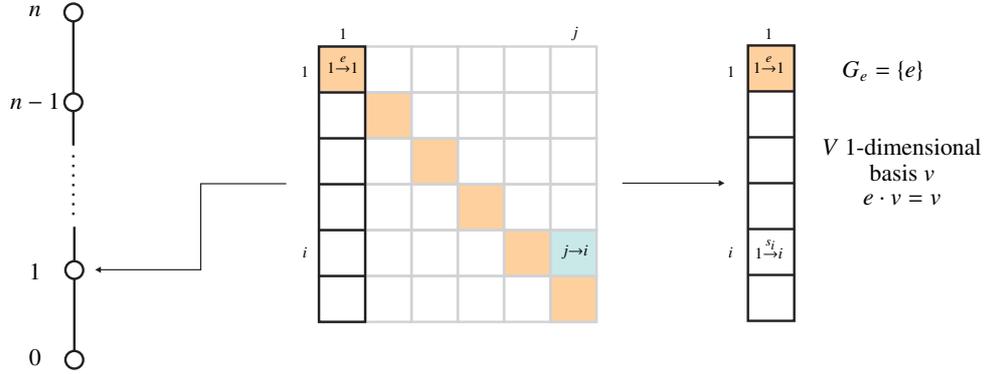 scaled 450}}
\rput(1.6,0.2){$0$}\rput(1.6,1.35){$1$}
\rput(1.6,3.6){$n-1$}\rput(1.6,4.8){$n$}
\rput(5.65,4.5){${\scriptstyle 1}$}
\rput(5.15,1.6){${\scriptstyle i}$}
\rput(5.15,4){${\scriptstyle 1}$}
\rput(8.7,4.5){${\scriptstyle j}$}
\rput(5.65,4.1){${\scriptstyle 1\stackrel{e}{\ra}1}$}
\rput(8.7,1.6){${\scriptstyle j\ra i}$}
\rput(11.25,4.1){${\scriptstyle 1\stackrel{e}{\ra}1}$}
\rput(11.25,1.65){${\scriptstyle 1\stackrel{s_i}{\ra}i}$}
\rput(5.6,0){
\rput(5.65,4.5){${\scriptstyle 1}$}
\rput(5.15,1.6){${\scriptstyle i}$}
\rput(5.15,4){${\scriptstyle 1}$}
}
\rput(12.75,4){$G_e=\{e\}$}
\rput(13,3){$V$ $1$-dimensional}
\rput(13,2.65){basis $v$}
\rput(13,2.3){$e\cdot v=v$}
}
\end{pspicture}
\caption{The trivial representation of the subgroup $G_e$, where
  $e:1\mapsto 1$, induces up to the partial reflectional representation
of $I_n$.}
  \label{fig:induction:trivial:to:partial:permutation}
\end{figure}

The $\JJ$-class containing the subgroup $G_e$ consists of all the bijections
with domain and image of size $1$ (Figure
\ref{fig:induction:trivial:to:partial:permutation}) and $\LL_e$ is the
column of all the maps with domain $1$. There is no choice for the
representatives $T$: they are the partial bijections $s_i:1\mapsto i$. 
The copy $V_i$ of $V$ has basis the vector $s_i\otimes v$, and so the space
$U$  of (\ref{eq:5}) is
$n$-dimensional with basis $\{s_i\otimes v\}_{1\leq i\leq n}$. 

For the $I_n$-action, we have $ts_i\in\LL_e$ when it has domain $1$, and this
is exactly when $i$ lies in the domain of $t$, in which case
$ts_i=s_{t(i)}=s_{t(i)}e$. Thus
$$
  t\cdot(s_i\otimes v)
=
\left\{
\begin{array}{ll}
s_{t(i)}\otimes v,  &\text{ if }i\in\dom(t),\\
0, &\text{ else.}
\end{array}
\right.
$$
Replacing $s_i\otimes v$ by $v_i$ we get the formula (\ref{eq:2}), and
so $U$ is the partial reflectional representation of $I_n$. 
\end{example}

In Examples
\ref{induction:example:trivial:to:trivial}-\ref{example:induction100} an
irreducible $G_e$-representation $V$ becomes an irreducible
$S$-representation $U$. We are not always so lucky:

\begin{example}[trivial $S_{\kern-.4mm 1}$ to mapping $T_n$]
\label{induction:example200}
This is very similar to Example \ref{example:induction100}. 
Let $S$ be the full transformation monoid $T_n$, with strategic
picture on the left of Figure
\ref{fig:induction:trivial:to:nonirreducible:Tn}. 

\begin{figure}
  \centering
\begin{pspicture}(0,0)(14,5.5)
\rput(0,0.25){
\rput(5.5,2.75){\BoxedEPSF{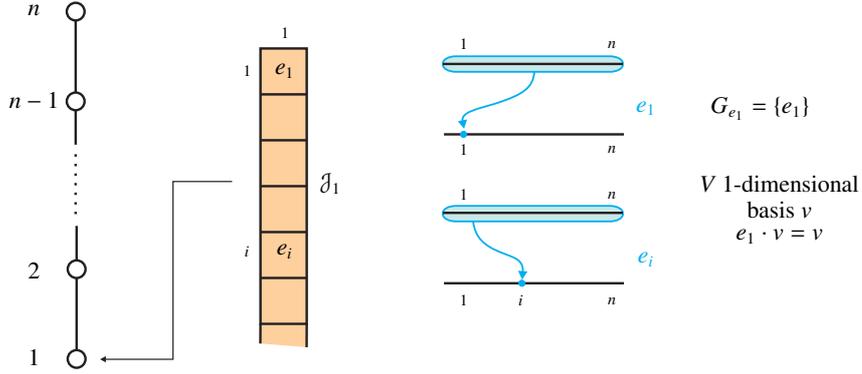 scaled 450}}
\rput(2.7,0.2){$1$}\rput(2.7,1.35){$2$}
\rput(2.7,3.6){$n-1$}\rput(2.7,4.8){$n$}
\rput(0.35,0){
\rput(5.65,4.5){${\scriptstyle 1}$}
\rput(5.15,1.6){${\scriptstyle i}$}
\rput(5.15,4){${\scriptstyle 1}$}
}
\rput(6,4){$e_1$}\rput(6,1.6){$e_i$}
\rput(6.6,2.5){$\JJ_1$}
\rput(8.35,4.35){${\scriptstyle 1}$}
\rput(10.3,4.35){${\scriptstyle n}$}
\rput(8.35,2.95){${\scriptstyle 1}$}
\rput(10.3,2.95){${\scriptstyle n}$}
\rput(8.35,2.35){${\scriptstyle 1}$}
\rput(10.3,2.35){${\scriptstyle n}$}
\rput(8.35,0.95){${\scriptstyle 1}$}
\rput(9.1,0.95){${\scriptstyle i}$}
\rput(10.3,0.95){${\scriptstyle n}$}
\rput(10.75,3.5){${\cyan e_1}$}
\rput(10.75,1.5){${\cyan e_i}$}
\rput(-0.5,-0.5){
\rput(12.75,4){$G_{e_1}=\{e_1\}$}
\rput(13,3){$V$ $1$-dimensional}
\rput(13,2.65){basis $v$}
\rput(13,2.3){$e_1\cdot v=v$}
}
}
\end{pspicture}
\caption{An irreducible $V$ doesn't necessarily give an irreducible
  $U$: the trivial representation of the subgroup $G_{e_1}$, where
  $e_1$ is the constant map $[n]\mapsto 1$, induces up to the mapping
  representation  
of $T_n$.}
  \label{fig:induction:trivial:to:nonirreducible:Tn}
\end{figure}

The $\JJ$-class at the bottom consists of the maps with image size $1$
-- the constant maps. There is a single $\LL$-class and $n$
$\RR$-classes, each containing the single constant map $e_i:[n]\mapsto
i$ for $1\leq i\leq n$. These are all idempotents, so that every
$\HH$-class in this $\JJ$-class is a maximal subgroup. (You can see this in Figure
\ref{fig:Tn_Hclasses_idempotents}, where every box of $\JJ_1$ is a maximal
subgroup). In anycase, there is no choice once again for $T$, which
must be the $\{e_i\}_{1\leq i\leq n}$. 

Let $e_1$ be the constant map $e:[n]\mapsto 1$ and $V$ be the
trivial $G_{e_1}$-representation with basis the vector $v$. For each $i$,
the space $V_i$ has basis the vector $e_i\otimes v$ and the $U$ of
(\ref{eq:5}) is $n$-dimensional with basis $\{e_i\otimes v\}_{1\leq
  i\leq n}$.

For any $t\in T_n$ we have $t\cdot e_i=e_{t(i)}=e_{t(i)}e_1$, so that
(\ref{eq:4}) becomes
$t\cdot(e_i\otimes v)=e_{t(i)}\otimes v$, and once again we have the
formula (\ref{eq:2}), after replacing $e_i\otimes v$ with $v_i$. The
space $U$ thus carries the mapping representation of $T_n$, which is
reducible, even though the seeding representation $V$ of the
subgroup $G_{e_1}\cong S_{\kern-.4mm 1}$ is irreducible. 

Before leaving the example we observe something for
later on: the $\RR$-class $\RR_e$ contains just the single element
$\{e_1\}$ with $e_1\cdot(e_i\otimes v)=e_1\otimes v$ for all $i$. Suppose
that $u\in U$ is a vector that is annihilated by $e_1$, i.e.
$$
u=\sum \lambda_i (e_i\otimes v)\text{ with }e_1\cdot u=0.
$$
Then
$$
e_1\cdot u=0
\Leftrightarrow
\biggl(\sum \lambda_i\biggr)(e_1\otimes v)=0
\Leftrightarrow 
\sum \lambda_i=0,
$$
so that, after replacing $e_i\otimes v$ with $v_i$,
the set of such annihilated vectors is the hyperplane $W$ consisting
of the $w=\sum \lambda_iv_i$ where $\sum\lambda_i=0$. 

These annihilated vectors thus form a subrepresentation of $U$. Even
more is true: $U$ is the $n$-dimensional mapping representation and
$W\subset U$ is $(n-1)$-dimensional, so the quotient representation
$U/W$ is $1$-dimensional, hence irreducible.
In particular, $W$ is a maximal subrepresentation of $U$. 
\end{example}

\paragraph{Returning to generalities,} let $V$ be 
a representation of the maximal
subgroup $G_e$ and $U$ the space given in (\ref{eq:5}). As in the
example just done, consider the vectors in $U$ that are annihilated by
the elements of the $\RR$-class $\RR_e$:
\begin{equation}
  \label{eq:6}
 \text{Ann}_{e}(U):=\{u\in U: s\cdot u=0\text{ for all }s\in\RR_e\}.
\end{equation}

\begin{definitionnumberless}[induced representations]
\label{definition:induced:representation}
Let $V$ be a representation of the maximal subgroup $G_e$ of $S$ and
$U$ be the $S$-representation given by  (\ref{eq:5}) and (\ref{eq:4}). Then
the $S$-representation induced by $V$ is the quotient
\begin{equation}
  \label{eq:7}
  V\uparrow S:=U/\text{Ann}_{e}(U)
\end{equation}
\end{definitionnumberless}

As in Example
  \ref{induction:example200}, if $V$ is irreducible, then
$\text{Ann}_{e}(U)$ is a maximal subrepresentation of
  $U$, and $V\uparrow S$ is irreducible. The construction depends
  only on the $\JJ$-class of $e$: if $e$ and $f$ lie in the same
  $\JJ$-class then the resulting 
  induced representations are isomorphic.

\begin{example}
\label{example:inversesemigroup:annihilator}
We can verify some of these general claims in the setting of $I_n$. In the
Exercise following, we do this for an inverse monoid of
the form $\sgl$. 

Suppose that $|X|=m$ and $V$ is
  a representation of the maximal subgroup $G_e\cong\Symm$, with $e$ the partial
  identity $X\ra X$, and let $U$ be the $I_n$-representation described
  in (\ref{eq:5})-(\ref{eq:4}). 
We show first that the annihilator
  $\text{Ann}_{e}(U)$ is trivial. 

To see this, let $T=\{s_i\}$ be the
  transversal used for the induction and
  $u=\sum_i s_i\otimes v_i\in \text{Ann}_{e}(U)$. 
Fix an $s_j:X\ra Y\in T$ with $s_j^*:Y\ra X$ the semigroup
  inverse of $s_j$ -- see Figure
  \ref{fig:induction:zeroannihilator}. Then for any $i$ we have:
$$
s_j^*s_i\in\LL_e
\Leftrightarrow 
\dom(s_j^*s_i)=X
\Leftrightarrow
\im(s_i)=\dom(s_j^*)
\Leftrightarrow
\im(s_i)=Y
\Leftrightarrow
s_i=s_j.
$$
Thus, on the one hand, by (\ref{eq:4}):
$$
s_j^*\cdot u=s_j^*\cdot(s_j\otimes v_j)=e\otimes v_j,
$$
while on the other, $s_j^*\in\RR_e$ and $u\in \text{Ann}_{e}(U)$ gives
$s_j^*\cdot u=0$. The conclusion 
is that $v_j=0$ and hence $s_j\otimes v_j=0$. Letting $j$ vary we see
that $u=0$ and hence $\text{Ann}_{e}(U)=0$ as claimed. 

\begin{figure}
  \centering
\begin{pspicture}(0,0)(14,5.5)
\rput(3.5,0.25){
\rput(3.5,2.5){\BoxedEPSF{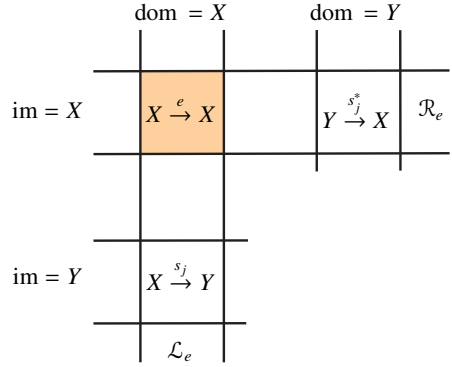 scaled 450}}
\rput(2.3,3.65){$X\stackrel{e}{\ra}X$}\rput(4.6,3.65){$Y\stackrel{s_j^*}{\ra}X$}
\rput(2.3,1.45){$X\stackrel{s_j}{\ra}Y$}
\rput(0.55,3.65){$\im=X$}\rput(0.55,1.45){$\im=Y$}
\rput(2.3,4.95){$\dom=X$}\rput(4.6,4.95){$\dom=Y$}
\rput(5.6,3.65){$\RR_e$}\rput(2.3,0.45){$\LL_e$}
}
\end{pspicture}
\caption{$\text{Ann}_{e}(U)=0$ in $I_n$.}
  \label{fig:induction:zeroannihilator}
\end{figure}

The annihilator is then certainly an $I_n$-subrepresentation of
$U$, albeit for trivial reasons! Suppose now that $V$ is an
\emph{irreducible\/} $\Symm$-representation. We 
claim that the annihilator is a maximal subrepresentation, or
equivalently, that $U$ is an irreducible $I_n$-representation. Let $u$
be a non-zero vector in $U=\bigoplus_T s_i\otimes V$ with
$u=\sum_i s_i\otimes v_i$ and $s_j\otimes v_j\not=0$ for some $j$ (so
that in particular, $v_j\not=0$). Let $W=I_n\cdot u\subset U$ be the
set (hence subspace -- Exercise) of all images of $u$ under the elements of
$I_n$. 

We claim that $W=U$. We have $s_j^*\cdot u\in W$, where as above
$$
s_j^*\cdot u=e\otimes v_j
$$
with $0\not=e\otimes v_j\in V$. On the one hand, we have
$\Symm\cdot(e\otimes v_j)\subseteq W$ (as $\Symm\subset I_n$), while on
the other $\Symm\cdot(e\otimes v_j)$ is a non-zero
$\Symm$-subrepresentation of the irreducible $\Symm$-representation
$V$. The conclusion is that $\Symm\cdot(e\otimes v_j)=V$, and hence
$V\subseteq W$. For any $s_i\in T$ we have by the definition of the action
in (\ref{eq:4}) that 
$$
s_i\cdot(e\otimes v)=s_i\otimes v
$$
and so $s_i\cdot V=s_i\otimes V$. Thus, as soon as we have $V\subseteq
W$ then we have $s_i\otimes V\subseteq W$ for each $s_i$, and thus
$U\subseteq W$. 
This proves the claim, and so the $I_n$-representation $U$ is irreducible. 

To summarise: if $V$ is an irreducible $\Symm$-representation then the
induced representation $V\uparrow I_n$ is the space $U$, and this in
turn is an irreducible $I_n$-representation.
\end{example}

\begin{vexercise}
\label{exercise:sgl:ann}
Let $S$ be an inverse monoid of the form $S=\sgl$ and $V$ an
irreducible representation of the maximal subgroup $G_e$. Let $U$ be the
$S$-representation given by (\ref{eq:5})-(\ref{eq:4}). 
  \begin{enumerate}
   \item Show that $\text{Ann}_{e}(U)=0$ (\emph{hint :\/} prove
     the following fact first: if the finite group $G$ acts on the
     lattice $L$, and if $a,b$ lie in the same $G$-orbit with $a\leq
     b$, then $a=b$; in other words, distinct elements of $L$ in the same
     $G$-orbit are not comparable.)
   \item Mimic the argument above for $I_n$ to show that $U$ is irreducible.
  \end{enumerate}
\end{vexercise}

\paragraph{Here is how induction works for an inverse monoid of the form
$S=S\kern-1pt(G,L)$.}
\label{example:induction:sgl}
Let $e=\id_a$ be an idempotent and $V$ a
representation of the maximal subgroup $G_a$ given in  (\ref{eq:17}). 

The $\HH$-classes in $\LL_e$ are parametrised by the $G$-orbit of $a$,
say $G\cdot a=\{a,b,\ldots\}$; let
$\aa=\id,\bb,\ldots$ be elements of $G$ such that 
$$
\aa:a\mapsto a, \bb:a\mapsto b,\ldots
$$
We then take our transversal $T$ to be
$\aa_a,\bb_a,\ldots$
In light of 
Exercise \ref{exercise:sgl:ann},
the
induced representation is carried by the space
$$
V\uparrow S=\bigoplus_{\bb_a\in T}b\otimes V
$$
where (after simplifying notation a little) $b\otimes V=\{b\otimes
v:v\in V\}$ and, as usual, the vector space 
operations happen in the ``$v$'' coordinate, with the ``$b\,\otimes$"
just a label for the copy of $V$.

Suppose that $s=g_c$ is some element of $S$. To understand the action
in (\ref{eq:4}) we need to compute products like $s\,\bb_a$: as
$s\,\bb_a=g_c\bb_a=(g\bb)_{\bb^{-1}\cdot c\wedge a}$, we have that
$s\,\bb_a$ lies in $\LL_e$ exactly when
$$
\bb^{-1}\cdot c\wedge a=a
\Leftrightarrow
a\leq \bb^{-1}\cdot c
\Leftrightarrow
b\leq c
$$
(recall that $\bb\in G$ sends $a$ to $b$). Moreover, if $b\leq c$, then
$s\,\bb_a$ lies in the $\HH$-class $\HH_d$, where $d=g\cdot b$. The
description (\ref{eq:14}) gives the element
$\delta_d^{-1}s\,\bb_a=(\delta^{-1}g\bb)_a$ of $G_a$ and
$$
s\,\bb_a=\delta_a\cdot (\delta^{-1}g\bb)_a
$$
Thus, for $b\otimes v$ an element of $V\uparrow S$ and for
$g_c$ an element of $S=S\kern-1pt(G,L)$ we have the action:
\begin{equation}
  \label{eq:18}
    g_c\cdot(b\otimes v)
=
\left\{
\begin{array}{ll}
d\otimes h\cdot v,  &\text{ if }b\leq c\\
0, &\text{ if }b\not\leq c,
\end{array}
\right.
\end{equation}
where $d=g\cdot b$ and $h=(\delta^{-1}g\bb)_a$ with $\delta:a\mapsto
d$ one of the elements of $G$ chosen above. 

We can say more. If $s$ lies in the $\JJ$-class $\JJ_c$ with
$\JJ_a\not\leq\JJ_c$, then for any $b$ in the $G$-orbit of $a$ we have
$b\not\leq c$ in $L$, hence $s\cdot(b\otimes v)=0$, and so $s\cdot
V\uparrow S=0$. On the other hand, if $s=\id_c\in\JJ_c$ with
$a\leq c$ and if $a\otimes v\not=0$
in $V\uparrow S$, then $s\cdot(a\otimes v)=a\otimes v$, and so $s\cdot
V\uparrow S\not=0$. The conclusion is that $s\cdot V\uparrow S\not= 0$
precisely for those $s$ lying in the $\JJ$-classes that are
$\geq\JJ_a$ in the $\JJ$-class poset. In particular, the apex of $V\uparrow S$, for
$V$ a representation of $G_a$, is $\JJ_a$.

\begin{example}
\label{example:Tn:induction}
  We return to $T_n$ and Example \ref{induction:example200} where
  $e_1:[n]\mapsto 1$ 
  is our idempotent, $V$ is the trivial
  representation of $G_{e_1}\cong S_{\kern-.4mm 1}$ and $U$ the mapping
  representation of $T_n$. We saw at the end 
of Example \ref{induction:example200}
that $\text{Ann}_{e}(U)$
  is the hyperplane $W$ in $U$ consisting of the $w=\sum \lambda_iv_i$
  with $\sum\lambda_i=0$.  The induced representation $V\uparrow T_n=U/W$
  is thus $1$-dimensional, and as $v_i-v_j\in W$ 
  we have $v_i+W=v_j+W$ for all $i$ and $j$. Taking $v_1+W$ to be the
  basis vector for $V\uparrow T_n$, we have for any $t\in T_n$ that:
$$
t\cdot(v_1+W)=t(v_1)+W=v_{t(1)}+W=v_1+W,
$$
so that $V\uparrow T_n$ is the trivial representation.
\end{example}


\section{The Clifford-Munn correspondence}
\label{section:clifford:munn}

Induction creates irreducible representations of a (finite regular)
monoid out of irreducible representations of its maximal
subgroups. With a little care in the accounting, this process gives a
1-1 correspondence between the irreducible $S$-representations and the
irreducibles of a certain collection of maximal subgroups. 
This bijection 
is called the
\emph{Clifford-Munn correspondence\/}. 

The bijection comes about by showing that reduction
is the inverse of
induction; for us, this is the principal purpose of reduction.
The apex of an $S$-representation $V$
tells us the ``right'' maximal subgroup to reduce to.

\begin{figure}
  \centering
\begin{pspicture}(0,0)(14,6)
\rput(0,0.25){
\rput(0,0.25){
\rput(4.5,1.65){\BoxedEPSF{fig9a.eps scaled 250}}
\rput(4.5,0.25){$\JJ$-classes}
}
\rput(2,3){\BoxedEPSF{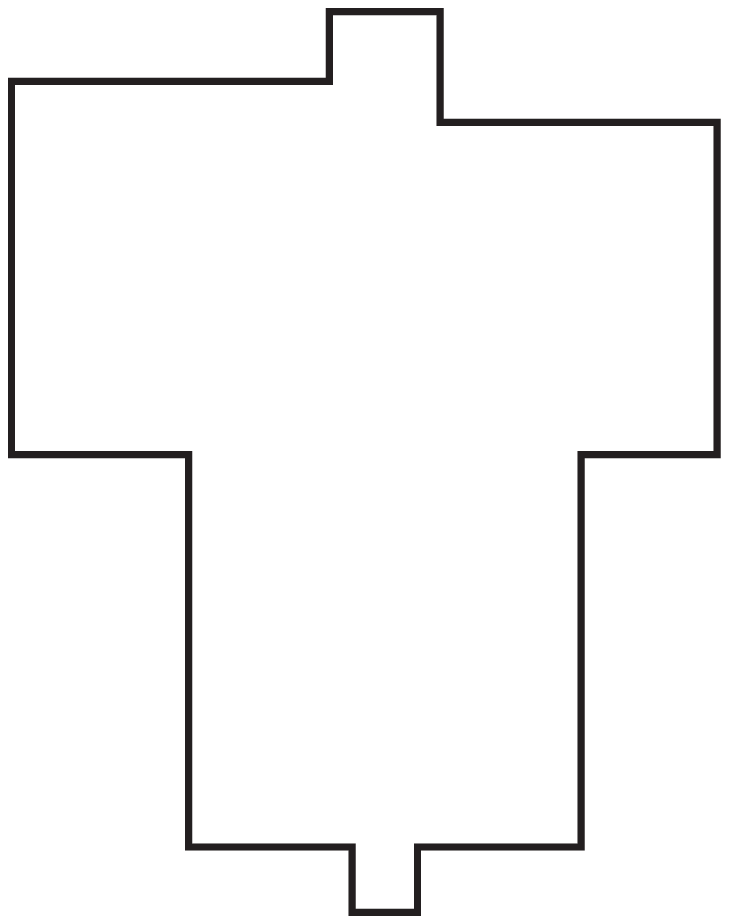 scaled 450}}
\rput(2,3.6){$\irr(S)$}
\psline[linewidth=0.5pt]{<->}(3.8,3.5)(5.2,3.5)
\rput(4.5,3.65){${\scriptstyle\text{partition}}$}
\rput(7,3){\BoxedEPSF{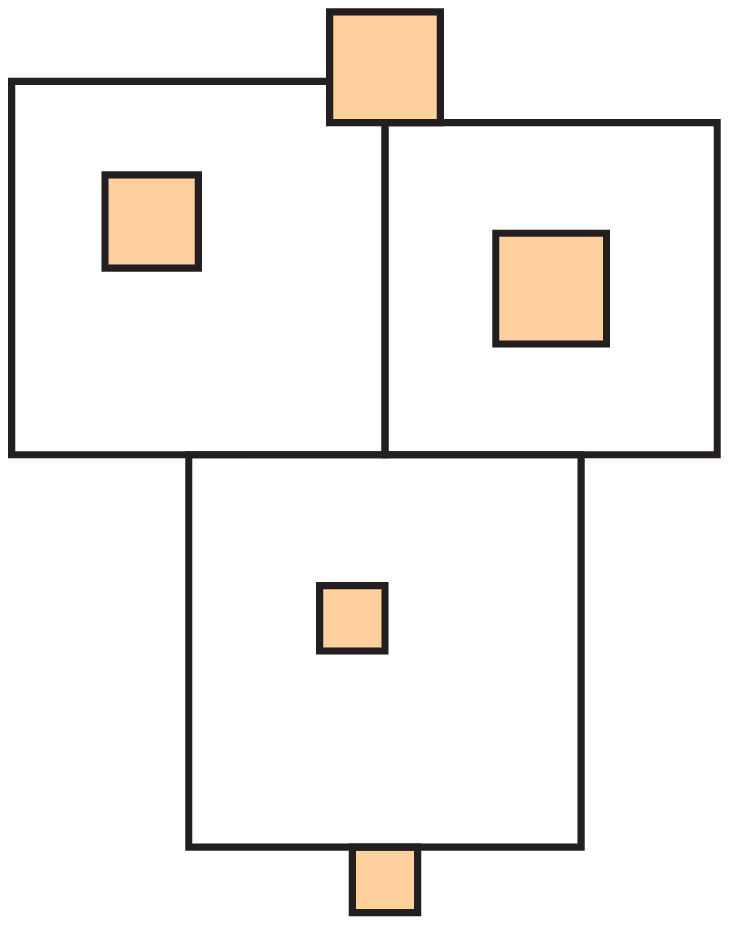 scaled 450}}
\rput(7.1,4.8){${\scriptstyle e_1}$}
\rput(6.05,4.1){${\scriptstyle e_2}$}
\rput(7.85,3.775){${\scriptstyle e_3}$}
\rput(6.965,2.275){${\scriptstyle e_4}$}
\rput(7.1,1.1){${\scriptstyle e_5}$}
\rput(6.2,3.4){${\scriptstyle\irr_{e_2}(S)}$}
\rput(7.8,3.3){${\scriptstyle\irr_{e_3}(S)}$}
\rput(7,1.675){${\scriptstyle\irr_{e_4}(S)}$}
\psline[linewidth=0.5pt]{<->}(8.75,3.5)(10.15,3.5)
\rput(9.45,3.65){${\scriptstyle\text{1-1}}$}
\rput(9.45,3.35){${\scriptstyle\text{correspondence}}$}
\rput(12.1,3){\BoxedEPSF{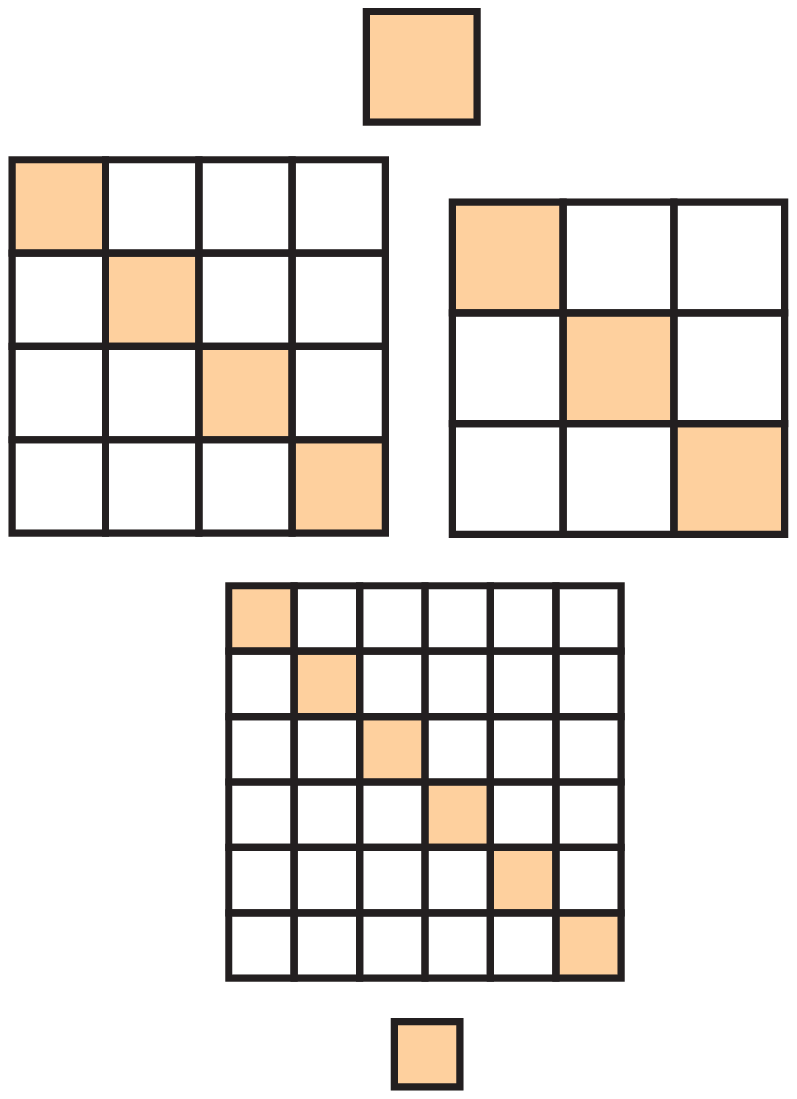 scaled 450}}
\rput(12.2,5.2){${\scriptstyle e_1}$}
\rput(11,4.15){${\scriptstyle e_2}$}
\rput(13.1,3.8){${\scriptstyle e_3}$}
\rput(12.1,2.1){${\scriptstyle e_4}$}
\rput(12.25,0.7){${\scriptstyle e_5}$}
\rput(11,5){${\scriptstyle\irr(G_{e_2})}$}
\rput(13.5,4.8){${\scriptstyle\irr(G_{e_3})}$}
\rput(0.1,-1){\rput{90}(11,3){${\scriptstyle\irr(G_{e_4})}$}}
}
\end{pspicture}
\caption{Schematic of the Clifford-Munn correspondence: the
  irreducibles of $S$ \emph{(left)} partitioned into their various
  apexes \emph{(middle)} which in turn are in 1-1 correspondence with
  the irreducibles of the corresponding maximal subgroups \emph{(right)}.}
  \label{fig:CM:correspondence}
\end{figure}

Figure \ref{fig:CM:correspondence} illustrates the correspondence, where as
usual, the strategic picture of $S$ drives the whole process. Let $\irr(S)$ be
the set of isomorphism classes of irreducible $S$-representations and
$E=\{e_i\}$ be a set of 
idempotents in 1-1 corrrespondence with the $\JJ$-classes of $S$. 
For $e\in E$ let
$$
\irr_e(S)=\{V\in\irr(S):\JJ_V=\JJ_e\},
$$
be the irreducible $S$-representations $V$ whose apex $\JJ_V$ is the
$\JJ$-class $\JJ_e$ containing $e$. Every irreducible $V$ has a
uniquely determined apex -- 
the set of  $S$-irreducibles $\irr(S)$ is thus partitioned into the
$\irr_e(S)$ as $e$ ranges over $E$. 
Finally, let $\irr(G_e)$ be the irreducible representations of the
maximal subgroup $G_e$. 

\begin{theorem}[Clifford-Munn correspondence]
\label{theorem:CM}
For a fixed $e\in E$, the maps:
$$
\begin{pspicture}(0,0)(14,1)
\rput(0,-0.35){
\rput(5.5,0.75){$\irr_e(S)$}
\psline[linewidth=0.5pt]{->}(6.2,0.9)(7.8,0.9)
\rput(7,1.1){${\scriptstyle V\rightarrow V\downarrow G_e}$}
\rput(8.5,0.75){$\irr(G_e)$}
\psline[linewidth=0.5pt]{<-}(6.2,0.6)(7.8,0.6)
\rput(7,0.4){${\scriptstyle V\uparrow S \leftarrow V}$}
}
\end{pspicture}
$$
are mutual inverses, inducing a bijection
$\irr(S)\rightleftarrows\bigcup_{e\in T}\irr(G_e)$.
\end{theorem}

\paragraph{} We will prove the correspondence in the context of the symmetric inverse
monoid $I_n$ when $k=\C$. 
Exercise \ref{exercise:CMP} at the end of the section asks for a 
proof for an inverse monoid of the form $S=\sgl$.

Fix then, in $I_n$, the $\JJ$-class
$\JJ_m$ for some $0\leq m\leq n$ and the idempotent
$e=\id:[m]\ra [m]$. The maximal subgroup
$G_e$ is isomorphic to $\Symm$ and consists of all partial bijections $[m]\ra[m]$. As
we have a nice total order on the $\JJ$-classes, we write
$\irr_m(I_n)$ for $\irr_e(I_n)$. 

\paragraph{The map $\irr_m(I_n)\rightarrow\irr(\Symm)$ given by $V\mapsto
  V\downarrow\Symm$:} as $V\in\irr_m(I_n)$, it is irreducible
with apex $\JJ_m$, and hence $V\downarrow\Symm=eV$ is an irreducible
$\Symm$-representation by 
\S\ref{section:reduction}. Thus $V\downarrow\Symm\in\irr(\Symm)$. 

\paragraph{The map $\irr(\Symm)\rightarrow\irr_m(I_n)$ given by $V\mapsto
  V\uparrow I_n$:} first, we show that this is indeed a map. For
$V\in\irr(\Symm)$, we saw in Example 
\ref{example:inversesemigroup:annihilator} that 
$\text{Ann}_{e}(U)=0$ where $U$ is the $I_n$-representation given in 
(\ref{eq:5})-(\ref{eq:4})
and that $V\uparrow I_n=U$ is
irreducible. Thus $V\uparrow I_n\in\irr(I_n)$; we need it to be in
$\irr_m(I_n)$, i.e. to have apex the $\JJ$-class $\JJ_m$. The
following essentially repeats the more general arguments immediately
preceding Example \ref{example:Tn:induction}, but in a concrete
setting.  

The $\LL$-class $\LL_e$ consists of all the partial bijections with
domain $[m]$. If $Y=\{i_1,\ldots,i_m\}$ is some subset of size $m$,
then let $s_Y:[m]\ra Y$ 
\begin{figure}
  \centering
\begin{pspicture}(0,0)(14,3.5)
\rput(0,0.25){
\rput(-3.5,0){
\rput(7,1.5){\BoxedEPSF{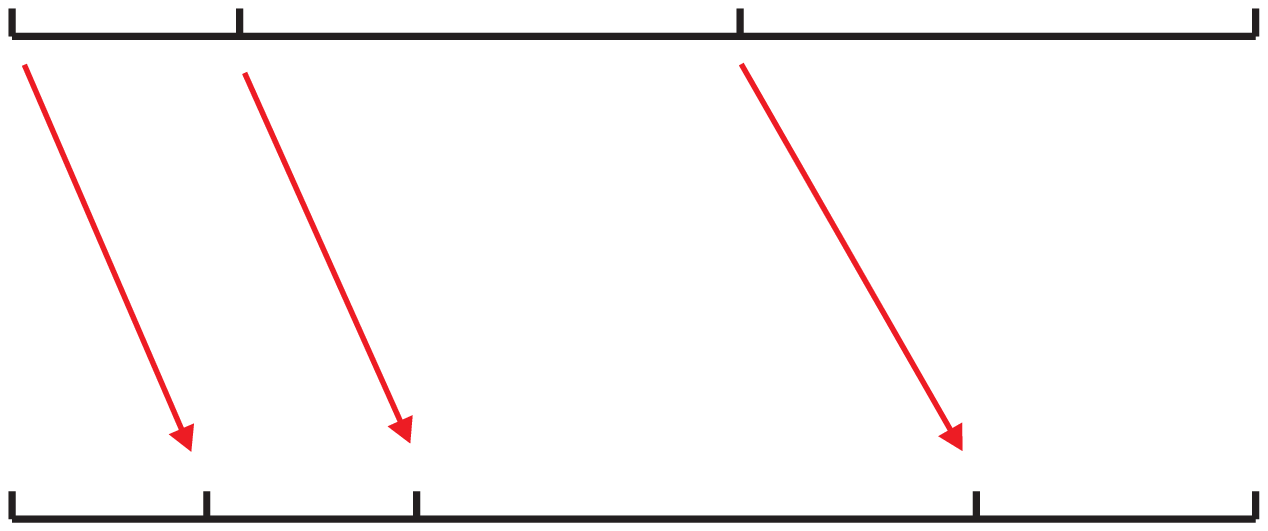 scaled 450}}
\rput(4.2,2.9){$1$}\rput(5.2,2.9){$2$}
\rput(7.5,2.9){$m$}\rput(9.85,2.9){$n$}
\rput(4.2,0.1){$1$}\rput(5,0.1){$i_1$}
\rput(6,0.1){$i_2$}\rput(8.6,0.1){$i_m$}
\rput(9.85,0.1){$n$}
\rput(9,1.5){${\red s_Y}$}
}
\rput(10.5,1.5){\BoxedEPSF{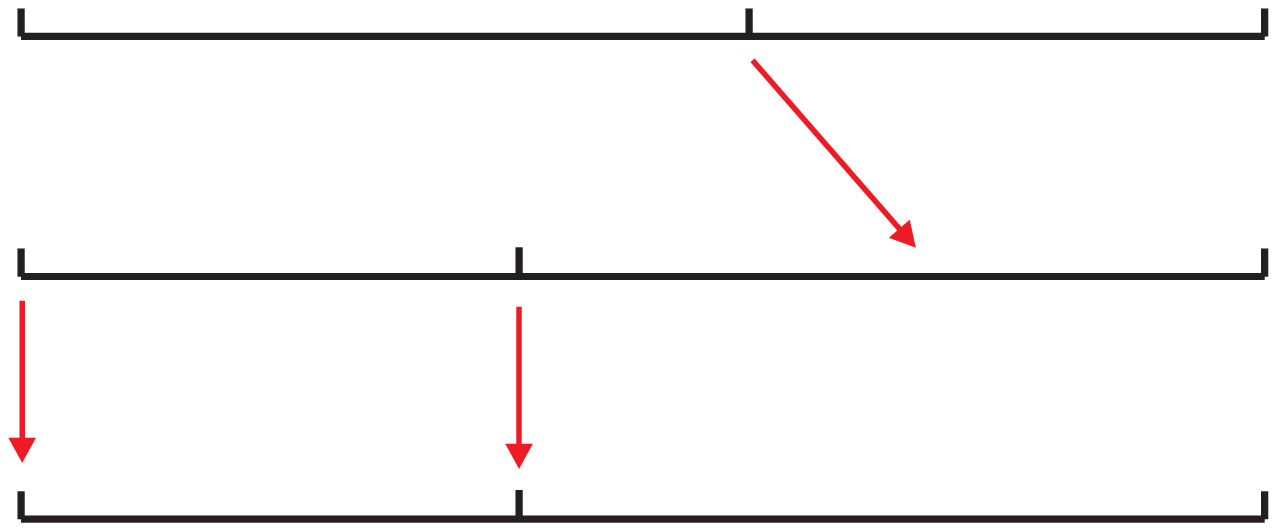 scaled 450}}
\rput(7.7,2.9){$1$}
\rput(11,2.9){$m$}\rput(13.35,2.9){$n$}
\rput(7.7,1.8){$1$}
\rput(10,1.8){$\ell$}\rput(13.35,1.8){$n$}
\rput(13,2){${\red s_Y}$}
\rput(13,0.9){${\red f}$}
\rput(7.7,0.1){$1$}
\rput(10,0.1){$\ell$}\rput(13.35,0.1){$n$}
}
\end{pspicture}
\caption{The $s_Y\in\LL_e$ in a transversal $T$ \emph{(left)\/} and
  $fs_Y$ \emph{(right)\/}.}
  \label{fig:transversal}
\end{figure}
be the map $s_Y:j\mapsto i_j$ given on the left of Figure
\ref{fig:transversal}. We take $T=\{s_Y\}$ to be the transversal used in
the induction process, for $Y$ ranging 
over all $m$-subsets of $[n]$. Thus
$$
V\uparrow I_n=U=\bigoplus_{s_Y} V_Y
$$
where $V_Y$ is the vector space consisting of the vectors $s_Y\otimes
v$ for $v\in V$. The $I_n$-action on $U$ is given by (\ref{eq:4}). 

We claim the following: if $f$ is an idempotent, then
$f(V\uparrow I_n)\not=0$ exactly when $f$ lies in a $\JJ$-class $\JJ_\ell$
with $\JJ_m\leq\JJ_\ell$.
Moreover, $e(V\uparrow I_n)$ is itself isomorphic, as an
$\Symm$-representation, to $V$. 

We choose $f$ conveniently in its $\JJ$-class $\JJ_\ell$:
$f=\id:[\ell]\ra[\ell]$. We have $\JJ_\ell<\JJ_m$
exactly when $\ell<m$, in which case the right part of Figure
\ref{fig:transversal} shows that $m\not\subseteq\dom(fs_Y)$ for any
$m$-subset $Y$. Hence $fs_Y\not\in\LL_e$ (the partial bijections with
domain $[m]$) for any $Y$, and so by (\ref{eq:4})
$$
f\cdot(s_Y\otimes v)=0
$$
for all $Y$ and all $v$. Thus $f(V\uparrow I_n)=0$ when
$\JJ_\ell<\JJ_m$.
If now $f=e$ then we have
$$
es_Y\in\LL_e\Leftrightarrow\dom(e s_Y)=[m]\Leftrightarrow
Y=[m]\Leftrightarrow s_Y=e
$$
in which case 
$$
e\cdot(s_Y\otimes v)\not=0\Leftrightarrow s_Y\otimes v=e\otimes v.
$$
The map $e\otimes v\mapsto v$ is then an isomorphism of vector
spaces $e(V\uparrow I_n)\ra V$, and for any $g\in\Symm$ the diagram
$$
\begin{pspicture}(0,0)(14,2)
\rput(5.75,2){$e\otimes v$}
\rput(8.25,2){$e\otimes gv $}
\rput(5.75,0){$v$}
\rput(8.25,0){$gv$}
\psline[linewidth=0.5pt]{->}(6.2,2)(7.7,2)
\psline[linewidth=0.5pt]{->}(6,0)(7.9,0)
\psline[linewidth=0.5pt]{->}(5.75,1.8)(5.75,0.2)
\psline[linewidth=0.5pt]{->}(8.25,1.8)(8.25,0.2)
\rput(7,2.2){${\scriptstyle g(-)}$}
\rput(7,0.2){${\scriptstyle g(-)}$}
\end{pspicture}
$$
commutes. Thus, the $\Symm$-representations $e(V\uparrow I_n)$ and $V$ are
isomorphic as claimed. 
Finally, Exercise \ref{exercise:reduction:I_n} part 1 gives that $f(V\uparrow
I_n)\not=0$ when $\JJ_\ell>\JJ_m$. This establishes all our claims. 

In particular, $\JJ_m$ is the apex of the $I_n$-representation
$V\uparrow I_n$, and so $V\uparrow I_n$ is indeed in $\irr_m(I_n)$. 

\paragraph{The composition
  $\irr(\Symm)\rightarrow\irr_m(I_n)\rightarrow\irr(\Symm)$:} We have
just seen, for $V$ an irreducible $\Symm$-representation, 
that $e(V\uparrow I_n)\cong V$. Thus $(V\uparrow
I_n)\downarrow\Symm\cong V$, and the composition is the identity.

\paragraph{The composition
$\irr_m(I_n)\rightarrow\irr(\Symm)\rightarrow\irr_m(I_n)$:} we now
show that $(V\downarrow\Symm)\uparrow I_n\cong V$ when $V$ is
an irreducible $I_n$-representation with apex $\JJ_V=\JJ_m$. The strategy is
to reconstruct the representation $(V\downarrow\Symm)\uparrow I_n$ inside $V$. 

We have already the idempotent $e=\id:[m]\ra[m]$ and the transversal
$T=\{s_Y\}$ in Figure \ref{fig:transversal}
for the $m$-sized subsets $Y$ of
$[n]$. 

Consider now the subspaces $(s_Ye)V$ of $V$ for the various 
$Y$. Then:
\begin{description}
\item[--] Each vector space $(s_Ye)V$ is isomorphic to $eV$: the
  linear map
  $eV\ra(s_Ye)V$ given by 
  $ev\mapsto (s_Ye)v$ has inverse the map $(s_Ye)v\mapsto
  s_Y^*(s_Ye)v=e^2v=ev$, and so is an isomorphism. 
\item[--] The sum $\sum_Y (s_Ye)V$ of these spaces is direct: for
  which we need to show that for a fixed subset $Y$, the intersection
  \begin{equation}
    \label{eq:11}
s_YeV\cap \sum_{Z\not=Y}s_ZeV
  \end{equation}
is the zero space.
We have just seen that $s_Y^*$ gives an isomorphism $s_YeV\ra eV$, hence
maps the subspace $s_YeV\cap \sum_{Z\not=Y}s_ZeV$ of $s_YeV$ isomorphically onto
its image in $eV$. But
\begin{equation}
  \label{eq:10}
s_Y^*\biggl(s_YeV\cap \sum_{Z\not=Y}s_ZeV\biggr)
\subseteq
s_Y^*s_YeV\cap s_Y^*\biggl(\sum_{Z\not=Y}s_ZeV\biggr)
=
eV \cap \sum_{Z\not=Y}s_Y^*s_ZeV  
\end{equation}
where $Z\not= Y$ gives that the domain of $s_Y^*s_Ze$ has size
strictly less than $m$, and so $s_Y^*s_Ze$ lies in a $\JJ$-class lower
down the strategic picture than $\JJ_m$ does. As $\JJ_m$ is the
apex of $V$ we have $s_Y^*s_ZeV=0$ for all $Z$, so that the right hand
side of (\ref{eq:10}) is $0$, and hence (\ref{eq:11}) is too.
\item[--] Restricting the $S$-action on $V$ to the subspace
  $\bigoplus_Y s_YeV$: if $t\in I_n$ then there are two
  possibilities for the product $ts_Y$. Either:
  \begin{description}
  \item[(i).]$ts_Y\in\LL_e$, in which case by (\ref{eq:1}), there is
    a $g\in G_e$ and an $m$-subset $Z$ such that $ts_Y=s_Zg$; or
  \item[(ii).] $ts_Y\not\in\LL_e$, and since this $\LL$-class consists
    of all partial bijections with domain $[m]$, and
    $\dom(ts_Y)\subseteq [m]$, we have that $\dom(ts_Y)$ is a proper
    subset of $[m]$. In particular $ts_Y$ lies in a $\JJ$-class
    lower down the strategic picture than $\JJ_m$.
  \end{description}
The $S$-action on $\bigoplus_Y s_Y eV$ is therefore given by
$$
  t\cdot(s_Ye)\cdot v
=
\left\{
\begin{array}{ll}
(s_Ze)\cdot (g\cdot v),&\text{ if }ts_Y\in\LL_e,\text{ or}\\
0, &\text{ else.}
\end{array}
\right.
$$
\end{description}
We conclude, first of all, that the subspace $\bigoplus_Y s_YeV$ is
in fact a subrepresentation of $V$; moreover $\bigoplus_Y s_YeV$ contains, by
taking $Y=[m]$, the subspace $eV\not=0$. Thus $\bigoplus_Y s_YeV$ is
a non-zero subrepresentation of the irreducible representation $V$, hence 
$$
\bigoplus_Y s_YeV=V.
$$
Finally, if $ts_Y\in\LL_e$ then the diagram
$$
\begin{pspicture}(0,0)(14,2.5)
\rput(0,0.25){
\rput(5.5,2){$(s_Ye)\cdot v$}
\rput(5.5,0){$s_Y\otimes e\cdot v$}
\rput(8.5,2){$(s_Ze)\cdot (g\cdot v)$}
\rput(8.5,0){$s_Z\otimes (eg)\cdot v$}
\psline[linewidth=0.5pt]{->}(6.2,2)(7.4,2)
\psline[linewidth=0.5pt]{->}(6.3,0)(7.5,0)
\psline[linewidth=0.5pt]{->}(5.5,1.8)(5.5,0.2)
\psline[linewidth=0.5pt]{->}(8.5,1.8)(8.5,0.2)
\rput(7,2.2){${\scriptstyle t(-)}$}
\rput(7,0.2){${\scriptstyle t(-)}$}
}
\end{pspicture}
$$
commutes (it trivially commutes if $ts_Y\not\in\LL_e$). Thus
$(V\downarrow\Symm)\uparrow I_n\cong V$ as $I_n$-representations, and the 
composition $\irr_m(I_n)\rightarrow\irr(\Symm)\rightarrow\irr_m(I_n)$
is the identity map. 

This completes the proof of the Clifford-Munn correspondence when $S=I_n$.

\begin{vexercise}
\label{exercise:CMP}
Mimic the proof above for an inverse monoid $S$ of the form $S=\sgl$
(\emph{hint:\/} much of the proof can be found scattered among what we
have already said).
\end{vexercise}


\begin{example}[The irreducibles of $I_n$]
\label{induction:specht}
We are finally in a position to describe the irreducible
representations over $\C$ of the symmetric inverse monoid $I_n$. 
By Theorem \ref{theorem: Munn:Oganesyan}, every $I_n$-representation
over $\C$ is a direct sum of these. 
As we will be
doing things this way in \S\ref{section:sexy:example} -- and this is
sort of a dry run at it -- 
we will use the $S\kern-1pt(\Symn,L)$ description of $I_n$ that we saw at the end of
\S\ref{section:semigroup_basics}, where $L$ is the
lattice of subsets of $[n]$. This allows us to follow the recipe for induction
given at the end of \S\ref{section:induction}. 

Fix an $m$ in the range $0\leq m\leq n$, hence a $\JJ$-class
corresponding to the $\Symn$-orbit on $L$ consisting of the subsets
of $[n]$ having size
$m$. Let $a=\{1,2,\ldots,m\}$ and $G_a$ the maximal
subgroup containing the idempotent $\id_a$. The elements of $G_a$ are
the $g_a$ where $g\in\Symn$ is such that $g\cdot a=a$ (rather than
being
the bijections $a\rightarrow a$
as they would
be in the ``usual'' way of describing $I_n$). Finally, let $\lambda$ be a
partition of $m$ and $S^\lambda$ be the Specht representation spanned by the
$v_T$ in (\ref{eq:19}) as $T$ ranges over the tableau of shape
$\lambda$. 

We will describe the representation $S^\lambda\uparrow I_n$. 
The Clifford-Munn correspondence tells us that the $S^\lambda\uparrow
I_n$, as both $\lambda$ and $m$ vary in $\lambda \vdash m$, form a
complete and non-redundant list of the $I_n$-irreducibles over $\C$. 

If
$b=\{i_1,\ldots,i_m\}$ is a subset of $[n]$ of size $m$, then let
$\beta$ be an element of $\Symn$ that sends $j\in a$ to $i_j\in b$. We then
take the transversal $T$ needed for induction to be the resulting $\beta_a$ as $b$ ranges
over the subsets of size $m$. 

If $T$ is a tableau of shape $\lambda$ filled with entries from $a$, then $\bb\cdot
T$ is a tableau of shape $\lambda$ filled with entries from $b$. Let $S^{\lambda,b}$
be a copy of $S^\lambda$, spanned by the 
$$
b\otimes v_T=\sum_{h\in c_{\bb\cdot T}}\text{sign}(h)\, h\cdot\{\bb\cdot
T\},
$$
as $T$ varies over the tableau (on $a$), and 
where $c_{\bb\cdot T}$ are those elements of the symmetric group on
the set $b$ preserving the columns of $T$. The vector $b\otimes v_T$
is just the vector $v_T$, but with every occurence of $j\in a$ in a tabloid
replaced by $i_j\in b$, and $S^{\lambda,b}$ is the space spanned by
the $b\otimes v_T$.

The representation $S^\lambda\uparrow I_n$ acts on the space
$$
S^\lambda\uparrow I_n
=
\bigoplus_{|b|=m}
S^{\lambda,b}
$$
To see how, fix an $s=g_c\in I_n$. We saw at the end of \S\ref{section:induction}
that the apex of
$S^\lambda\uparrow I_n$
is the $\JJ$-class containing the maximal subgroup $G_a$ that we started with, so
if $|c|< m$ we get $s\cdot S^\lambda\uparrow
I_n=0$. 

On the other hand, by (\ref{eq:18}), 
if $c$ has size at least $m$, then it will \emph{not\/} kill
those summands $S^{\lambda,b}$ for which $b\subseteq c$. In this case
$s\,\bb_a$ lies in the $\HH$-class labelled by the subset $d=g\cdot
b$, so that for $b\otimes v_T$ spanning $S^{\lambda,b}$ we get 
$$
s\cdot(b\otimes v_T)=d\otimes h\cdot v_T
$$
where $h=(\delta^{-1}g\,\bb)_a\in G_a$.
\end{example}


\section{A sexy example}
\label{section:sexy:example}

For the purposes of these notes, ``sexy'' will mean a certain family
of Renner monoids. These encode much of the structure of algebraic
monoids, and are ubiquitous in nature. 

We first set the examples up in the form $\sgl$ from Section
\ref{section:semigroup_basics}. As usual $G$ is the 
symmetric group $\Symn$, but the lattice is one we haven't seen
before. Let $L_0$ consist of the \emph{ordered\/} partitions of $[n]$,
i.e. the \emph{tuples\/} $\Lambda=(\Lambda_1,\ldots,\Lambda_p)$ with
$\{\Lambda_1,\ldots,\Lambda_p\}$ a partition of $[n]$. Partially order
the ordered partitions via 
$(\Lambda_1,\ldots,\Lambda_p)\leq (\Delta_1,\ldots,\Delta_q)$ if and
only if
\begin{align*}
  \label{eq:16}
 &\text{\textendash}\text{ each }\Lambda_i\subseteq\text{some
   }\Delta_j,\text{ and}\\
 &\text{\textendash}\text{ if }\Lambda_i\subseteq\Delta_j\text{ and for
   }i<k\text{ we have }\Lambda_k\subseteq\Delta_\ell,\text{ then }j\leq\ell.
\end{align*}
$L_0$ then has maximum element the ordered partition $([n])$ with a
single block and \emph{minimal\/} elements the ordered partitions
where every block has size one; these minima are in 1-1 correspondence
with the permutations of $[n]$. 

Formally adjoin a minimum
$\0$ to $L_0$ to get the lattice $L$. The $\Symn$-action on $L$ is the
usual $g\cdot
(\Lambda_1,\ldots,\Lambda_p)=(g\cdot\Lambda_1,\ldots,g\cdot\Lambda_p)$
together with $g\cdot\0=\0$.

\paragraph{A short diversion on where the example comes from.} A
linear algebraic group $\ams{G}$, over an algebraically closed field
$k$, is an affine algebraic variety over $k$,
together with a morphism $\varphi:\ams{G}\times \ams{G}\rightarrow \ams{G}$
of varieties, such that the product $gh:=\varphi(g,h)$ gives
$\ams{G}$ the structure of a group. Generalising this idea, a
\emph{linear algebraic monoid\/} $\ams{M}$ arises when
$\varphi:\M\times\M\ra\M$
gives $\ams{M}$ the structure of a monoid. 

The canonical examples are $\ams{G}=\gl_nk$, the group of invertible
matrices over $k$, and $\ams{M}=\m_nk$, the monoid of all
$n\times n$ matrices over $k$ (both under multiplication). In fact
$\gl_nk$ is the group of units 
of the monoid $\m_nk$, and indeed for sensible $\ams{M}$ the group 
of units $\ams{G}$ is an algebraic group with Zariski closure
$\overline{\ams{G}}=\ams{M}$. 

There is a standard construction of algebraic monoids that 
starts with a sensible
algebraic group $\ams{G}_0$ and a sensible representation 
$f:\ams{G}_0\rightarrow GL(V)$. The resulting algebraic
monoid is then $\ams{M}=\ov{k^\times f(\ams{G}_0)}\subset\m_{m}k$ with
group of units 
$\ams{G}=k^\times f(\ams{G}_0)$. For example, if $\ams{G}_0=\sl_n, \so_n$ and
$\sp_n$ and $f$ is the natural representation of $\ams{G}_0$,
then the resulting $\ams{M}$ are the \emph{classical monoids\/}: the general linear monoids
$\m_n=\ov{k^\times \sl_n}$, the 
orthogonal monoids $\mso_n=\ov{k^\times \so_n}$ and the symplectic
monoids $\msp_n=\ov{k^\times \sp_n}$.  

Associated to a (reductive) algebraic group $\ams{G}$ is a finite group --
called the Weyl group -- that encodes much of the structure of
$\ams{G}$; for a  (reductive) algebraic monoid $\ams{M}$ there
is a finite inverse monoid $R$ --
called the \emph{Renner monoid\/} -- that plays an analogous role.
For example, the Weyl group of $\gl_nk$ is the symmetric group $\Symn$
and the Renner monoid of $\m_nk$ is the symmetric inverse monoid
$I_n$. In general the group of units of the Renner monoid $R$ of
$\ams{M}$ is the Weyl group $W$ of the algebraic group of units $\ams{G}$ of $\ams{M}$.

If $R$ is the Renner monoid of the algebraic monoid $\ams{M}$ having
group of units the algebraic group $\ams{G}$, then $R$ is an inverse
monoid of the form
$\swl$: the group $W$ is the Weyl group of 
$\ams{G}$ and the lattice $L$ turns out to be the face lattice of a convex
polytope. 

To see what this means, a polytope $P$ in $\R^m$ is the convex hull of a finite set
of points. It has $r$-dimensional faces, for $-1\leq r\leq m$, with
the $0$-dimensional faces being the vertices, $1$-dimensional faces
the edges, and so on, with $P$ itself the unique $m$-dimensional
face; for formal reasons (mainly so that we get a lattice below) we
take the empty set $\varnothing$ to be 
the unique face of dimension $-1$. The \emph{face lattice\/} of $P$ 
consists of the faces ordered by inclusion; it is a lattice with meet
$\ss\wedge\tau$ the intersection and join $\ss\vee\tau$ the smallest
face containing both $\ss$ and $\tau$.

The Renner monoid of $\m_nk$ has the form $\swl$ where $W$ is the
symmetric group $\Symn$ and $L$ is the face lattice of an
$(n-1)$-dimensional simplex. If $[n]=\{1,\ldots,n\}$ are the labels of
the vertices of the simplex, 
then $L$ is the lattice of subsets of $[n]$ ordered by
inclusion, and the $\Symn$-action on $L$ is the usual one. This is
the description of $I_n$ we gave at the end of
\S\ref{section:semigroup_basics}. 

Now to the example we are interested in: let $\ams{G}_0=\sl_n$ and
$V_0$ be the natural
module for $\ams{G}_0$. Let $\bigwedge^p V_0$ be the $p$-th exterior
power of $V_0$
and finally
$$
V=\bigotimes_{p=1}^{n-1}\bigwedge^p V_0,\text{ with }
\dim V:=m=\prod_{p=1}^{n-1}\binom{n}{p}.
$$
If $f:\ams{G}_0\rightarrow GL(V)$ is the corresponding representation
then let $\ams{M}=\ov{k^\times f(\ams{G}_0)}$ and let $R$ be
the Renner monoid of $\ams{M}$. Then $R\cong\swl$ with $W$ the
symmetric group $\Symn$ and $L$ the face lattice of the
$(n-1)$-dimensional permutohedron. This is the polytope in $\R^n$
obtained by taking the convex hull of the $n!$ points arising 
from all permutations of the coordinates of the point
$(1,2,\ldots,n)$. As all these points lie in the hyperplane with
equation $x_1+x_2+\cdots+x_n=1+2+\cdots+n$, the polytope is actually
$(n-1)$-dimensional.  

\begin{figure}
  \centering
\begin{pspicture}(0,0)(14,7)
\rput(2.5,3){\BoxedEPSF{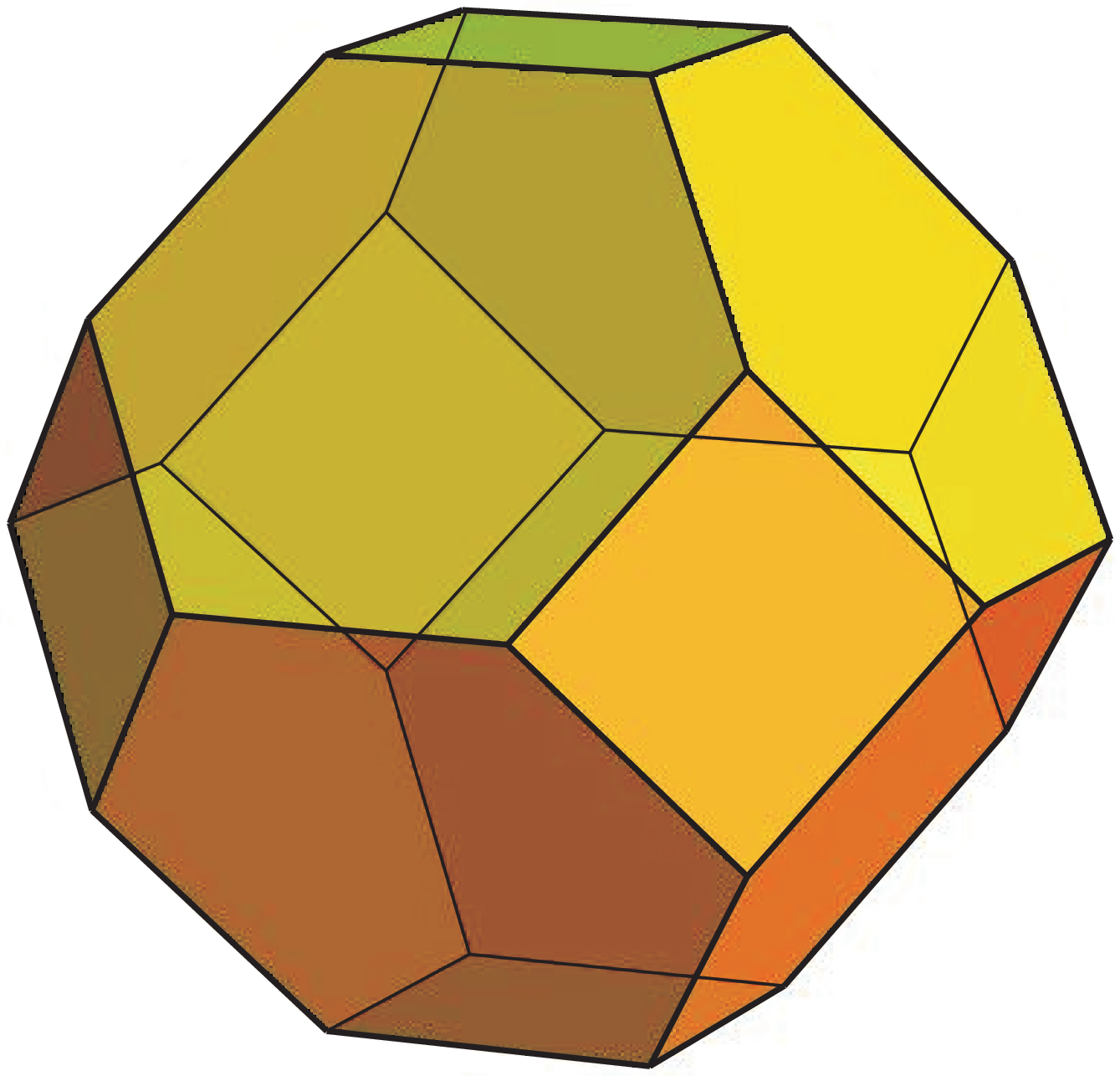 scaled 250}}
\rput(0.75,0.5){
\rput(9.15,3){\BoxedEPSF{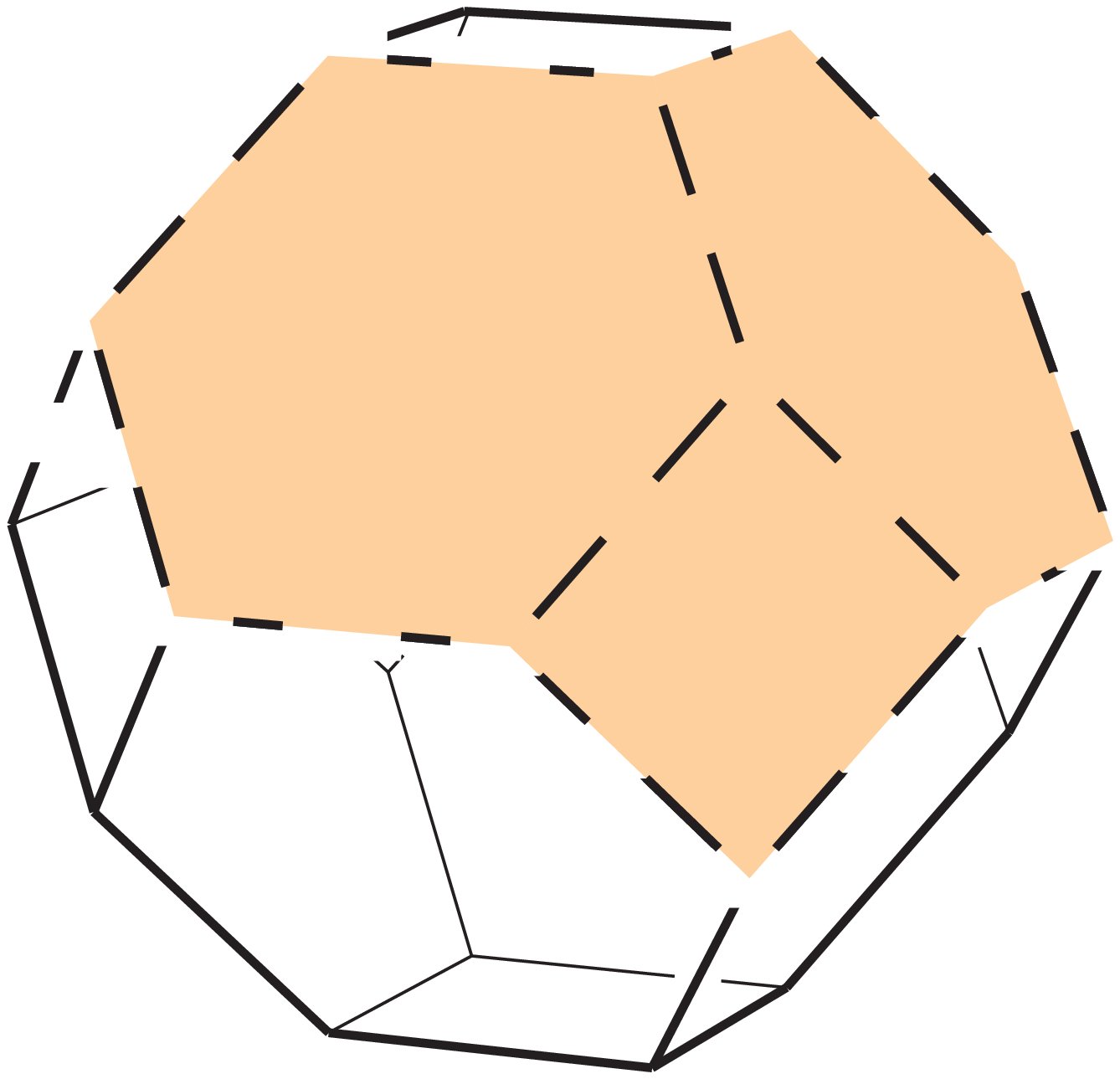 scaled 450}}
\rput(10.05,3.9){$\scriptscriptstyle{{\blue 1,2,3,4}}$}
\rput(9.55,5.5){$\scriptscriptstyle{{\blue 1,3,2,4}}$}
\rput(11.5,4.5){$\scriptscriptstyle{{\blue 1,4,3,2}}$}
\rput(8.75,2.4){$\scriptscriptstyle{{\blue 2,1,3,4}}$}
\rput(6.9,2.55){$\scriptscriptstyle{{\blue 2,3,1,4}}$}
\rput(6.45,4.2){$\scriptscriptstyle{{\blue 3,2,1,4}}$}
\rput(7.75,5.65){$\scriptscriptstyle{{\blue 3,1,2,4}}$}
\rput(10.05,1.15){$\scriptscriptstyle{{\blue 2,1,4,3}}$}
\rput*(9.5,0.1){$\scriptscriptstyle{{\blue 2,4,1,3}}$}
\rput*(10.2,0.5){$\scriptscriptstyle{{\blue 4,2,1,3}}$}
\rput(11.3,2.6){$\scriptscriptstyle{{\blue 1,2,4,3}}$}
\rput*(7.7,0.3){$\scriptscriptstyle{{\blue 2,4,3,1}}$}
\rput*(6.4,1.5){$\scriptscriptstyle{{\blue 2,3,4,1}}$}
\rput*(6,3.1){$\scriptscriptstyle{{\blue 3,2,4,1}}$}
\rput(12,3){$\scriptscriptstyle{{\blue 1,4,2,3}}$}
\rput(10.25,5.8){$\scriptscriptstyle{{\blue 1,3,4,2}}$}
\rput*(8.5,5.9){$\scriptscriptstyle{{\blue 3,1,4,2}}$}
\rput*(11.45,2){$\scriptscriptstyle{{\blue 4,1,2,3}}$}
%
\rput(9.8,4.7){$\scriptscriptstyle{{\red 1,23,4}}$}
\rput(9.4,3.15){$\scriptscriptstyle{{\red 12,3,4}}$}
\rput(7.85,2.5){$\scriptscriptstyle{{\red 2,13,4}}$}
\rput(6.7,3.45){$\scriptscriptstyle{{\red 23,1,4}}$}
\rput(7.05,4.9){$\scriptscriptstyle{{\red 3,12,4}}$}
\rput(8.65,5.6){$\scriptscriptstyle{{\red 13,2,4}}$}
\rput(9.35,1.85){$\scriptscriptstyle{{\red 2,1,34}}$}
\rput(10.7,3.25){$\scriptscriptstyle{{\red 1,2,34}}$}
\rput(10.65,1.9){$\scriptscriptstyle{{\red 12,4,3}}$}
\rput*(8.5,0.2){$\scriptscriptstyle{{\red 2,4,13}}$}
\rput*(7,0.9){$\scriptscriptstyle{{\red 2,34,1}}$}
\rput*(10.9,1.4){$\scriptscriptstyle{{\red 4,12,3}}$}
\rput(11.75,3.75){$\scriptscriptstyle{{\red 1,4,23}}$}
\rput(10.9,5.15){$\scriptscriptstyle{{\red 1,34,2}}$}
\rput*(6.2,2.35){$\scriptscriptstyle{{\red 23,4,1}}$}
\rput*(9.8,0.7){$\scriptscriptstyle{{\red 2,14,3}}$}
\rput*(6.7,2.1){$\scriptscriptstyle{{\red 2,3,14}}$}
\rput(6.15,3.55){$\scriptscriptstyle{{\red 3,2,14}}$}
\rput(8.15,3.9){$123,4$}
\rput(10.05,2.5){$12,34$}
\rput(10.8,4.2){$1,234$}
}
\end{pspicture}
\caption{The $3$-permutodehron, after identifying the hyperplane 
$x_1+x_2+x_2+x_4=10\subset\R^4$ with $\R^3$ \emph{(left)\/} 
and part of the lattice of ordered partitions overlaid on the faces
\emph{(right)\/}, with the blocks of the partitions separated by commas.}
  \label{fig:permutohedron}
\end{figure}

The face lattice of the permutohedron is isomorphic to the lattice $L$
of ordered
partitions, with $\0$ adjoined, described at the beginning of the section. Figure
\ref{fig:permutohedron} shows the $n=4$ case. 

\paragraph{} To describe the irreducible representations over $\C$ of
our Renner monoid, we use the $\sgl$ description from the beginning
of the section, and
start by drilling down a little more into
the structure of the monoid, following
\S\ref{section:semigroup_basics}. 

First, we have our usual ambiguity with the elements of $\sgl$,
where $g_a=h_b$ when $a=b$ and $g^{-1}h\cdot c=c$ for all $c\leq
a$. In this case it turns out to disappear.  If $a$ is the ordered partition
$(\Lambda_1,\ldots,\Lambda_p)$ and 
$c$ is a minimal element $\not=\0$ with the property that $c\leq a$,
then $c$ has the form 
$$
c=(\{x_{11},\ldots,x_{1q_1}\},\ldots,\{x_{p1},\ldots,x_{pq_p}\})
$$
where 
$\Lambda_1=\{x_{11},\ldots,x_{1q_1}\},\ldots,\Lambda_p=\{x_{p1},\ldots,x_{pq_p}\}$. 
If $k$ is an
element of $\Symn$ with $k\cdot c=c$ then $k=\id$. Thus $g_a=h_b$ iff
$a=b$ and $g=h$. 

\begin{figure}
  \centering
\begin{pspicture}(0,0)(4,6)
\rput(2,3){\BoxedEPSF{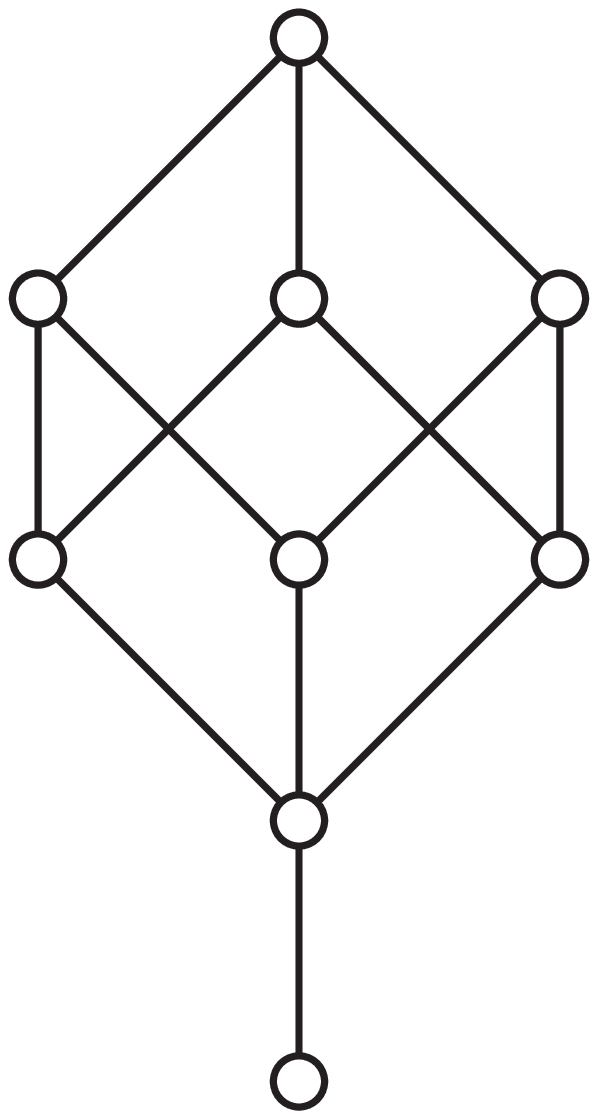 scaled 450}}
\rput(1.45,1.8){${\scriptstyle (1,1,1,1)}$}
\rput(0.3,3){${\scriptstyle (2,1,1)}$}
\rput(1.5,3){${\scriptstyle (1,2,1)}$}
\rput(2.7,3){${\scriptstyle (1,1,2)}$}
\rput(0.4,4.2){${\scriptstyle (3,1)}$}
\rput(1.6,4.2){${\scriptstyle (2,2)}$}
\rput(2.8,4.2){${\scriptstyle (1,3)}$}
\rput(2,5.7){${\scriptstyle (4)}$}
\rput(2,0.3){${\scriptstyle\0}$}
\end{pspicture}
\caption{The $\JJ$-class poset of our Renner monoid $R$ when $n=4$,
  corresponding to the poset of compositions of $4$.}
  \label{fig:permutohedron2}
\end{figure}

If $(\Lambda_1,\ldots,\Lambda_p)$ is an ordered partition with
$\lambda_i=|\Lambda_i|$, then the ordered tuple
$(\lambda_1,\ldots,\lambda_p)$ is called a \emph{composition\/} of
$n$: namely, the $\lambda_i$ are totally ordered with
$\sum\lambda_i=n$. Call the composition the \emph{type\/} of the
ordered partition. Two ordered partitions are then in the same
$\Symn$-orbit when they have the same type, and the $\JJ$-class poset
has elements the compositions ordered by
$(\lambda_1,\ldots,\lambda_p)\leq (\mu_1,\ldots,\mu_q)$ whenever
$$
(\lambda_1,\ldots,\lambda_p)=
(\lambda_{11},\ldots,\lambda_{1m_1},\ldots,\lambda_{p1},\ldots,\lambda_{p,m_p})
$$
with $\mu_i=\lambda_{i1}+\ldots+\lambda_{i,p_i}$. Figure 
\ref{fig:permutohedron2} shows this poset when $n=4$.

Fix a composition $(\lambda_1,\ldots,\lambda_p)$ and consider the
$\JJ$-class of ordered partitions of this type. 
If $a$ is one of them, then the
maximal subgroup is $G_a=\{g_a:g\cdot a=a\}$, and this is isomorphic to
the Young subgroup $S_{\kern-1pt\lambda_1}\times\cdots\times
S_{\kern-1pt\lambda_p}$ of $\Symn$. We saw above that our usual
ambiguity in expressing elements vanishes in $R$; this is why there
is no need to form a quotient when describing $G_a$. 

Let
$a=(\Lambda_1,\ldots,\Lambda_p)$ be the ordered partition of type
$(\lambda_1,\ldots,\lambda_p)$
given by:
\begin{equation}
  \label{eq:20}
\Lambda_1=\{1,\ldots,\lambda_1\},
\ldots,
\Lambda_p=\{n-\lambda_p+1,\ldots,n\},  
\end{equation}
We now describe the irreducible representations of $R$ that arise by
inducing up those of the maximal subgroup $G_a\cong
S_{\kern-1pt\lambda_1}\times\cdots\times 
S_{\kern-1pt\lambda_p}$.
Varying the composition produces a complete list of the irreducibles
over $\C$ of the Renner monoid $R$. We borrow one more fact from the
representation theory of finite groups: if $\{V_i\}_{i\in I}$ and
$\{U_j\}_{j\in J}$ are the irreducibles of the groups $G$ and $H$, then the
$\{V_i\otimes U_j\}_{I\times J}$ are the irreducibles of $G\times H$,
with $V_i\otimes U_j$ a $(G\times H)$-representation
via the action $(g,h)\cdot v\otimes u=g\cdot v\otimes h\cdot u$. 

Now fix partitions $\mu_1\vdash\lambda_1,\ldots,\mu_p\vdash\lambda_p$ and
consider the irreducible $(S_{\kern-1pt\lambda_1}\times\cdots\times
S_{\kern-1pt\lambda_p})$-representation
\begin{equation}
  \label{eq:21}
S^{\mu_1}\otimes\cdots\otimes S^{\mu_p}  
\end{equation}
where $S^{\mu_i}$ is the Specht representation of $S_{\kern-1pt\lambda_i}$
corresponding to the partition $\mu_i\vdash\lambda_i$. The
representation (\ref{eq:21})
is spanned by the vectors 
$$
v_{T_1}\otimes\cdots\otimes v_{T_p}
$$
defined in (\ref{eq:19}) and
as the $T_i$ range over the tableau of shape $\mu_i$ filled with the numbers $\Lambda_i$
in (\ref{eq:20}). To describe $S^{\mu_1}\otimes\cdots\otimes
S^{\mu_p}\uparrow R$, let $b=(\Delta_1,\ldots,\Delta_p)$ be another ordered
partition of type $(\lambda_1,\ldots,\lambda_p)$ and let
$\beta\in\Symn$ be such that $\beta:\Lambda_i\mapsto\Delta_i$ in an
order preserving way, i.e. if $x<y\in\Lambda_i$ then
$\beta(x)<\beta(y)\in\Delta_i$. Let $S^{\mu_1,\,\bb}\otimes\cdots\otimes
S^{\mu_p,\,\bb}$ be a copy of (\ref{eq:21}) with spanning vectors of
the form:
$$
\bb\otimes
v_{T_1}\otimes\cdots\otimes v_{T_p},
$$
defined to be $v_{T_1}\otimes\cdots\otimes v_{T_p}$, but with every
occurence of $j$ in a tabloid replaced by $\bb(j)$. (Warning: this
vector is linear in the $v_{T_i}$ coordinates only; the
``$\bb\,\otimes$", as usual, is just notation that comes along for the
ride). 

The representation 
$S^{\mu_1}\otimes\cdots\otimes
S^{\mu_p}\uparrow R$
is carried by the space
$$
S^{\mu_1}\otimes\cdots\otimes
S^{\mu_p}\uparrow R
=
\bigoplus_{b} 
S^{\mu_1,\,\bb}\otimes\cdots\otimes
S^{\mu_p,\,\bb}
$$
with the direct sum over the ordered partitions of type
$(\lambda_1,\ldots,\lambda_p)$. Let $s=g_c\in R$ with
$g\in\Symn$ and $c$ the ordered partition $(\Omega_1,\ldots,\Omega_q)$ of type
$(\omega_1,\ldots,\omega_q)$. If
$(\lambda_1,\ldots,\lambda_p)\not\leq(\omega_1,\ldots,\omega_q)$ then
$$
s\cdot (S^{\mu_1}\otimes\cdots\otimes
S^{\mu_p}\uparrow R)=0.
$$
Otherwise, when 
$(\lambda_1,\ldots,\lambda_p)\leq(\omega_1,\ldots,\omega_q)$
we have $s\cdot (S^{\mu_1,\,\bb}\otimes\cdots\otimes
S^{\mu_p,\,\bb})\not=0$ when $b=(\Delta_1,\ldots,\Delta_p)\leq
(\Omega_1,\ldots,\Omega_q)$ and in this case
$$
s\cdot(\bb\otimes
v_{T_1}\otimes\cdots\otimes v_{T_p})
=\delta\otimes h\cdot (v_{T_1}\otimes\cdots\otimes v_{T_p})
$$
where $d=g\cdot b$ and 
$h=(\delta^{-1}g\,\bb)_a\in G_a\cong S_{\kern-1pt\lambda_1}\times\cdots\times
S_{\kern-1pt\lambda_p}$.


\section*{Notes and References}


There are numerous books that deal with semigroup representations,
starting with the classic \cite{MR0132791}*{Chapter 5}; more modern sources are 
\cite{MR2460611,MR3525092}. The reader who has got this far will see
large overlap with \cite{MR2460611}, making \cite{MR3525092} a good
next step. The original papers of Clifford \cite{MR0006551} and Munn
\cite{MR0066355}-\cite{MR0172953}  are still very readable, as is the later reworking by
Rhodes \cite{MR1142387}.

\paragraph{Semigroups.} The standard reference on semigroups is
\cite{Howie95}, where we have followed Chapters 1, 2 and 5; see
also \cite{MR0132791,Lawson98,MR2460611}. The three running examples are very
much in the style of \cite{MR2460611}. For the reader who is wondering
about the ``full'' transformation semigroup, there is a partial
version $PT_n$, which is a sort-of-amalgam of $I_n$ and $T_n$; see
\cite{MR2460611}*{Chapter 2}. 

The restriction to finite
regular monoids is purely to make things cleaner. An expert (who
shouldn't be reading this anyway) can make the appropriate
adjustments, especially in Section \ref{section:clifford:munn}. One
convenience that results is that the relation $\langle \LL,\RR\rangle$,
usually called $\DD$ by semigroup theorists, coincides with $\JJ$. So
all mention of $\DD$ (which gives the eggbox
pictures) has been
merged with $\JJ$ (which gives the partial order
on the eggboxes). Figure \ref{fig:Tn_Hclasses_idempotents} is adapted
from a picture by James East.  

The inverse monoids $\sgl$ appear in \cite{Everitt10}*{Section 9.2} as monoids of
partial permutations, although they are implicit in the
literature. Their purpose in \cite{Everitt10} is to shed light on the
\emph{factorisable\/} 
inverse monoids: these are monoids $S$ with the property that
$S=EG=GE$, where $G$ is the group of units of $S$ and $E$ the
idempotents -- see \cite {Chen74, FitzGerald10}. 
Exercise \ref{exercise:sgl:order} can be done by counting the $g_a$,
but bearing in mind the 
ambiguity; another way, more natural in this context, is to count up
the entries in the boxes in the strategic picture. The monoid of
uniform block permutations of Example \ref{exercise:partition} first
appears in \cite{FitzGerald03}. The picture of the Hasse diagram for
the partition lattice $\Pi(4)$ in Figure \ref{fig:S(G,L)} is based on
one by Tilman Piesk \cite{Baez15}.

\paragraph{Representations.} There are many books on group
representation theory; we have followed the notation and style
of \cite{MR1153249}*{Part I}. In particular the approach is elementary,
\emph{aka\/} ``module-free''. In this section the representations are
over an arbitrary field $k$; one moral to be extracted at the end is
that in dealing with semigroup representations in characteristic $p>0$,
one needs to be just as careful, if not 
more careful, than one does in group representation theory. The
emphasis thus moves to $k=\C$ in later sections. Another omission is
the theory of semigroup characters, which is well developed for the
running examples. 

The restriction to monoids (rather than semigroups) and
monoid homomorphisms removes null representations from consideration
-- this makes a number of statements 
less cluttered.

The standard reference on reflection groups is \cite{Humphreys90}. A
finite reflection group (acting on a real vector space) can be boiled
down to a very concise piece of combinatorial data called a Coxeter
symbol. Starting from a Coxeter symbol one can construct a
representation of the reflection group, called the reflectional
representation; a fundamental result in the theory of reflection
groups is that the reflectional representation is
irreducible. Starting from the type $A$ Coxeter symbol:
$$
\begin{pspicture}(0,0)(14,1)
\rput(7,0.5){\BoxedEPSF{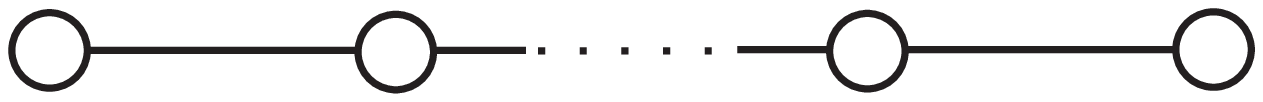 scaled 350}}
\rput(4.9,0.9){$s_1$}\rput(6.1,0.9){$s_2$}
\rput(7.8,0.9){$s_{n-2}$}\rput(9.1,0.9){$s_{n-1}$}
\end{pspicture}
$$
this process gives the reflectional representation of Example
\ref{example:permutation:representation}. The elementary argument showing that
this is irreducible was supplied by Michael
Bate. 

Munn \cite{MR0081910} extends the cycle notation for permutations in
$\Symn$ to elements $s\in I_n$ in the following neat way: for
$x\in[n]$, repeated application of $s$ either results in a \emph{cycle\/}:
$x,s(x),s^2(x),\ldots,$ $s^{k+1}(x)=x$, in which case we write
$(x,s(x),\ldots,s^k(x))$ as usual; or, $s^{k}(x)$ is the first
iteration of $S$ that does not lie in
the domain of $s$, so that no more applications of $s$ can be made. In
this case we have a \emph{link\/} $[x,s(x),\ldots,s^k(x)]$. Any $s\in
I_n$ can then be written uniquely as a juxtaposition of disjoint cycles and
links; the element $[1,2,3]\in I_3$ on the right of Figure
\ref{fig:reflection_rep_Sn} is an example. 
Reflection monoids appear in \cite{Everitt10}, where 
$I_n$ is a Boolean monoid of type $A$.

The formulation of semisimplicity suffers a little from the
module-free approach, where it is cleaner to talk in terms of the
semisimplicity of the semigroup algebra $kS$. 
We have also avoided the notion of decomposability: the mapping
representation of $T_n$ is thus indecomposable but not irreducible,
\emph{even in characteristic $0$\/}. One imagines that this is the
aspect of the whole thing that group theorists find most distressing.
Theorem
\ref{theorem:JordanHolder} is standard -- we have followed
\cite{MR1984740}*{Theorem 6.1.15}; Theorem \ref{theorem:Mashke} similarly
(see e.g. \cite{MR1984740}*{Theorem 3.1.14}); 
Theorem \ref{theorem:
  Munn:Oganesyan} is 
less well known, except to the cognoscenti; 
see \cite{MR0132791,MR3525092}. For $I_n$ and
$T_n$ see also \cite{MR2460611}*{Section 11.5}.

\paragraph{Interlude: the symmetric group.} A standard introductory
text to all aspects of the representations of $\Symn$ is
\cite{MR1824028}; for the Young tableau of this section we have
followed \cite{MR1464693}*{Section 7.2}; see also
\cite{MR1153249}*{Chapter 4}. The irreducibles of $\Symn$ are more
commonly called Specht \emph{modules\/} rather than representations;
as we are not mentioning modules, we hope the change of nomenclature
is not too discombobulating. The representation 
$S^{\BoxedEPSF{tableau1a.eps
      scaled 80}}$
in Figure \ref{fig:symmetric:three} is $\Symthree$ as the symmetries 
--
obtained by permuting its three vertices
--
of
the equilateral triangle. In
general 
$S^{\BoxedEPSF{tableau1.eps scaled 100}}$ is the representation of $\Symn$ acting as the
symmetries of the regular $(n-1)$-simplex; it is another incarnation
of the reflectional representation of $\Symn$ mentioned in Example
\ref{example:permutation:representation} and 
in the notes
to the previous section. The number of irreducible representations of $\Symn$
over $\C$ is equal to the number $p(n)$ of partitions $\lambda\vdash n$;
there is no known closed formula for $p(n)$, but many
weird and wonderful properties are known. To choose just one, there is
the
generating function
$$
\sum_{n=0}^{\infty} p(n)x^n
=
\prod_{k=1}^{\infty}\left(\frac{1}{1-x^k}\right)
$$
Exercise 
\ref{exercise:exterior:powers;reflectional} is
\cite{MR1153249}*{Exercise 4.6}.

\paragraph{Reduction.} We have generally followed
\cite{MR2529864}. The philosophy of the Clifford-Munn correspondence described in
Section \ref{section:clifford:munn}
is that knowledge 
of group representations yields knowledge of semigroup
representations. The passage from groups to
semigroups is the induction construction of 
Section \ref{section:induction}.
The current section is thus a little more perfunctory, as
reduction -- for us -- is merely the inverse construction, its
principal purpose being to establish the Clifford-Munn bijection.
The usual terminology is ``restriction'' in
much of the literature, but we have gone for reduction on two counts:
it is first of all a double restriction -- in that an action of $S$ on
$V$ is being restricted to both a subgroup of
$S$ and a subspace of $V$ -- and secondly, \emph{re}duction seems a
more satisfying counterpoint to 
\emph{in}duction.

\paragraph{Induction.} The section is based on
\cite{MR2460611}*{Chapter 11}. Like there we adopt an elementary
approach; for example, in module-theoretic terms the representation $U$ is
$kS\otimes_{kG_e}V$; the notation $s_i\otimes v$ for the elements of
the copy $V_i$ of $V$ is a nod to this. That the construction is
independent of the transversal $T$ is \cite{MR2460611}*{Theorem 
11.3.1(ii)}. The general picture is from \cite{MR2529864}*{Theorem 7}. The justification
in Example \ref{example:inversesemigroup:annihilator} that $U$ is an
irreducible $I_n$-representation closely follows \cite{MR2460611}*{Theorem 11.3.1}.

\paragraph{The Clifford-Munn correspondence.} Again we have followed
\cite{MR2460611, MR2529864} for the general picture. The irreducibles
of the symmetric inverse monoid in Example \ref{induction:specht} are a
venerable topic. Munn \cite{MR0081910} took a character-theoretic
approach while Grood \cite{MR1887083} constructed the ``Specht'' representations
for $I_n$ from scratch, and seemingly without reference to the
Clifford-Munn correspondence. Our approach  follows
\cite{Majed}, where this and representations of other Boolean reflection monoids are
described. 

\paragraph{A sexy example.} For algebraic groups
and Weyl groups see \cite{Humphreys75} and for algebraic monoids and
Renner monoids, the
books of Putcha and Renner \cite{Putcha88, Renner05}. A beautiful
expository article is \cite{Solomon95}; the example in this section is
taken from
\cite{Solomon95}*{Example 5.7}.

The meaning of ``sensible'', when talking about algebraic
groups and monoids, depends on the context. If $\M$ 
is irreducible, meaning its underlying variety is irreducible, then
the units $\G$ are a connected algebraic group with
$\overline{\G}=\M$. If $\G_0$ is connected semisimple and the
representation $f:\G_0\rightarrow GL(V)$ is rational with finite
kernel, then we have the construction for 
$\ams{M}=\overline{k^\times f(\ams{G}_0)}$ described. 

The Renner monoid of $\m_nk$ is isomorphic to the symmetric inverse
monoid; in this incarnation, 
$I_n$ is called the \emph{Rook monoid\/} and consists of the
$n\times n$ matrices, with $0,1$-entries, such that each row and
column contains \emph{at most\/} one $1$. The name comes about as 
the matrices can be identified with
$n\times
n$ chessboards, with rooks in the positions occupied by $1$'s, 
and with the property that no two rooks are attacking each
other. Warning: the Renner monoid is not in
general a submonoid of $\M$, much as the Weyl group is not in general a subgroup
of $\G$; both $\gl_nk$ and $\m_nk$ are a little special in this
way. 

Good references for polytopes are \cite{Grunbaum03, MR1311028} where
one can also find the combinatorial description of the face polytope
of a permutohedron in terms of ordered partitions.

\vspace{1em}


%
%
%

\end{document}